\documentclass{amsbook}
\usepackage[utf8]{inputenc}
\usepackage{a4}
\usepackage{graphics}
\usepackage{color}
\usepackage{amssymb}
\usepackage[all]{xy}
\usepackage{amsmath}
\usepackage{stmaryrd}

\CompileMatrices
\newtheoremstyle{thm}
  {}
  {}
  {\itshape}
  {0pt}
  {\scshape}
  {.}
  { }
  {}

\newtheorem{theorem}{Theorem}[subsection]
\newtheorem{lemma}[theorem]{Lemma}
\newtheorem{proposition}[theorem]{Proposition}
\newtheorem{corollary}[theorem]{Corollary}
\theoremstyle{definition}

\newtheorem{definition}[theorem]{Definition}
\newtheorem{example}[theorem]{Example}

\newtheorem{construction}[theorem]{Construction}

\newtheorem{remark}[theorem]{Remark}

\theoremstyle{remark}

\numberwithin{equation}{section}
\def\Chi{{\mathbb X}}

\def\div{{\rm div}}

\def\quot{/\!\!/}
\def\mal{\! \cdot \!}

\def\reg{{\rm reg}}
\def\rq#1{\widehat{#1}}
\def\t#1{\widetilde{#1}}
\def\b#1{\overline{#1}}
\def\bangle#1{\langle #1 \rangle}

\def\KK{{\mathbb K}}
\def\TT{{\mathbb T}}
\def\ZZ{{\mathbb Z}}

\def\QQ{{\mathbb Q}}
\def\PP{{\mathbb P}}

\def\LL{{\mathbb L}}
\def\Of{{\mathcal{O}}}

\def\Quot{\operatorname{Quot}}
\def\CDiv{\operatorname{CDiv}}
\def\WDiv{\operatorname{WDiv}}
\def\PDiv{\operatorname{PDiv}}

\def\irr{{\rm irr}}
\def\id{{\rm id}}

\def\ord{{\rm ord}}

\def\Cl{\operatorname{Cl}}

\def\Pic{\operatorname{Pic}}

\def\GL{{\rm GL}}

\def\codim{{\rm codim}}
\def\Supp{{\rm Supp}}
\def\Spec{{\rm Spec}}
\def\Proj{{\rm Proj}}

\def\cone{{\rm cone}}

\def\lin{{\rm lin}}

\def\pr{{\rm pr}}

\newcounter{itemnumber}

\makeindex
\begin{document}
\title{Cox Rings}

%
%
\author[I.~V.~Arzhantsev]{Ivan Arzhantsev} 
\address{Department of Algebra, 
Faculty of Mechanics and Mathematics, Moscow State
  University, Leninskie Gory 1, Moscow, 119991, Russia}
\email{arjantse@mccme.ru}
\author[U.~Derenthal]{Ulrich Derenthal} 
\address{Mathematisches Institut,
Ludwig-Maximilians-Universität Mün\-chen,
Theresienstraße 39,
80333 München,
Germany}
\email{ulrich.derenthal@mathematik.uni-muenchen.de}
\author[J.~Hausen]{J\"urgen Hausen} 
\address{Mathematisches Institut, Universit\"at T\"ubingen,
Auf der Morgenstelle 10, 72076 T\"ubingen, Germany}
\email{juergen.hausen@uni-tuebingen.de}
\author[A.~Laface]{Antonio Laface} 
\address{
Departamento de Matem\'atica, 
Universidad de Concepci\'on, 
Casilla 160-C, 
Concepci\'on, Chile}
\email{alaface@udec.cl}


\maketitle


This is the first chapter of an 
introductory text under construction;
further chapters are available on the 
author's web pages.
Our aim is to provide an elementary access
to Cox rings and their applications in 
algebraic and arithmetic geometry.
We are grateful to Victor Batyrev and 
Alexei Skorobogatov for helpful remarks 
and discussions.
Any comments and suggestions on this draft 
will be highly appreciated.


\vfill

\begin{center}
Version  \today
\end{center}

\tableofcontents

\cleardoublepage


\chapter{Basic concepts}
\label{chap:crbasics}

In this chapter we introduce 
the Cox ring and, more generally, 
the {\em Cox sheaf\/} and its 
geometric counterpart, the 
{\em characteristic space\/}. 
Moreover, algebraic and geometric 
aspects are discussed. 
Section~\ref{sec:gradalg} is devoted
to commutative algebras graded by monoids.
In Section~\ref{sec:gradtoract},
we recall the correspondence between 
actions of quasitori (also called 
diagonalizable groups) 
on affine varieties and 
affine algebras graded by abelian 
groups, 
and we provide the necessary background 
on good quotients.
Section~\ref{sec:divalgs} is a first
step towards Cox rings.
Given a normal variety~$X$
and a finitely generated subgroup 
$K \subseteq \WDiv(X)$ of the group of 
Weil divisors, we consider the associated 
{\em sheaf of divisorial algebras\/}
\begin{eqnarray*}
\mathcal{S}
& = & 
\bigoplus_{D \in K} \mathcal{O}_X(D).
\end{eqnarray*}
We present criteria for local finite generation
and consider the relative spectrum.
A first result 
says that $\Gamma(X,\mathcal{S})$
is a unique factorization domain if
$K$ generates the divisor class group $\Cl(X)$.
Moreover, we characterize divisibility in 
the ring $\Gamma(X,\mathcal{S})$ in terms of 
divisors on $X$.
In Section~\ref{sec:coxrings}, the Cox sheaf 
of a normal variety $X$ with finitely generated 
divisor class group $\Cl(X)$ is introduced; 
roughly speaking it is given as 
\begin{eqnarray*}
\mathcal{R}
& = & 
\bigoplus_{[D] \in \Cl(X)} \mathcal{O}_X(D).
\end{eqnarray*}
The Cox ring then is the corresponding ring of 
global sections.
In the case of a free divisor class group 
well-definiteness is straightforward.
The case of torsion needs some effort,
the precise way to define $\mathcal{R}$ then 
is to take the quotient of an appropriate sheaf 
of divisorial algebras by a certain ideal sheaf.
Basic algebraic properties and divisibility theory
of the Cox ring are investigated in 
Section~\ref{sec:cralgprops}.
Finally, in Section~\ref{sec:univtorsors}, 
we study the characteristic space, i.e., the relative 
spectrum $\rq{X} = \Spec_X \mathcal{R}$ of the Cox 
sheaf. 
It comes with an action of the 
{\em characteristic quasitorus\/} 
$H = \Spec \, \KK[\Cl(X)]$ and a 
good quotient $\rq{X} \to X$.
We relate geometric properties of $X$ to 
properties of this action and give a 
characterization of the characteristic space 
in terms of Geometric Invariant Theory.


\section{Graded algebras}
\label{sec:gradalg}

\subsection{Monoid graded algebras}
We recall basic notions on algebras
graded by abelian monoids.
In this subsection, 
$R$ denotes a commutative ring with unit 
element. 

\begin{definition}
\index{graded algebra!}%
\index{grading!}%
\index{algebra! graded}%
Let $K$ be an abelian monoid.
A {\em $K$-graded $R$-algebra} 
is an associative, commutative 
$R$-algebra $A$ with unit that 
comes with a direct sum decomposition
\begin{eqnarray*}
A 
& = & 
\bigoplus_{w \in K} A_w
\end{eqnarray*}
into $R$-submodules $A_w \subseteq A$ such that 
$A_w \cdot A_{w'} \subseteq A_{w+w'}$ holds for any 
two elements $w,w' \in K$.
\end{definition}

We will also speak of a $K$-graded $R$-algebra
as a monoid graded algebra or just as a graded algebra.
In order to compare $R$-algebras $A$ and $A'$,
which are graded by different abelian monoids 
$K$ and $K'$, we work with the following notion 
of a morphism.

\begin{definition}
\index{graded algebra! morphism of}%
\index{graded homomorphism}%
\index{morphism! of graded algebras}%
A {\em morphism\/} from a $K$-graded algebra $A$ to
a $K'$-graded algebra $A'$ is a pair
$(\psi, \t{\psi})$, where 
$\psi \colon A \to A'$ is a homomorphism
of $R$-algebras, 
$\t{\psi} \colon K \to K'$ is a homomorphism
of abelian monoids and 
\begin{eqnarray*}
\psi(A_w) & \subseteq & A_{\t{\psi}(w)}
\end{eqnarray*}
holds for every $w \in K$.
In the case $K = K'$ and $\t{\psi} = \id_{K}$,
we denote a morphism of graded algebras
just by $\psi \colon A \to A'$ and also 
refer to it as a {\em ($K$-)graded homomorphism}.
\end{definition}

\begin{example}
Given an abelian monoid $K$ and 
$w_1, \ldots, w_r \in K$,
the polynomial ring $R[T_1,\ldots, T_r]$ 
can be made into a $K$-graded $R$-algebra 
by setting 
\begin{eqnarray*}
R[T_1, \ldots, T_r]_w
& := &
\left\{
\sum_{\nu \in \ZZ^r_{\ge 0}} a_\nu T^{\nu}; \; 
a_{\nu} \in R, \, 
\nu_1w_1 + \ldots + \nu_rw_r = w 
\right\}.
\end{eqnarray*}
This $K$-grading is determined by  
$\deg(T_i) = w_i$ for 
$1 \le i \le r$.
Moreover, $R[T_1, \ldots, T_r]$ comes 
with the natural $\ZZ^r_{\ge 0}$-grading 
given by 
\begin{eqnarray*}
R[T_1, \ldots, T_r]_{\nu}
& := &
R \cdot T^{\nu},
\end{eqnarray*}
and we have a canonical morphism 
$(\psi,\t{\psi})$ from $R[T_1, \ldots, T_r]$
to itself,
where $\psi = \id$ and 
$\t{\psi} \colon \ZZ^r_{\ge 0} \to K$
sends $\nu$ to $\nu_1w_1 + \ldots + \nu_rw_r$.
\end{example}

For any abelian monoid $K$, we denote by 
$K^\pm$ the associated group of 
differences and by 
$K_{\QQ} := K^\pm \otimes_{\ZZ} \QQ$
the associated rational vector space.
Note that we have canonical maps 
$K \to K^\pm \to K_{\QQ}$, where the first 
one is injective if $K$ admits cancellation 
and the second one is injective if $K^\pm$ 
is free.
\index{algebra! integral}%
\index{integral algebra!}%
By an integral $R$-algebra, we mean an 
$R$-algebra $A$ without zero-divisors.

\begin{definition}
\index{weight monoid!}%
\index{weight group!}%
\index{weight cone!}%
Let $A$ be an integral $K$-graded $R$-algebra.
The {\em weight monoid} of $A$ 
is the submonoid
$$ 
S(A)
\ := \ 
\{w \in K; \; A_w \ne 0\}
\ \subseteq \ 
K.
$$
The {\em weight group\/} of $A$ is the
subgroup $K(A) \subseteq K^\pm$ generated 
by  $S(A) \subseteq K$.
The {\em weight cone\/} of $A$ is the convex 
cone $\omega(A) \subseteq K_{\QQ}$ 
generated by $S(A) \subseteq K$.
\end{definition}

We recall the construction of 
the algebra associated 
to an abelian monoid; it defines 
a covariant 
functor from the category of 
abelian monoids to the category of 
monoid graded algebras.

\begin{construction}
\label{ex:semigroupalgebra}
\index{monoid algebra!}%
Let $K$ be an abelian monoid. 
As an $R$-module, the associated 
{\em monoid algebra\/} 
over $R$ is given by
\begin{eqnarray*}
R[K]
& := &
\bigoplus_{w \in K} R \mal \chi^w
\end{eqnarray*}
and its multiplication is defined 
by $\chi^w \cdot \chi^{w'} := \chi^{w+w'}$.
If  $K'$ is a further abelian monoid
and $\t{\psi} \colon K \to K'$ 
is a homomorphism,
then we have a homomorphism
$$ 
\psi := R[\t{\psi}] \colon R[K] \ \to \ R[K'],
\qquad
\chi^w \ \mapsto \ \chi^{\t{\psi}(w)}.
$$
The pair $(\psi, \t{\psi})$
is a morphism from the $K$-graded algebra 
$R[K]$ to the $K'$-graded algebra $R[K']$,
and this assignment is functorial.
\end{construction}

Note that the monoid algebra $R[K]$ has $K$ as 
its weight monoid, and $R[K]$ is finitely 
generated over $R$ if and only if the monoid
$K$ is finitely generated.
In general, if a $K$-graded algebra $A$ 
is finitely generated over $R$, 
then its weight monoid is finitely 
generated and its weight cone is 
polyhedral.

\begin{construction}[Trivial extension]
Let $K \subseteq K'$ be an inclusion of 
abelian monoids and~$A$ a $K$-graded 
$R$-algebra.
Then we obtain an $K'$-graded $R$-algebra~$A'$ 
by setting
$$ 
A' 
\ := \ 
\bigoplus_{u \in K'} A'_u,
\qquad\qquad
A'_u 
\ := \ 
\begin{cases}
A_u & \text{if } u \in K,
\\
\{0\} & \text{else}.
\end{cases} 
$$
\end{construction}

\begin{construction}[Lifting]
\label{constr:lifting}
Let $G \colon \t{K} \to K$ be a homomorphism 
of abelian monoids and~$A$ a $K$-graded 
$R$-algebra.
Then we obtain a $\t{K}$-graded $R$-algebra 
$$ 
\t{A}
\ := \ 
\bigoplus_{u \in \t{K}} \tilde A_u,
\qquad\qquad
\t{A}_u 
\ := \ 
A_{G(u)}.
$$
\end{construction}

\begin{definition}
\index{homogeneous ideal!}%
\index{ideal! homogeneous}%
Let $A$ be a $K$-graded $R$-algebra.
An ideal $I \subseteq A$ is called
{\em ($K$-)homogeneous\/} if it
is generated by ($K$-)homogeneous
elements.
\end{definition}

An ideal $I \subseteq A$ of a $K$-graded 
$R$-algebra $A$ is homogeneous if and only 
if it has a direct sum decomposition
$$ 
I 
\ = \ 
\bigoplus_{w \in K} I_w,
\qquad\qquad
I_w \ := \ I \cap A_w.
$$

\begin{construction}
Let $A$ be a $K$-graded $R$-algebra
and $I \subseteq A$ a homogeneous ideal.
Then the factor algebra $A/I$ is 
$K$-graded by
$$ 
A/I
\ = \
\bigoplus_{w \in K} (A/I)_w
\qquad\qquad
 (A/I)_w \ := \  A_w + I.
$$
Moreover, for each homogeneous component 
$(A/I)_w \subseteq A/I$, one has a 
canonical isomorphism of $R$-modules 
$$ 
A_w/I_w \ \to \ (A/I)_w,
\qquad
f + I_w \ \mapsto \ f + I.
$$
\end{construction}

\begin{construction}
\index{coarsened grading!}%
\index{grading! coarsened}%
Let $A$ be a $K$-graded $R$-algebra,
and $\t{\psi} \colon K \to K'$ be a homomorphism
of abelian monoids.
Then one may consider $A$ as a $K'$-graded algebra with
respect to the {\em coarsened grading\/}
$$
A=\bigoplus_{u \in K'} A_u, 
\qquad \qquad 
A_u:=\bigoplus_{\t{\psi}(w)=u} A_w. 
$$
\end{construction}

\begin{example}
Let $K = \ZZ^2$ and consider the 
$K$-grading of $R[T_1, \ldots, T_5]$
given by $\deg(T_i) = w_i$, where
$$ 
w_1 = (-1,2),
\quad
w_2 = (1,0),
\quad
w_3 = (0,1),
\quad
w_4 = (2,-1),
\quad
w_5 = (-2,3).
$$
Then the polynomial 
$T_1T_2 + T_3^2 + T_4T_5$
is $K$-homogeneous of
degree $(0,2)$, and thus we have a 
$K$-graded factor algebra
\begin{eqnarray*}   
A 
& = & 
R[T_1, \ldots, T_5]/\bangle{T_1T_2 + T_3^2 + T_4T_5}.
\end{eqnarray*}
The standard $\ZZ$-grading of the algebra $A$ with 
$\deg(T_1)=\ldots=\deg(T_5)=1$ may be obtained 
by coarsening via the homomorphism 
$\t{\psi} \colon \ZZ^2 \to \ZZ$, $(a,b) \mapsto a+b$.   
\end{example}

\begin{proposition}
\label{prop:integral}
Let $A$ be a $\ZZ^r$-graded $R$-algebra 
satisfying $ff' \ne 0$ for any two 
non-zero homogeneous $f,f' \in A$.
Then the following statements hold.
\begin{enumerate}
\item
The algebra $A$ is integral.
\item
If $gg'$ is homogeneous for 
$0 \ne g,g' \in A$, then 
$g$ and $g'$ are homogeneous.
\item
Every unit $f \in A^*$ is homogeneous.
\end{enumerate}
\end{proposition}

\begin{proof}
Fix a lexicographic ordering on $\ZZ^r$.
Given two non-zero $g,g' \in A$, write 
$g = \sum f_u$ and $g' = \sum f'_u$ 
with homogeneous $f_u$ and $f'_u$.
Then the maximal (minimal) component
of $gg'$ is $f_wf'_{w'} \ne 0$, where
$f_w$ and $f'_{w'}$ are the maximal (minimal)
components of $f$ and $f'$ respectively.
The first two assertions follow.
For the third one observe that $1 \in A$ is 
homogeneous (of degree zero). 
\end{proof}



\subsection{Veronese subalgebras}
We introduce Veronese subalgebras of monoid 
graded algebras and present statements relating 
finite generation of the algebra to finite generation 
of a given Veronese subalgebra and vice versa.
Again, $R$ is a commutative ring with 
unit element.

\begin{definition}
Given an abelian monoid $K$
admitting cancellation,
a $K$-graded $R$-algebra $A$ and 
a submonoid $L \subseteq K$, 
one defines the associated
\index{Veronese subalgebra!}%
{\em Veronese subalgebra\/}
to be
$$
A(L)
\ := \
\bigoplus_{w \in L} A_w
\ \subseteq \ 
\bigoplus_{w \in K} A_w
\ = \ 
A.
$$
\end{definition}

\begin{proposition}
\label{prop:fgsubalgebra}
Let $K$ be an abelian monoid admitting 
cancellation and~$A$ a finitely 
generated $K$-graded $R$-algebra.
Let $L \subseteq K$ be a 
finitely generated submonoid.
Then the associated Veronese subalgebra 
$A(L)$ is finitely generated over $R$.
\end{proposition}

\begin{lemma}
\label{lem:intersectmonoid}
Let $K$ be an abelian monoid 
admitting cancellation 
and let $L,M \subseteq K$ be
finitely generated submonoids.
Then $L \cap M$ is finitely generated.
\end{lemma}

\begin{proof}
Consider the embedding $K \subseteq K^{\pm}$ into 
the group of differences and define a 
homomorphism $\alpha \colon \ZZ^r \to K^{\pm}$
with $L,M \subseteq \alpha(\ZZ^r)$.
Then $\alpha^{-1}(L)$ and 
$\alpha^{-1}(M)$
are finitely generated monoids;
indeed, if  $w_i = \alpha(v_i)$, where $1 \le i \le k$,
generate $L$ and 
$u_1, \dots, u_s$ is a basis for $\ker(\alpha)$,
then $\alpha^{-1}(L)$ is generated by 
$v_1, \dots, v_k$ and  $\pm u_1, \dots, \pm u_s$. 

To prove the assertion, it
suffices to show that the intersection
$\alpha^{-1}(L) \cap \alpha^{-1}(M)$
is finitely generated. 
In other words, we may assume that $K = \ZZ^r$ 
holds.
Then $L$ and $M$ generate convex polyhedral 
cones $\tau$ and $\sigma$ in $\QQ^r$, 
respectively. 
Consider $\omega := \tau \cap \sigma$ 
and the tower of algebras
$$ 
\QQ
\ \subseteq \ 
\QQ[L \cap M] 
\ \subseteq \ 
\QQ[\omega \cap \ZZ^r].
$$
Gordan's Lemma~\cite[Theorem~7.6]{MiSt}
shows that $\QQ[\omega \cap \ZZ^r]$ 
is finitely generated over~$\QQ$.
Moreover, for every $v \in \tau \cap \sigma$,
some positive integral multiple $k \mal v$ 
belongs to $L \cap M$.
Thus, $\QQ[\omega \cap \ZZ^r]$ is integral 
and hence finite over $\QQ[L \cap M]$.
The Artin-Tate Lemma~\cite[page~144]{Ei} 
tells us that 
$\QQ[L \cap M]$ is finitely generated over 
$\QQ$.
Consequently, the weight monoid $L \cap M$ 
of $\QQ[L \cap M]$ is finitely generated.
\end{proof}

\begin{proof}[Proof of Proposition~\ref{prop:fgsubalgebra}]
According to Lemma~\ref{lem:intersectmonoid},
we may assume that $L$ is contained in the 
weight monoid $S(A)$.
Moreover, replacing $K$ with its 
group of differences, we may assume 
that $K$ is a group.
Fix homogeneous generators 
$f_1, \ldots, f_r$ for $A$
and set $w_i := \deg(f_i)$.
Then we have an epimorphism
$$ 
\alpha 
\colon
\ZZ^r \ \to \ K, 
\qquad
e_i \ \mapsto \ w_i.
$$
Moreover, set 
$B := R[T_1, \ldots, T_r]$
and endow it with the natural 
$\ZZ^r$-grading.
Then we obtain a morphism 
$(\pi, \alpha)$ of graded 
$R$-algebras from $B$ to $A$,
where $\pi$ is the epimorphism
defined by
$$ 
\pi \colon 
B \ \to \ A,
\qquad
T_i \ \mapsto \ f_i.
$$
The inverse image $\alpha^{-1}(L) \subseteq \ZZ^r$
is a finitely generated monoid.
By Lemma~\ref{lem:intersectmonoid}, the 
intersection $M := \alpha^{-1}(L) \cap \ZZ_{\ge 0}^r$
is finitely generated and hence generates
a polyhedral convex cone $\sigma = \cone(M)$ 
in $\QQ^r$.
Consider the tower of $R$-algebras
$$ 
R
\ \subseteq \ 
R[M]
\ \subseteq \ 
R[\sigma \cap \ZZ^r].
$$
The $R$-algebra $R[\sigma \cap \ZZ^r]$ 
is finitely generated by
Gordan's Lemma~\cite[Theorem~7.6]{MiSt},
and it is integral and thus finite over $R[M]$.
The Artin-Tate Lemma~\cite[page~144]{Ei} 
then shows that $R[M]$ is finitely generated 
over $R$.
By construction, $\pi \colon B \to A$ maps 
$R[M] \subseteq B$ onto $A(L) \subseteq A$.
This implies finite generation of $A(L)$.
\end{proof}

\begin{proposition}
\label{prop:verofg2allfg}
Suppose that $R$ is noetherian.
Let $K$ be a finitely generated abelian group, 
$A$ a $K$-graded integral $R$-algebra and 
$L \subseteq K$ be a submonoid such that 
for every $w \in S(A)$ there exists an 
$n \in \ZZ_{\ge 1}$ with $nw \in L$.
If the Veronese subalgebra $A(L)$ is 
finitely generated over~$R$, then also 
$A$ is finitely generated over~$R$.
\end{proposition}

\begin{proof}
We may assume that $K$ is generated by $S(A)$.
A first step is to show that the quotient field
of $A$ is a finite extension of that of $A(L)$.
Fix generators  $u_1, \ldots, u_r$ for $K$.
Then we may write $u_i = u_i^+ - u_i^-$ with 
$u_i^{\pm} \in S(A)$.
Choose nontrivial elements 
$g_i^{\pm} \in A_{u_i^{\pm}}$.
With $f_i := g_i^+/g_i^-$ we have
\begin{eqnarray*}
\Quot(A) 
& = & 
\Quot(A(L))(f_1, \ldots, f_r).
\end{eqnarray*}
By our assumption, $A$ is contained in the integral
closure of $A(L)$ in $\Quot(A)$.
Applying~\cite[Proposition~5.17]{AtMa}
we obtain that $A$ is a submodule of some 
finitely generated $A(L)$-module.
Since $R$ and hence $A(L)$ is noetherian, 
$A$ is finitely generated as a module over 
$A(L)$ and thus as an algebra over $R$.
\end{proof}

Putting Propositions~\ref{prop:fgsubalgebra}
and~\ref{prop:verofg2allfg} together,
we obtain the following well known statement
on gradings by abelian groups.

\begin{corollary}
\label{cor:verofg}
Let $R$ be noetherian, $K$ a finitely generated 
abelian group,
$A$ an integral $K$-graded $R$-algebra and 
$L \subseteq K$ a subgroup of finite index.
Then the following statements are equivalent.
\begin{enumerate}
\item
The algebra $A$ is finitely generated over $R$.
\item
The Veronese subalgebra $A(L)$ is finitely  generated over $R$.
\end{enumerate}
\end{corollary}

\begin{proposition}
\label{prop:squeeze}
Suppose that $R$ is noetherian.
Let $L,K$ be abelian monoids admitting cancellation
and $(\varphi,F)$ be a morphism from an $L$-graded 
$R$-algebra $B$ to an integral 
$K$-graded $R$-algebra $A$.
Assume that the weight monoid of $B$ is finitely 
generated
and $\varphi \colon B_u \to A_{F(u)}$ is an 
isomorphism for every $u \in L$.
Then finite generation of $A$ implies 
finite generation of $B$.
\end{proposition}

\begin{proof}
We may assume that $K$ is an abelian group.
In a first step we treat the case $L = \ZZ^r$
without making use of finite 
generation of $S(B)$.
Since $A$ is integral, there are no 
$\ZZ^r$-homogeneous zero divisors in $B$ and
thus $B$ is integral as well, see 
Proposition~\ref{prop:integral}.
Set $L_0 := \ker(F)$. By the elementary divisors
theorem there is a basis $u_1, \ldots, u_r$
for $\ZZ^r$ and $a_1, \ldots, a_s \in \ZZ_{\ge 1}$ 
such that
$a_1u_1, \ldots, a_su_s$ is a basis for $L_0$.
Let $L_1 \subseteq \ZZ^r$ be the sublattice 
spanned by $u_{s+1}, \ldots, u_r$.
This gives Veronese subalgebras
$$ 
B_0 \ := \ \bigoplus_{u \in L_0} B_u,
\qquad
B_1 \ := \ \bigoplus_{u \in L_1} B_u,
\qquad
C \ := \ \bigoplus_{u \in L_0\oplus L_1} B_u.
$$
Note that $\varphi$ maps $B_1$ isomorphically onto 
the Veronese subalgebra of $A$ defined by $F(L_1)$.
In particular, $B_1$ is finitely generated.
Moreover, $C$ is generated by $B_1$ and the 
(unique, invertible) elements 
$f_i^{\pm 1} \in B_{\pm a_iu_i}$ 
mapping to $1 \in A_0$.
Thus, the Veronese subalgebra $C \subseteq B$ 
is finitely generated.
Since $L_0\oplus L_1$ is of finite index in 
$\ZZ^r$, also~$B$ is finitely generated, 
see Corollary~\ref{cor:verofg}.

We turn to the general case.
Let $u_1, \ldots, u_r \in L$ generate the 
weight monoid of~$B$.
Consider the homomorphism 
$G \colon \ZZ^r \to L^{\pm}$
to the group of differences  
sending the $i$-th canonical basis vector 
$e_i \in \ZZ^r$ to $u_i \in L$ and the 
composition $G' := F^{\pm} \circ G$,
where $F^{\pm} \colon L^{\pm} \to K$ 
extends $F \colon L \to K$.
Regarding $B$ as $L^{\pm}$-graded,
$G$ and $G'$ define us lifted algebras 
$\t{B}$ and $\t{A}$, 
see~Construction~\ref{constr:lifting}, 
fitting into a commutative diagram of 
canonical morphisms
$$ 
\xymatrix{
{\t{B}}
\ar[r]
\ar[d]
&
{\t{A}}
\ar[d]
\\
B 
\ar[r]_{\varphi}
&
A
}
$$
The canonical morphism from $\t{A}$ to $A$ 
is as required in the first step
and thus $\t{A}$ is finitely generated.
The weight monoid $M$ of $\t{B}$ is 
generated by the kernel of $G$ and preimages
of generators of the weight monoid of $B$;
in particular, $M$ is finitely generated.
Moreover, $\t{B}$ maps isomorphically 
onto the Veronese subalgebra of $\t{A}$ 
defined by $M \subseteq \ZZ^r$;
here we use that $\varphi \colon B_u \to A_{F(u)}$ 
is an isomorphism for every $u \in L$.
By Proposition~\ref{prop:fgsubalgebra}, 
the algebra $\t{B}$ is finitely generated.
Finally, $\t{B}$ maps onto $B$ which 
gives finite generation of $B$.
\end{proof}


\section{Gradings and quasitorus actions}
\label{sec:gradtoract}

\subsection{Quasitori}
\index{diagonalizable group!}%
\index{group! diagonalizable}%
We recall the functorial correspondence
between finitely generated abelian groups
and quasitori
(also called diagonalizable groups).
Details can be found in the standard textbooks
on algebraic groups, see for 
example~\cite[Section~8]{Bo},
\cite[Section~16]{Hu},
\cite[Section~3.2.3]{OnVi} or~\cite[Section~2.5]{Sp}.

\index{algebraic group!}%
\index{group! algebraic}%
\index{affine algebraic group!}%
\index{group! affine algebraic}%
We work in the category of algebraic varieties
defined over an algebraically closed field
$\KK$ of characteristic zero.
Recall that an {\em (affine) algebraic group\/} is 
an (affine) variety $G$ together with a group 
structure such that the group laws
$$ 
G \times G \to G, \quad (g_1,g_2) \mapsto g_1g_2,
\qquad\qquad
G \to G, \quad g \mapsto g^{-1}
$$
are morphisms of varieties. 
\index{algebraic group! morphism of}%
\index{morphism! of algebraic groups}%
A {\em morphism\/}
of two algebraic groups $G$ and $G'$ 
is a morphism $G \to G'$ of the underlying varieties
which moreover is a homomorphism of groups.

\index{character!}%
A {\em character\/} of an algebraic group 
$G$ is a morphism $\chi \colon G \to \KK^*$ 
of algebraic groups, where $\KK^*$ is the multiplicative
group of the ground field $\KK$.
\index{character group!}%
The {\em character group\/} of $G$ is the 
set $\Chi(G)$ of all characters of 
$G$ together with pointwise multiplication.
Note that $\Chi(G)$ is an abelian group,
and, given any morphism $\varphi \colon G \to G'$ 
of algebraic groups, 
one has a pullback homomorphism
$$ 
\varphi^* \colon \Chi(G') \ \to \ \Chi(G),
\qquad
\chi' \ \mapsto \ \chi' \circ \varphi.
$$

\begin{definition}
\label{def:quasitorusdef}
\index{quasitorus!}%
\index{torus!}%
A {\em quasitorus}
is an affine algebraic group $H$ 
whose algebra of regular functions 
$\Gamma(H,\mathcal{O})$ is generated 
as a $\KK$-vector space by the characters
$\chi \in \Chi(H)$.
A {\em torus\/} is a connected quasitorus.
\end{definition}

\begin{example}
\index{standard $n$-torus!}%
The {\em standard $n$-torus\/} $\TT^n := (\KK^*)^n$ 
is a torus in the sense of~\ref{def:quasitorusdef}.
Its characters are precisely the Laurent monomials
$T^{\nu} = T_1^{\nu_1} \cdots T_n^{\nu_n}$, 
where $\nu \in \ZZ^n$, and its 
algebra of regular functions is the Laurent polynomial 
algebra
$$ 
\Gamma(\TT^n,\mathcal{O})
\ = \ 
\KK[T_1^{\pm 1}, \ldots,T_n^{\pm 1}]
\ = \ 
\bigoplus_{\nu \in \ZZ^n} \KK \cdot T^{\nu}
\ = \ 
\KK[\ZZ^n].
$$
\end{example}

We now associate to any finitely generated 
abelian group $K$ in a functorial way a
quasitorus, namely $H := \Spec\, \KK[K]$;
the construction will show that $H$ is
the direct product of a standard torus and 
a finite abelian group.

\begin{construction}
Let $K$ be any finitely generated abelian group.
Fix generators $w_1, \ldots, w_r$ of $K$ such 
that the epimorphism $\pi \colon \ZZ^r \to K$, 
$e_i \mapsto w_i$ has the kernel 
\begin{eqnarray*}
\ker(\pi)
& = & 
\ZZ a_1e_1 \oplus \ldots \oplus \ZZ a_se_{s}
\end{eqnarray*}
with $a_1, \ldots, a_s \in \ZZ_{\ge 1}$. 
Then we have the following exact sequence of 
abelian groups
$$ 
\xymatrix{
0
\ar[r]
&
{\ZZ^s}
\ar[rr]^{e_i \mapsto a_ie_i}
&&
{\ZZ^r}
\ar[rr]^{e_i \mapsto w_i}
&&
K
\ar[r]
&
0}
$$
Passing to the respective spectra of group
algebras we obtain with $H := \Spec \, \KK[K]$
the following sequence of morphisms  
$$ 
\xymatrix{
1
\ar@{<-}[r]
&
{\TT^s}
\ar@{<-}[rr]^{(t_1^{a_1}, \ldots, t_s^{a_s}) \mapsfrom t}
&&
{\TT^r}
\ar@{<-}[rr]^{\imath}
&&
H
\ar@{<-}[r]
&
1}
$$
The ideal of $H \subseteq \TT^r$ is 
generated by $T_i^{a_i}-1$, where 
$1 \le i \le s$. 
Thus $H$ is a closed subgroup of $\TT^r$
and the sequence is an exact sequence 
of quasitori; note that
$$
H 
\ \cong \
C(a_1)  \times  \ldots \times C(a_s) \times \TT^{r-s},
\qquad\qquad
C(a_i) 
\ := \ 
\{\zeta \in \KK^*; \; \zeta^{a_i} = 1\}.
$$
The group structure on $H = \Spec \, \KK[K]$
does not depend on the choices made:
the multiplication map is given 
by its comorphism
$$ 
\KK[K]
\ \to \
\KK[K] \otimes_{\KK} \KK[K],
\qquad
\chi^w \ \mapsto \chi^w \otimes \chi^w,
$$
and the neutral element of $H = \Spec \, \KK[K]$ 
is the ideal $\bangle{\chi^w-1; \; w \in K}$.
Moreover, every homomorphism 
$\psi \colon K \to K'$ defines a morphism
$$ 
\Spec \, \KK[\psi] \colon \Spec \, \KK[K'] \ \to \  \Spec \, \KK[K].
$$
\end{construction}

\begin{theorem}
We have contravariant exact functors
being essentially inverse to each other:
\begin{eqnarray*}
\{\text{finitely generated abelian groups}\}
& \longleftrightarrow &
\{\text{quasitori}\}
\\[1ex]
K 
& \mapsto & 
\Spec \, \KK[K]
\\
\psi 
& \mapsto & 
\Spec \, \KK[\psi],
\\[1ex]
\Chi(H)
& \mapsfrom &
H,
\\
\varphi^*
& \mapsfrom &
\varphi.
\end{eqnarray*}
Under these equivalences, 
the free finitely generated abelian groups 
correspond to the tori.
\end{theorem}

This statement includes in particular the 
observation that closed subgroups 
as well as homomorphic images 
of quasitori are again quasitori.
Note that homomorphic images of tori are again tori,
but every quasitorus occurs as a closed subgroup
of a torus.

\index{rational representation!}%
\index{representation! rational}%
Recall that a {\em rational representation\/}
of an affine algebraic group $G$ is a 
morphism $\varrho \colon G \to \GL(V)$
to the group $\GL(V)$
of linear automorphisms of 
a finite dimensional $\KK$-vector 
space~$V$.
In terms of representations,
one has the following characterization 
of quasitori, 
see e.g.~\cite[Theorem~2.5.2]{Sp}.

\begin{proposition}
\label{prop:quasitorus2repres}
An affine algebraic group $G$ is 
a quasitorus if and only if 
any rational representation 
of $G$ splits into one-dimensional 
subrepresentations.
\end{proposition}

\subsection{Affine quasitorus actions}
Again we work over an algebraically closed field
$\KK$ of characteristic zero.
Recall that one has contravariant equivalences
between affine algebras, i.e.~finitely generated 
$\KK$-algebras without nilpotent elements,
and affine varieties:
$$ 
A \ \mapsto \ \Spec \, A,
\qquad\qquad
X \ \mapsto \ \Gamma(X, \mathcal{O}).
$$
We first specialize these correspondences to 
graded affine algebras and affine varieties with 
quasitorus action;
here ``graded'' means graded by a finitely 
generated abelian group.
Then we look at basic concepts such as 
orbits and isotropy groups from both sides.

\index{$G$-variety!}%
\index{$G$-variety! morphism of}%
\index{morphism! equivariant}%
\index{equivariant morphism!}%
A {\em $G$-variety\/} is a
variety $X$ together with a morphical action
$G \times X \to X$ of an affine algebraic group $G$.
A {\em morphism} from a $G$-variety $X$ to $G'$-variety
$X'$ is a pair $(\varphi,\t{\varphi})$, where 
$\varphi \colon X \to X'$ is a morphism of varieties 
and $\t{\varphi} \colon G \to G'$ is a morphism 
of algebraic groups such that we have 
$$ 
\varphi(g \mal x) \ = \ \t{\varphi}(g) \mal \varphi(x)
\qquad
\text{for all } 
(g,x) \ \in \ G \times X.
$$
If $G'$ equals $G$ and $\t{\varphi}$ is the identity, 
then we refer to this situation by calling 
$\varphi \colon X \to X'$ a {\em $G$-equivariant\/} morphism.

\begin{example}
\index{diagonal action!}
Let $H$ be a quasitorus.
Any choice of characters
$\chi_1, \ldots, \chi_r \in \Chi(H)$
defines a {\em diagonal $H$-action\/} on 
$\KK^r$ by
\begin{eqnarray*}
h \mal z 
& := & 
(\chi_1(h)z_1, \ldots, \chi_r(h)z_r).
\end{eqnarray*} 
\end{example}

We now associate in functorial manner 
to every affine algebra graded by a finitely 
generated abelian group an affine variety 
with a quasitorus action.

\begin{construction}
Let $K$ be a finitely generated abelian group
and $A$ a $K$-graded affine algebra. 
Set $X = \Spec \, A$.
If $f_i \in A_{w_i}$, $i=1,\ldots,r$, generate $A$, 
then we have a closed embedding
$$ 
X \ \to \ \KK^r,
\qquad
x \ \mapsto \ (f_1(x), \ldots, f_r(x)),
$$
and $X \subseteq \KK^r$ is invariant under 
the diagonal action of $H = \Spec \, \KK[K]$
given by the characters $\chi^{w_1}, \ldots, \chi^{w_r}$.
Note that for any $f \in A$ homogeneity is 
characterized by 
\begin{eqnarray*}
f \in A_w
& \iff & 
f (h \mal x) = \chi^w(h) f(x)
\text{ for all } h \in H, \, x \in X.
\end{eqnarray*}
This shows that the induced $H$-action on $X$ 
does not depend on the embedding into
$\KK^r$: its comorphism is given by 
$$ 
A \ \to \ \KK[K] \otimes_{\KK} A,
\qquad 
A_w \ni f_w \ \mapsto \ \chi^w \otimes f_w \in \KK[K]_w \otimes_{\KK} A_w.
$$
This construction is functorial:
given a morphism $(\psi,\t{\psi})$ from a 
$K$-graded affine algebra $A$ to 
$K'$-graded affine algebra $A'$,
we have a morphism $(\varphi, \t{\varphi})$ from 
the associated $H'$-variety $X'$ to the 
$H$-variety $X$, where 
$\varphi = \Spec \, \psi$ and 
$\t{\varphi} = \Spec \, \KK[\t{\psi}]$.
\end{construction}

\index{representation! rational}%
\index{rational representation!}%
For the other way round, i.e.,  from 
affine varieties $X$ with action of a 
quasitorus~$H$ to graded affine algebras, 
the construction relies on the fact that
the representation of $H$ on 
$\Gamma(X, \mathcal{O})$ is {\em rational\/}, 
i.e., a union of finite dimensional 
rational subrepresentations, 
see~\cite[Proposition~2.3.4]{Sp}
and~\cite[Lemma~2.5]{KKLV} for 
non-affine $X$.
Proposition~\ref{prop:quasitorus2repres}
then shows that it splits into one-dimensional 
subrepresentations.

\begin{construction}
Let a quasitorus $H$ act on a
not necessarily affine variety~$X$.
Then $\Gamma(X, \mathcal{O})$ becomes a rational 
$H$-module by
\begin{eqnarray*} 
(h \mal f) (x) 
& :=  & 
f (h \mal x).
\end{eqnarray*} 
The decomposition of $\Gamma(X, \mathcal{O})$
into one-dimensional subrepresentations 
makes it into a $\Chi(H)$-graded algebra:
$$
\Gamma(X, \mathcal{O})
\ = \ 
\bigoplus_{\chi \in \Chi(H)} \Gamma(X, \mathcal{O})_\chi, 
\qquad 
\Gamma(X, \mathcal{O})_\chi
\ := \ 
\{f \in \Gamma(X, \mathcal{O}); \;  f(h \, \mal \, x) = \chi(h) f(x) \}. 
$$
Again this construction is functorial.
If $(\varphi,\t{\varphi})$ is a morphism from an 
$H$-variety $X$ to an $H'$-variety $X'$,
then $(\varphi^*,\t{\varphi}^*)$ is a morphism of 
the associated graded algebras.
\end{construction}

\begin{theorem}
\label{prop:gradalg2afftorac}
We have contravariant functors
being essentially inverse to each other:
\begin{eqnarray*}
\{\text{graded affine algebras}\}
& \longleftrightarrow &
\{\text{affine varieties with quasitorus action}\}
\\[1em]
A
& \mapsto & 
\Spec \, A,
\\
(\psi,\t{\psi})
& \mapsto & 
(\Spec\, \psi, \Spec\, \KK[\t{\psi}])
\\[1em]
\Gamma(X,\mathcal{O})
& \mapsfrom &
X,
\\
(\varphi^*,\t{\varphi}^*)
& \mapsfrom &
(\varphi,\t{\varphi}).
\end{eqnarray*}
Under these equivalences the graded homomorphisms 
correspond to the equivariant morphisms.
\end{theorem}

We use this equivalence of categories to 
describe some geometry of a quasitorus 
action in algebraic terms. 
The first basic observation is the following.

\begin{proposition}
\label{prop:homogideal}
Let $A$ be a $K$-graded affine algebra
and consider the action of 
$H = \Spec \, \KK[K]$ on $X = \Spec \, A$.
Then for any closed subvariety 
$Y \subseteq X$ and its vanishing ideal
$I \subseteq A$, the following statements
are equivalent.
\begin{enumerate}
\item
The variety $Y$ is $H$-invariant.
\item
The ideal $I$ is homogeneous.
\end{enumerate}
Moreover, if one of these equivalences holds,
then one has a commutative diagram of $K$-graded 
homomorphisms
$$ 
\xymatrix{
{\Gamma(X,\mathcal{O})}
\ar@{<->}[rr]^{\cong}
\ar[d]_{f \mapsto f_{\vert Y}}
&&
A
\ar[d]^{f \mapsto f + I}
\\
{\Gamma(Y,\mathcal{O})} 
\ar@{<->}[rr]_{\cong}
&&
A/I
}
$$
\end{proposition}

We turn to orbits and isotropy groups.
First recall the following fact 
on general algebraic group actions,
see e.g.~\cite[Section~II.8.3]{Hu}.

\begin{proposition}
Let $G$ be an algebraic group,
$X$ a $G$-variety, and let $x \in X$.
Then the isotropy group $G_x \subseteq G$ 
is closed,
the orbit $G \mal x \subseteq X$ is locally closed, 
and one has a commutative diagram 
of equivariant morphisms of $G$-varieties
$$ 
\xymatrix{
& 
G
\ar[dl]_{\pi}
\ar[dr]^{g \mapsto g \cdot x} 
&
\\
G/G_x
\ar[rr]^{\cong}_{gG_x  \mapsto g \cdot x}
&&
{G \cdot x}
}
$$
Moreover, the orbit closure $\b{G \mal x}$ 
is the union of $G \mal x$ and orbits 
of strictly lower dimension and it contains 
a closed orbit.
\end{proposition}

\goodbreak

\begin{definition}
\label{def:orbmonoid}
\index{orbit monoid!}%
\index{orbit group!}%
Let $A$ be a $K$-graded affine algebra
and consider the action of 
$H = \Spec \, \KK[K]$ on $X = \Spec \, A$.
\begin{enumerate}
\item
The {\em orbit monoid\/} of $x \in X$
is the submonoid $S_x \subseteq K$ generated
by all $w \in K$ that admit a function $f \in A_w$ 
with $f(x) \ne 0$.
\item
The {\em orbit group\/} of $x \in X$ is 
the subgroup $K_x \subseteq K$ generated
by the orbit monoid $S_x \subseteq K$.
\end{enumerate}
\end{definition}

\begin{proposition}
\label{prop:isogrorbitlattice}
Let $A$ be a $K$-graded affine algebra,
consider the action of $H = \Spec \, \KK[K]$ 
on $X = \Spec \, A$ and let $x \in X$.
Then there is a commutative diagram
with exact rows
$$ 
\xymatrix{
0
\ar[r]
&
K_x
\ar[r]
\ar[d]^{\cong}
&
K
\ar[r]
\ar[d]_{w \mapsto \chi^w}^{\cong}
&
K/K_x
\ar[r]
\ar[d]^{\cong}
&
0
\\
0
\ar[r]
&
{\Chi(H/H_x)}
\ar[r]_{\pi^*}
&
{\Chi(H)}
\ar[r]_{\imath^*}
&
{\Chi(H_x)}
\ar[r]
&
0
}
$$
where $\imath \colon H_x \to H$
denotes the inclusion of the isotropy group
and $\pi \colon H \to H/H_x$ the projection.
In particular, we obtain 
$H_x \cong \Spec\, \KK[K/K_x]$.
\end{proposition}

\begin{proof}
Replacing $X$ with $\b{H \cdot x}$
does not change~$K_x$.
Moreover, take a homogeneous $f \in A$ 
vanishing along 
$\b{H \cdot x} \setminus H \cdot x$ 
but not at $x$. Then replacing 
$X$ with $X_f$ does not affect $K_x$.
Thus, we may assume that $X = H \mal x$ holds. 
Then the weight monoid of the $H$-variety
$H \mal x$ is $K_x$ and by the commutative diagram
$$ 
\xymatrix{
& 
H
\ar[dl]_{\pi}
\ar[dr]^{h \mapsto h \cdot x} 
&
\\
H/H_x
\ar[rr]_{\cong}
&&
{H \cdot x}
}
$$
we see that $\pi^*(\Chi(H/H_x))$ consists precisely
of the characters $\chi^w$ with $w \in K_x$, which
gives the desired diagram.
\end{proof}

\begin{proposition}
\label{prop:orbmon2orbcl}
Let $A$ be a $K$-graded affine algebra,
consider the action of $H = \Spec \, \KK[K]$ 
on $X = \Spec \, A$ and let $x \in X$.
Then the orbit closure $\b{H \mal x}$
comes with an action of $H/H_x$,
and there is an isomorphism 
$\b{H \mal x} \cong \Spec \, \KK[S_x]$ 
of $H/H_x$-varieties.
\end{proposition}

\begin{proof}
Write for short $Y := \b{H \mal x}$ 
and $V := H \mal x$. 
Then $V \subseteq Y$ is an affine open subset,
isomorphic to $H/H_x$,
and we have a commutative diagram
$$ 
\xymatrix{
{\Gamma(V,\mathcal{O})}
\ar[rr]^{\cong}
&&
{\KK[K_x]}
\\
{\Gamma(Y,\mathcal{O})}
\ar[rr]_{\cong}
\ar[u]^{f \mapsto f_{\vert V}}
&&
{\KK[S_x]}
\ar[u]
}
$$
of graded homomorphisms,
where the horizontal arrows send 
a homogeneous $f$ of degree 
$w$ to $f(x)\chi^{w}$.
The assertion is part of this.
\end{proof}

\begin{proposition}
Let $A$ be an integral $K$-graded affine algebra
and consider the action of $H = \Spec \, \KK[K]$ 
on $X = \Spec \, A$.
Then there is a nonempty invariant open subset
$U \subseteq X$ with 
$$ 
S_x \ = \ S(A),
\qquad
K_x \ = \ K(A)
\qquad
\text{for all }
x \in U.
$$
\end{proposition}

\begin{proof}
Let $f_1, \ldots, f_r$ be homogeneous generators for $A$.
Then the set $U \subseteq X$ obtained by removing the 
zero sets $V(X,f_i)$ from $X$ for $i=1,\ldots, r$
is as wanted.
\end{proof}

Recall that an action of a group $G$ on a set $X$ 
is said to be
\index{effective! action}%
\index{action! effective}%
{\em effective\/} if $g \cdot x = x$ for all $x \in X$ 
implies $g = e_G$.

\begin{corollary}
Let $A$ be an integral $K$-graded affine algebra
and consider the action of $H = \Spec \, \KK[K]$ 
on $X = \Spec \, A$.
Then the action of $H$ on $X$ is effective 
if and only if $K = K(A)$ holds.
\end{corollary}



\subsection{Good quotients}
We summarize the basic facts
on good quotients.
Everything takes place over an algebraically 
closed field $\KK$ of characteristic zero.
Besides varieties, we consider more generally 
possibly non-separated prevarieties.%
\index{prevariety!}
By definition, a {\em \mbox{($\KK$-)prevariety}\/} 
is a space $X$ with a sheaf $\mathcal{O}_X$ of 
$\KK$-valued functions covered 
by open subspaces $X_1, \ldots, X_r$,
each of which is an affine \mbox{($\KK$-)variety}.

\index{morphism! invariant}\index{invariant morphism!}%
\index{morphism! affine}\index{affine morphism!}%
Let~an algebraic group $G$ act 
on a prevariety $X$, where, here and later,
we always assume that this action is given by 
a morphism $G \times X \to X$.  
Recall that a morphism $\varphi \colon X \to Y$ 
is said to be
\index{morphism! invariant}\index{invariant morphism!}%
{\em $G$-invariant\/} 
if it is constant along the orbits. 
Moreover, a morphism  $\varphi \colon X \to Y$ is called%
\index{morphism! affine}\index{affine morphism!}
{\em affine\/} if for any open affine 
$V \subseteq Y$ the preimage
$\varphi^{-1}(V)$ is an affine variety.
\index{reductive group!}%
\index{group! reductive}%
When we speak of a {\em reductive\/} algebraic 
group,
we mean a not necessarily connected affine 
algebraic group $G$ such that every rational 
representation of $G$ splits into irreducible 
ones.

\begin{definition}
\index{good quotient!}%
\index{quotient! good}%
\index{geometric quotient!}%
\index{quotient! geometric}%
\label{def:goodquot}
Let $G$ be a reductive algebraic group $G$ 
act on a prevariety $X$.
A morphism $p \colon X \to Y$ of 
prevarieties is called 
a {\em good quotient\/} for this action  
if it has the following properties:
\begin{enumerate}
\item
$p \colon X \to Y$ is affine and $G$-invariant,
\item
the pullback $p^* \colon \mathcal{O}_Y \to (p_*\mathcal{O}_X)^G$
is an isomorphism.
\end{enumerate}
A morphism $p \colon X \to Y$ is called a 
{\em geometric quotient\/} if it is  a good quotient
and its fibers are precisely the $G$-orbits.
\end{definition}

\begin{remark}
Let $X = \Spec \, A$ be an affine $G$-variety
with a reductive algebraic group $G$. 
The finiteness theorem of Classical
Invariant Theory ensures that the algebra 
of invariants $A^G \subseteq A$ is finitely 
generated~\cite[Section~II.3.2]{Kr}.
This guarantees existence of a good quotient 
$p \colon X \to Y$, where $Y:= \Spec \, A^G$.
The notion of a good quotient is locally modeled 
on this concept, because 
for any good quotient $p' \colon X' \to Y'$ 
and any affine open $V \subseteq Y'$
the variety $V$ is isomorphic to $\Spec\,\Gamma(p'^{-1}(V), \Of)^G$, and the 
restricted morphism $p'^{-1}(V) \to V$ is the morphism just described. 
\end{remark}

\goodbreak

\begin{example}
\label{exam:1torusonk2}
Consider the $\KK^*$-action $t \cdot (z,w) = (t^az,t^bw)$ 
on $\KK^2$.
The following three cases are typical.
\begin{enumerate}
\item
We have $a=b=1$. Every $\KK^*$-invariant function is 
constant and the constant map $p \colon \KK^2 \to \{{\rm pt}\}$
is a good quotient.
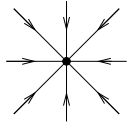
\begin{center}
\begin{picture}(0,0)%
\includegraphics{elliptic.pstex}%
\end{picture}%
\setlength{\unitlength}{829sp}%
\begingroup\makeatletter\ifx\SetFigFontNFSS\undefined%
\gdef\SetFigFontNFSS#1#2#3#4#5{%
  \reset@font\fontsize{#1}{#2pt}%
  \fontfamily{#3}\fontseries{#4}\fontshape{#5}%
  \selectfont}%
\fi\endgroup%
\begin{picture}(3666,6454)(1318,-6482)
\put(4051,-5011){\makebox(0,0)[lb]{\smash{{\SetFigFontNFSS{6}{7.2}{\familydefault}{\mddefault}{\updefault}{\color[rgb]{0,0,0}$p$}%
}}}}
\end{picture}%

\end{center}
\item
We have $a=0$ and $b=1$. 
The algebra of $\KK^*$-invariant functions
is generated by $z$ and  the map $p \colon \KK^2 \to \KK$,
$(z,w) \mapsto z$ is a good quotient.
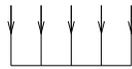
\begin{center}
\begin{picture}(0,0)%
\includegraphics{parabolic.pstex}%
\end{picture}%
\setlength{\unitlength}{829sp}%
\begingroup\makeatletter\ifx\SetFigFontNFSS\undefined%
\gdef\SetFigFontNFSS#1#2#3#4#5{%
  \reset@font\fontsize{#1}{#2pt}%
  \fontfamily{#3}\fontseries{#4}\fontshape{#5}%
  \selectfont}%
\fi\endgroup%
\begin{picture}(3846,4566)(1228,-4594)
\put(4051,-3211){\makebox(0,0)[lb]{\smash{{\SetFigFontNFSS{6}{7.2}{\familydefault}{\mddefault}{\updefault}{\color[rgb]{0,0,0}$p$}%
}}}}
\end{picture}%

\end{center}
\item
We have $a=1$ and $b=-1$. 
The algebra of $\KK^*$-invariant functions
is generated by $zw$ and  $p \colon \KK^2 \to \KK$,
$(z,w) \mapsto zw$ is a good quotient.
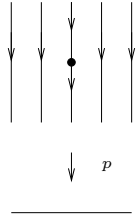
\begin{center}
\begin{picture}(0,0)%
\includegraphics{hyperbolic.pstex}%
\end{picture}%
\setlength{\unitlength}{829sp}%
\begingroup\makeatletter\ifx\SetFigFontNFSS\undefined%
\gdef\SetFigFontNFSS#1#2#3#4#5{%
  \reset@font\fontsize{#1}{#2pt}%
  \fontfamily{#3}\fontseries{#4}\fontshape{#5}%
  \selectfont}%
\fi\endgroup%
\begin{picture}(3846,6366)(1228,-6394)
\put(4051,-5011){\makebox(0,0)[lb]{\smash{{\SetFigFontNFSS{6}{7.2}{\familydefault}{\mddefault}{\updefault}{\color[rgb]{0,0,0}$p$}%
}}}}
\end{picture}%

\end{center}
Note that the general $p$-fiber is a single 
$\KK^*$-orbit, whereas $p^{-1}(0)$ consists
of three orbits and is reducible.
\end{enumerate}
\end{example}

\begin{example}
Let $A$ be a $K$-graded affine algebra.
Consider a homomorphism $\psi \colon K \to L$  
of abelian groups and the coarsified 
grading
$$ 
A \ = \ \bigoplus_{u \in L} A_u,
\qquad
A_u \ = \ \bigoplus_{w \in \psi^{-1}(u)} A_w.
$$
Then the diagonalizable group $H = \Spec\, \KK[L]$
acts on $X = \Spec\, A$, and for the algebra 
of invariants we have 
$$ 
A^H
\ = \ 
\bigoplus_{w \in \ker(\psi)} A_w.
$$
Note that in this special case, 
Proposition~\ref{prop:fgsubalgebra}
ensures finite generation of the algebra 
of invariants.
\end{example}

\begin{example}[Veronese subalgebras]
\index{Veronese subalgebra!}
Let $A$ be a $K$-graded 
affine algebra and $L \subseteq K$
a subgroup. 
Then we have the corresponding 
Veronese subalgebra
$$ 
B
\ = \ 
\bigoplus_{w \in L} A_{w}
\ \subseteq \ 
\bigoplus_{w \in K} A_{w}
\ = \ 
A.
$$
By the preceding example, the morphism 
$\Spec\, A \to \Spec\, B$ is a good 
quotient for the action of $\Spec \, \KK[K/L]$
on $\Spec\, A$.
\end{example}

We list basic properties
of good quotients. 
The key to most of the statements 
is the following central observation.

\begin{theorem}
\label{prop:GclosGsep}
Let a reductive algebraic group $G$ 
act on a prevariety $X$.
Then any good quotient  
$p \colon X \to Y$ has the following
properties.
\begin{enumerate}
\item
\index{$G$-closedness!}
$G$-closedness:
If $Z \subseteq X$ is $G$-invariant and closed,
then its image $p(Z) \subseteq Y$ is closed.
\item
\index{$G$-separation!}
$G$-separation:
If $Z,Z' \subseteq X$ are $G$-invariant, closed
and disjoint, then $p(Z)$ and $p(Z')$ are 
disjoint.
\end{enumerate}
\end{theorem}

\begin{proof}
Since  $p \colon X \to Y$ is affine
and the statements are local with respect to~$Y$, 
it suffices to prove them for affine $X$.
This is done in~\cite[Section~II.3.2]{Kr},
or~\cite[Theorems~4.6 and~4.7]{PoVi}.
\end{proof}

As an immediate consequence, one obtains 
basic information on the structure of 
the fibers of a good quotient.

\begin{corollary}
\label{prop:goodquotfibers}
Let a reductive algebraic group $G$ 
act on a prevariety $X$, and let 
$p \colon X \to Y$ be a good quotient.
Then $p$ is surjective and 
for any $y \in Y$ one has:
\begin{enumerate}
\item
There is exactly one closed $G$-orbit 
$G \mal x$ in the fiber $p^{-1}(y)$.
\item
Every orbit $G \mal x' \subseteq p^{-1}(y)$
has $G \mal x$ in its closure.
\end{enumerate}
\end{corollary}

The first statement means that a good quotient 
$p \colon X \to Y$ parametrizes the closed orbits 
of the $G$-prevariety $X$. 

\goodbreak 

\begin{corollary}
\label{cor:univgoodquot}
Let a reductive algebraic group $G$ 
act on a prevariety $X$, and let 
$p \colon X \to Y$ be a good quotient.
\begin{enumerate}
\item
The quotient space $Y$ carries the quotient topology 
with respect to the map $p \colon X \to Y$.
\item
For every $G$-invariant 
morphism of prevarieties
$\varphi \colon X \to Z$, there 
is a unique morphism $\psi \colon Y \to Z$ 
with $\varphi = \psi \circ p$.
\end{enumerate}
\end{corollary}

\begin{proof}
The first assertion follows from 
Theorem~\ref{prop:GclosGsep}~(i).
The second one follows
from Corollary~\ref{prop:goodquotfibers},
Property~\ref{def:goodquot}~(ii)
and the first assertion.
\end{proof}

\index{categorical quotient!}%
\index{quotient! categorical}%
A morphism $p \colon X \to Y$ with the last property 
is also called a {\em categorical quotient\/}.
The fact that a good quotient is categorical implies 
in particular, that the good quotient space is unique 
up to isomorphy. 
This justifies the notation $X \to X \quot G$ for 
good and $X \to X/G$ for geometric quotients.

\begin{proposition}
\label{prop:goodquotperm}
Let a reductive algebraic group $G$ 
act on a prevariety $X$, and let 
$p \colon X \to Y$ be a good quotient.
\begin{enumerate}
\item 
Let $V \subseteq Y$ be an open subset.
Then the restriction  $p \colon p^{-1}(V) \to V$
is a good quotient for the restricted $G$-action. 
\item
Let $Z \subseteq X$ be a closed $G$-invariant subset.
Then the restriction  $p \colon Z \to p(Z)$
is a good quotient for the restricted $G$-action. 
\end{enumerate}
\end{proposition}

\begin{proof}
The first statement is clear and the second one follows 
immediately from the corresponding statement on the 
affine case, see~\cite[Section~II.3.2]{Kr}.
\end{proof}


\begin{example}[The Proj construction]
\label{ex:projconstr}
\index{ideal! irrelevant}%
\index{irrelevant ideal!}%
\index{Proj construction!}%
Let $A = \oplus A_d$ be a  
$\ZZ_{\ge 0}$-graded affine algebra.
The {\em irrelevant ideal\/} in $A$ is
defined as
$$ 
A_{> 0}
\ := \ 
\bangle{f; \; f \in A_d \text{ for some } d > 0}
\ \subseteq \ 
A.
$$ 
For any homogeneous $f \in A_{>0}$ 
the localization 
$A_f$ is a $\ZZ$-graded affine algebra; 
concretely, the grading is given by
$$ 
A_f \ = \ \bigoplus_{d \in \ZZ} (A_f)_d,
\qquad
(A_f)_d
\ := \ 
\{h/f^l \in A_f; \; \deg(h) - l \deg(f) = d\}.
$$
In particular, we have the, again finitely 
generated, degree zero part of $A_f$; 
it is given by
$$ 
A_{(f)} 
\ := \
(A_f)_0
\ = \ 
\{h/f^l \in A_f; \; \deg(h) = l \deg(f)\}.
$$
Set $X := \Spec(A)$ and $Y_0 := \Spec(A_0)$,
and, for a homogeneous $f \in A_{>0}$, 
set 
$X_f := \Spec \, A_{f}$
and
$U_f := \Spec \, A_{(f)}$.
Then, for any two homogeneous $f,g \in A_{>0}$,
we have the commutative diagrams
$$
\xymatrix{
A_{f} 
\ar[r]
&
A_{fg}
&
A_{g}
\ar[l]
\\
A_{(f)}
\ar[r]
\ar[u] 
&
A_{(fg)}
\ar[u] 
&
A_{(g)}
\ar[l]
\ar[u] 
\\
& 
A_0
\ar[u] 
\ar[ul] 
\ar[ur] 
&
}
\qquad
\qquad
\qquad
\xymatrix{
X_{f} 
\ar[d]_{\pi_f}
&
X_{fg}
\ar[d]
\ar[l]
\ar[r]
&
X_{g}
\ar[d]^{\pi_g}
\\
U_{f} 
\ar[dr] 
&
U_{fg}
\ar[l]
\ar[r]
\ar[d]
&
U_{g}
\ar[dl] 
\\
& 
Y_0
&
}
$$
where the second one arises from the first 
one by applying the $\Spec$-functor.
The morphisms $U_{fg} \to U_f$
are open embeddings and 
gluing the $U_f$ gives the variety 
$Y = \Proj(A)$.
With the zero set $F := V(X,A_{>0})$
of the ideal $A_{>0}$,
we have canonical morphisms,
where the second one is projective:
$$ 
\xymatrix{
{X \setminus F}
\ar[r]^{\pi}
&
Y
\ar[r]
&
Y_0}.
$$
Geometrically the following happened.
The subset $F \subseteq X$ 
is precisely the fixed point set of 
the $\KK^*$-action on $X$ given by the 
grading.
Thus, $\KK^*$ acts with closed orbits 
on $W := X \setminus F$.
The maps $X_f \to U_f$ are geometric 
quotients, and glue together to 
a geometric quotient 
$\pi \colon W \to Y$.
Moreover, the $\KK^*$-equivariant
inclusion $W \subseteq X$
induces the morphism of 
quotients $Y \to Y_0$.
\end{example}


\section{Divisorial algebras}
\label{sec:divalgs}

\subsection{Sheaves of divisorial algebras}
We work over an algebraically closed field~$\KK$ of
characteristic zero. 
We will not only deal with varieties over $\KK$ 
but more generally with prevarieties.

Let $X$ be an irreducible prevariety.
\index{Weil divisor!}%
\index{divisor! Weil}%
The group of {\em Weil divisors\/} of $X$ 
is the free abelian 
group $\WDiv(X)$ generated by all%
\index{prime divisor!}%
\index{divisor! prime}
{\em prime divisors},
i.e., irreducible subvarieties $D \subseteq X$ 
of codimension one.%
\index{divisor! of a function}
To a non-zero rational function $f \in \KK(X)^*$ 
one associates a Weil divisor 
using its \index{order! of a function}%
{\em order\/} along  prime 
divisors~$D$; recall that,
if $f$ belongs to the local ring $\mathcal{O}_{X,D}$, 
then $\ord_D(f)$ is 
the length of the $\mathcal{O}_{X,D}$-module 
$\mathcal{O}_{X,D}/\bangle{f}$,
and otherwise one 
writes $f = g/h$ with $g,h \in \mathcal{O}_{X,D}$
and defines the order of $f$ to be the 
difference of the orders of $g$ and $h$.
The divisor of $f \in \KK(X)^*$ then is
\begin{eqnarray*}
\div(f)
& := &
\sum_{D \ \rm{prime}} \ord_D(f) \cdot D.
\end{eqnarray*}
The assignment $f \mapsto \div(f)$ is a homomorphism
$\KK(X)^* \to \WDiv(X)$, and its image 
$\PDiv(X) \subseteq \WDiv(X)$ is called the 
subgroup of {\em principal divisors}.%
\index{principal divisor!}%
\index{divisor! principal}%
\index{divisor class group!}
The {\em divisor class group\/} 
of~$X$ is the factor group
\begin{eqnarray*}
\Cl(X)
& := & 
\WDiv(X) \ / \ \PDiv(X).
\end{eqnarray*}

A Weil divisor 
$D = a_1 D_1 + \ldots + a_s D_s$ with prime 
divisors $D_i$ is called
\index{effective divisor!}%
\index{divisor! effective}%
{\em effective},
denoted as $D \ge 0$, if 
$a_i \ge 0$ holds for $i = 1, \ldots, s$.
To every divisor $D \in \WDiv(X)$, one associates 
a sheaf $\mathcal{O}_X(D)$ of $\mathcal{O}_X$-modules
by defining its sections over an open $U \subseteq X$ 
as
\begin{eqnarray*}
\Gamma(U,\mathcal{O}_X(D))
& := & 
\{f \in \KK(X)^*; \; (\div(f) + D)_{\vert U} \ge 0 \}
\ \cup \
 \{0\},
\end{eqnarray*}
where the restriction map $\WDiv(X) \to \WDiv(U)$
is defined for a prime divisor $D$ as
$D_{\vert U} := D \cap U$ if it intersects $U$ 
and $D_{\vert U} := 0$ otherwise.
Note that for any two functions 
$f_1 \in \Gamma(U,\mathcal{O}_X(D_1))$ 
and 
$f_2 \in \Gamma(U,\mathcal{O}_X(D_2))$ 
the product $f_1f_2$ belongs to
$\Gamma(U,\mathcal{O}_X(D_1+D_2))$.

\begin{definition}
\index{divisorial algebra! sheaf of}
\index{sheaf of divisorial algebras}
The {\em sheaf of divisorial algebras\/} 
associated to a subgroup
$K \subseteq \WDiv(X)$
is the sheaf of $K$-graded 
$\mathcal{O}_X$-algebras
$$ 
\mathcal{S} 
\ := \ 
\bigoplus_{D \in K} \mathcal{S}_D,
\qquad\qquad
\mathcal{S}_D
\ := \ 
\mathcal{O}_X(D),
$$
where the multiplication in $\mathcal{S}$ 
is defined by multiplying homogeneous 
sections in the field of functions $\KK(X)$. 
\end{definition}

\begin{example}
On the projective line $X = \PP_1$,
consider $D := \{\infty\}$,
the group 
$K := \ZZ D$, 
and 
the associated $K$-graded sheaf 
of algebras $\mathcal{S}$.
Then we have isomorphisms
$$ 
\varphi_n 
\colon \KK[T_0,T_1]_n \ \to \ \Gamma(\PP_1, \mathcal{S}_{nD}),
\qquad
f \ \mapsto \ f(1,z),
$$
where $\KK[T_0,T_1]_n \subseteq \KK[T_0,T_1]$ 
denotes the vector space of all polynomials 
homogeneous of degree~$n$.
Putting them together we obtain a graded 
isomorphism
$$ 
\KK[T_0,T_1] 
\ \cong \
\Gamma(\PP_1, \mathcal{S}).
$$
\end{example}

Fix a normal (irreducible) prevariety~$X$,
a subgroup $K \subseteq \WDiv(X)$ 
on the normal prevariety $X$
and let $\mathcal{S}$ be the 
associated divisorial algebra.
We collect first properties.

\begin{remark}
\label{rem:codim2divsheaf}
If  $V \subseteq U \subseteq X$ are open 
subsets such that $U \setminus V$ is of 
codimension at least two in $U$, then 
we have an isomorphism 
$$
\Gamma(U, \mathcal{S}) 
\ \to \
\Gamma(V, \mathcal{S}).
$$
In particular, the algebra $\Gamma(U, \mathcal{S})$
equals the algebra $\Gamma(U_\reg, \mathcal{S})$,
where $U_\reg \subseteq U$ denotes the set
of smooth points.
\end{remark}

\begin{remark}
\label{rem:locprinc2split}
Assume that $D_1, \ldots, D_s$ is 
a basis for $K \subseteq \WDiv(X)$
and suppose that $U \subseteq X$ 
is an open subset on which each $D_i$ 
is principal, say $D_i = \div(f_i)$.
Then, with $\deg(T_i) = D_i$ and 
$f_i^{-1} \in \Gamma(X,\mathcal{S}_{D_i})$, 
we have a graded isomorphism
$$ 
\Gamma(U, \mathcal{O})
\otimes_{\KK} 
\KK[T_1^{\pm 1}, \ldots, T_s^{\pm 1}]
\ \to \ 
\Gamma(U,\mathcal{S}),
\qquad
g \otimes T_1^{\nu_1}\cdots T_s^{\nu_s}
\ \mapsto \ 
g f_1^{-\nu_1}\cdots f_s^{-\nu_s}.
$$
\end{remark}

\begin{remark}
\label{RR}
If $K$ is of finite rank, say $s$,
then the algebra $\Gamma(X, \mathcal{S})$ 
of global sections
can be realized as a graded subalgebra 
of the Laurent polynomial algebra
$\KK(X)[T_1^{\pm 1}, \ldots, T_s^{\pm 1}]$.
Indeed, let $D_1, \ldots, D_s$ be a basis for 
$K$.
Then we obtain a monomorphism
$$ 
\Gamma(X, \mathcal{S}) 
\ \to \ 
\KK(X)[T_1^{\pm 1}, \ldots, T_s^{\pm 1}],
\qquad
\Gamma(X, \mathcal{S}_{a_1D_1 + \ldots + a_sD_s})
\ni f
\ \mapsto \ 
f T_1^{a_1} \cdots T_s^{a_s}.
$$
In particular, $\Gamma(X, \mathcal{S})$ 
is an integral ring and we have an embedding
of the associated quotient fields
\begin{eqnarray*}
\Quot(\Gamma(X, \mathcal{S})) 
& \to & 
\KK(X)(T_1, \ldots, T_s).
\end{eqnarray*}
For quasiaffine $X$, we have 
$\KK(X) \subseteq \Quot(\Gamma(X, \mathcal{S}))$
and for each variable $T_i$ 
there is a non-zero function 
$f_i \in \Gamma(X, \mathcal{S}_{D_i})$.
Thus, for quasiaffine $X$, one obtains
\begin{eqnarray*}
\Quot(\Gamma(X, \mathcal{S})) 
& \cong & 
\KK(X)(T_1, \ldots, T_s).
\end{eqnarray*}
\end{remark}

\index{divisor! support of}%
\index{support! of a divisor}%
The {\em support\/} $\Supp(D)$ of a Weil divisor 
$D = a_1 D_1 + \ldots + a_s D_s$ with prime 
divisors $D_i$ is the union of those $D_i$
with $a_i \ne 0$.
Moreover, for a Weil divisor $D$ on a normal 
prevariety $X$ and a non-zero section 
$f \in \Gamma(X, \mathcal{O}_X(D))$,
\index{$D$-divisor!}%
\index{$D$-localization!}%
we define the {\em $D$-divisor\/} 
and the {\em $D$-localization\/}
$$ 
\div_D(f)
\ := \ 
\div(f) + D
\ \in \ 
\WDiv(X),
\qquad
X_{D,f}
\ := \ 
X \setminus \Supp(\div_D(f))
\ \subseteq \ 
X.
$$
The $D$-divisor is always effective.
Moreover, given sections 
$f \in \Gamma(X, \mathcal{O}_X(D))$
and
$g \in \Gamma(X, \mathcal{O}_X(E))$, 
we have
$$ 
\div_{D+E}(fg)
\ = \
\div_{D}(f)
+
\div_{E}(g),
\qquad\qquad
f^{-1} 
\ \in \ 
\Gamma(X_{D,f}, \mathcal{O}_X(-D)).
$$

\begin{remark}
\label{rem:Sloc}
Let $D \in K$ and consider a non-zero 
homogeneous section
$f \in \Gamma(X,\mathcal{S}_D)$.
Then one has a canonical isomorphism of $K$-graded algebras
\begin{eqnarray*}
\Gamma(X_{D,f}, \mathcal{S}) 
& \cong & 
\Gamma(X,\mathcal{S})_f.
\end{eqnarray*}
Indeed, the canonical monomorphism 
$\Gamma(X,\mathcal{S})_f \to \Gamma(X_{D,f}, \mathcal{S})$
is surjective, because for any 
$g \in \Gamma(X_{D,f},\mathcal{S}_E)$,
we have $gf^m \in \Gamma(X,\mathcal{S}_{mD+E})$
with some $m \in \ZZ_{\ge 0}$.
\end{remark}

\subsection{The relative spectrum}
Again we work over an algebraically 
closed field $\KK$ of characteristic zero.
Let $X$ be a normal prevariety.
As any quasicoherent sheaf of 
$\mathcal{O}_X$-algebras,
the sheaf of divisorial algebras $\mathcal{S}$ 
associated to 
a group $K \subseteq \WDiv(X)$ 
of Weil divisors
defines in a natural way a geometric 
object, its relative spectrum
$\t{X} := \Spec_X \mathcal{S}$.
We briefly recall how to obtain it.

\begin{construction}
\label{constr:relspec}
\index{relative spectrum!}
\index{spectrum! relative}
Let $\mathcal{S}$ 
be any quasicoherent sheaf of 
reduced $\mathcal{O}_X$-algebras 
on a prevariety $X$, and
suppose that  $\mathcal{S}$ is locally 
of finite type, i.e., $X$ is covered 
by open affine subsets 
$X_1, \ldots, X_r \subseteq X$
with $\Gamma(X_i,\mathcal{S})$ finitely 
generated.
Cover each intersection $X_{ij} := X_i \cap X_j$ by 
open subsets $(X_{i})_{f_{ijk}}$, where 
$f_{ijk} \in \Gamma(X_i,\mathcal{O})$.
Set $\t{X}_i := \Spec \, \Gamma(X_i,\mathcal{S})$
and let $\t{X}_{ij} \subseteq \t{X}_i$ 
be the union of the open subsets
$(\t{X}_i)_{f_{ijk}}$.
Then we obtain commutative diagrams
$$
\xymatrix{
{\t{X}_i}
\ar[d]
& 
{\t{X}_{ij}}
\ar[l]
\ar@{<->}[r]^{\cong}
\ar[d]
& 
{\t{X}_{ji}}
\ar[r]
\ar[d]
&
{\t{X}_j}
\ar[d]
\\
X_i
& 
X_{ij}
\ar[l]
\ar@{=}[r]
& 
X_{ji}
\ar[r]
&
X_j
}
$$
This allows us to glue together the $\t{X}_i$ 
along the $\t{X}_{ij}$, and we obtain a 
prevariety $\t{X} = \Spec_X \, \mathcal{S}$ 
coming with a canonical morphism
$p \colon \t{X} \to X$.
Note that $p_*(\mathcal{O}_{\t{X}}) = \mathcal{S}$
holds.
In particular, $\Gamma(\t{X},\mathcal{O})$ equals
$\Gamma(X,\mathcal{S})$.
Moreover, the morphism $p$ is affine 
and $\t{X}$ is separated if $X$ is 
so. 
Finally, the whole construction does not depend
on the choice of the $X_i$.
\end{construction}

Before specializing this construction to the 
case of our sheaf of divisorial algebras $\mathcal{S}$
on~$X$, we provide two criteria for $\mathcal{S}$
being locally of finite type.
The first one is an immediate consequence of 
Remark~\ref{rem:Sloc}.

\begin{proposition}
\label{prop:globfin2locfin}
Let $X$ be a normal prevariety, 
$K \subseteq \WDiv(X)$ a finitely generated 
subgroup, and $\mathcal{S}$ the associated 
sheaf of divisorial algebras.
If $\Gamma(X,\mathcal{S})$ is finitely 
generated and $X$ is covered by affine open subsets 
of the form $X_{D,f}$, where $D \in K$ and 
$f \in \Gamma(X,\mathcal{S}_D)$, 
then $\mathcal{S}$
is locally of finite type.
\end{proposition}

\index{Cartier divisor!}%
\index{divisor! Cartier}%
A Weil divisor $D \in \WDiv(X)$ on a prevariety 
$X$ is called {\em Cartier\/} if it is locally a principal 
divisor, i.e., locally of the form $D = \div(f)$ 
with a rational function $f$.
\index{locally factorial! variety}%
The prevariety $X$ is {\em locally factorial\/}, i.e., 
all local rings $\mathcal{O}_{X,x}$ are unique 
factorization domains if and only if every Weil
divisor of $X$ is Cartier.
Recall that smooth prevarieties are locally factorial.
\index{$\QQ$-factorial! prevariety}%
More generally, a normal prevariety is called 
{\em $\QQ$-factorial\/} if for any 
Weil divisor some positive multiple is Cartier.

\begin{proposition}
\label{prop:qfact2locfin}
Let $X$ be a normal prevariety and 
$K \subseteq \WDiv(X)$ a finitely 
generated subgroup. 
If $X$ is $\QQ$-factorial, then the 
associated sheaf $\mathcal{S}$ of 
divisorial algebras is locally of finite type.
\end{proposition}

\begin{proof}
By $\QQ$-factoriality, the subgroup $K^0 \subseteq K$ 
consisting of all Cartier divisors is of finite index
in $K$.
Choose a basis $D_1, \ldots, D_s$ for $K$ 
such that with suitable $a_i > 0$ the 
multiples $a_iD_i$, where $1 \le i \le s$,
form a basis for $K^0$.
Moreover, cover $X$ by   
open affine subsets $X_1, \ldots, X_r \subseteq X$
such that for any $D \in K^0$ all 
restrictions $D_{\vert X_i}$ are principal. 
Let $\mathcal{S}^0$ be the sheaf of divisorial 
algebras associated to~$K^0$. 
Then $\Gamma(X_i, \mathcal{S}^0)$ is the Veronese
subalgebra of $\Gamma(X_i, \mathcal{S})$ defined 
by $K^0 \subseteq K$.
By Remark~\ref{rem:locprinc2split}, the algebra
$\Gamma(X_i, \mathcal{S}^0)$ is finitely generated.
Since $K^0 \subseteq K$ is of finite index, we 
can apply Proposition~\ref{prop:verofg2allfg} and 
obtain that $\Gamma(X_i, \mathcal{S})$ is finitely
generated.
\end{proof}

\begin{construction}
\label{constr:divrelspec}
Let $X$ be a normal prevariety, 
$K \subseteq \WDiv(X)$ a finitely 
generated subgroup and 
$\mathcal{S}$ the associated 
sheaf of divisorial algebras.
We assume that $\mathcal{S}$ is 
locally of finite type.
Then, in the notation of
Construction~\ref{constr:relspec},
the algebras $\Gamma(X_i,\mathcal{S})$ 
are $K$-graded.
This means that each affine variety $\t{X}_i$
comes with an action of the torus $H := \Spec\, \KK[K]$,
and, because of $\mathcal{S}_0 = \mathcal{O}_X$,
the canonical map $\t{X}_i \to X_i$ is a good quotient
for this action.
Since the whole gluing process is equivariant, 
we end up with an $H$-prevariety $\t{X} = \Spec_X \mathcal{S}$ 
and $p \colon \t{X} \to X$ is a good quotient 
for the $H$-action. 
\end{construction}

\begin{example}
Consider once more the projective line $X = \PP_1$,
the group $K := \ZZ D$, where $D := \{\infty\}$,
and the associated sheaf $\mathcal{S}$
of divisorial algebras.
For the affine charts $X_0 = \KK$ and 
$X_1 = \KK^* \cup \{\infty\}$
we have the graded isomorphisms
$$ 
\KK[T_0^{\pm 1},T_1] 
\ \to \
\Gamma(X_0, \mathcal{S}),
\qquad
\KK[T_0^{\pm 1},T_1]_n \ni f
\ \mapsto \ 
f(1,z) \in 
\Gamma(X_0, \mathcal{S}_{nD}),
$$ 
$$ 
\KK[T_0,T_1^{\pm 1}] 
\ \to \
\Gamma(X_1, \mathcal{S}),
\qquad
\KK[T_0,T_1^{\pm 1}]_n \ni f
\ \mapsto \ 
f(z,1) \in 
\Gamma(X_1, \mathcal{S}_{nD}).
$$ 
Thus, the corresponding spectra are 
$\KK^2_{T_0}$ and $\KK^2_{T_1}$.
The gluing takes place along 
$(\KK^*)^2$ and gives
$\t{X} = \KK^2 \setminus \{0\}$.
The action of $\KK^* = \Spec\, \KK[K]$
on $\t{X}$ is the usual 
scalar multiplication.
\end{example}

The above example fits into the more 
general context of sheaves of divisorial 
algebras associated to groups generated by 
a very ample divisor, i.e., the pullback
of a hyperplane with respect to an
embedding into a projective space.

\begin{example}
Suppose that $X$ is projective and
$K = \ZZ D$ holds with a
very ample divisor $D$ on $X$.
Then $\Gamma(X,\mathcal{S})$
is finitely generated and thus we have
the affine cone 
$\b{X} := \Spec \, \Gamma(X,\mathcal{S})$
over $X$.
It comes with a $\KK^*$-action and an 
attractive fixed point $\b{x}_0 \in \b{X}$,
i.e., $\b{x}_0$ lies in the closure of any
$\KK^*$-orbit.
The relative spectrum 
$\t{X} = \Spec_X \, \mathcal{S}$
equals $\b{X} \setminus \{\b{x}_0\}$.
\end{example}

\goodbreak

\begin{remark}
\label{rem:principalbundle}
In the setting of~\ref{constr:divrelspec},
let $U \subseteq X$ be an open subset such 
that all divisors $D \in K$ are principal 
over $U$.
Then there is a commutative diagram
of $H$-equivariant morphisms
$$ 
\xymatrix{
p^{-1}(U) 
\ar[rr]^{\cong}
\ar[dr]_{p}
& & 
H \times U
\ar[dl]^{\pr_U}
\\
& 
U 
& 
}
$$
where $H$ acts on $H \times U$ by multiplication
on the first factor.
In particular, if $K$ consists of Cartier divisors,
e.g. if $X$ is locally factorial, 
then  $\t{X} \to X$ is a locally trivial $H$-principal
bundle.
\end{remark}

\begin{proposition}
\label{prop:smallmap}
Situation as in Construction~\ref{constr:divrelspec}.
The prevariety $\t{X}$ is normal.
Moreover, for any closed $A \subseteq X$ 
of codimension at least two, 
$p^{-1}(A) \subseteq \t{X}$
is as well of codimension at least two.
\end{proposition}

\begin{proof}
For normality, we have to show 
that for every affine open $U \subseteq X$ 
the algebra $\Gamma(p^{-1}(U),\mathcal{O})$
is a normal ring.
According to Remark~\ref{rem:codim2divsheaf},
we have 
\begin{eqnarray*}
\Gamma(p^{-1}(U),\mathcal{O})
& = & 
\Gamma(p^{-1}(U_{\reg}),\mathcal{O}).
\end{eqnarray*}
Using Remark~\ref{rem:locprinc2split}, 
we see that the latter ring is normal.
The supplement is then an immediate 
consequence of Remark~\ref{rem:codim2divsheaf}.
\end{proof}

\subsection{Unique factorization in the global ring}
Here we investigate divisibility properties 
of the ring of 
global sections of the sheaf of divisorial 
algebras $\mathcal{S}$ 
associated to a subgroup $K \subseteq \WDiv(X)$ 
on a normal prevariety $X$.
The key statement is the following.

\begin{theorem}
\label{thm:fact1smooth}
Let $X$ be a smooth prevariety, 
$K \subseteq \WDiv(X)$ a finitely generated subgroup,
$\mathcal{S}$ the associated sheaf of divisorial algebras
and $\t{X} = \Spec_X \, \mathcal{S}$.
Then the following statements are equivalent.
\begin{enumerate}
\item
The canonical map $K \to \Cl(X)$ is surjective.
\item
The divisor class group $\Cl(\t{X})$ is trivial.
\end{enumerate}
\end{theorem}

We need a preparing observation concerning 
the pullback of Cartier divisors.
Recall that for any dominant morphism $\varphi \colon X \to Y$ 
of normal prevarieties, there is a pullback of Cartier 
divisors: if a Cartier divisor $E$ on $Y$ is locally 
given as $E = \div(g)$, then the pullback divisor 
$\varphi^*(E)$ is the Cartier divisor locally defined by 
$\div(\varphi^*(g))$.

\begin{lemma}
\label{smoothpull}
Situation as in Construction~\ref{constr:divrelspec}.
Suppose that $D \in K$ is Cartier 
and consider a non-zero section 
$f \in \Gamma(X,\mathcal{S}_D)$.
Then one has
\begin{eqnarray*}
p^*(D) 
& = & 
\div(f)  - p^*(\div(f)),
\end{eqnarray*}
where on the right hand side $f$ is firstly
viewed as a homogeneous function on $\t{X}$, 
and secondly as a rational function on $X$.
In particular, $p^*(D)$ is principal.
\end{lemma}

\begin{proof}
On suitable open sets $U_i \subseteq X$,
we find defining equations $f_i^{-1}$ for $D$ 
and thus may write $f = h_i f_i$, where 
$h_i \in \Gamma(U_i,\mathcal{S}_0) =  
\Gamma(U_i,\mathcal{O})$ and
$f_i \in \Gamma(U_i,\mathcal{S}_D)$. 
Then, on $p^{-1}(U_i)$, we have $p^*(h_i) = h_i$
and the function $f_i$ 
is homogeneous of degree $D$ and invertible.
Thus, we obtain
\begin{eqnarray*} 
p^*(D)
& = &
p^*(\div(f)+D) - p^*(\div(f))
\\ 
& = & 
p^*(\div(h_i))  - p^*(\div(f))
\\
& = & 
\div(h_i)  - p^*(\div(f))
\\
& = &  
\div(h_if_i)  - p^*(\div(f))
\\
& = & 
\div(f)  - p^*(\div(f)).
\end{eqnarray*}
\end{proof}

We are almost ready for proving the Theorem.
Recall that, given an action of an algebraic 
group $G$ on a normal prevariety $X$, we obtain
an induced action of $G$ on the group of 
Weil divisors by sending a prime divisor 
$D \subseteq X$ to $g \mal D \subseteq X$.
In particular, we can speak about invariant
Weil divisors.

\begin{proof}[Proof of Theorem~\ref{thm:fact1smooth}]
Suppose that~(i) holds.
It suffices
to show that every effective divisor 
$\t{D}$ on $\t{X}$ 
is principal.
We work with the action of the torus
$H = \Spec \, \KK[K]$ on $\t{X}$.
Choosing an $H$-linearization of $\t{D}$, 
see~\cite[Section~2.4]{KKLV},
we obtain a representation of $H$ on 
$\Gamma(\t{X},\mathcal{O}_{\t{X}}(\t{D}))$
such that for any section 
$\t{f} \in \Gamma(\t{X},\mathcal{O}_{\t{X}}(\t{D}))$
one has
\begin{eqnarray*}
\div_{\t{D}}(h \cdot \t{f})
& = & 
h \cdot \div_{\t{D}}(\t{f}).
\end{eqnarray*}
Taking a non-zero $\t{f}$, which is homogeneous 
with respect to this representation, 
we obtain that $\t{D}$ is linearly equivalent 
to the $H$-invariant divisor $\div_{\t{D}}(\t{f})$.
This reduces the problem to the case 
of an invariant divisor~$\t{D}$; 
compare also~\cite[Theorem~4.2]{Anders}.
Now, consider any invariant prime divisor $\t{D}$ 
on $\t{X}$.
Let $D := p(\t{D})$ be the image under 
the good quotient $p \colon \t{X} \to X$.
Remark~\ref{rem:principalbundle} 
gives $\t{D} = p^*(D)$.
By assumption, $D$ is linearly equivalent 
to a divisor $D' \in K$.
Thus, $\t{D}$ is linearly 
equivalent to $p^*(D')$, which in turn
is principal by Lemma~\ref{smoothpull}.

Now suppose that~(ii) holds. 
It suffices to show that any effective 
$D \in \WDiv(X)$ is linearly equivalent 
to some $D' \in K$.
The pullback $p^*(D)$ is the divisor of 
some function $f \in \Gamma(\t{X},\mathcal{O})$.
We claim that $f$ is $K$-homogeneous.
Indeed 
$$
F \colon H \times \t{X} \ \to \ \KK,
\qquad
(h,x) \ \mapsto \ f(h \mal x)/f(x)
$$
is an invertible function. 
By Rosenlicht's Lemma~\cite[Section~1.1]{KKV}, 
we have $F(h,x) = \chi(h) g(x)$
with $\chi \in \Chi(H)$ and 
$g \in \Gamma(\t{X}, \mathcal{O}^*)$. 
Plugging $(1,x)$ into $F$ yields 
$g=1$ and, consequently, $f(h \mal x) = \chi(h) f(x)$ 
holds.
Thus, we have $f \in \Gamma(X,\mathcal{S}_{D'})$ 
for some $D' \in K$.
Lemma~\ref{smoothpull} gives
$$
p^*(D) 
\ = \
\div(f)
\ = \ 
p^*(D') + p^*(\div(f)),
$$ 
where in the last term, $f$ is regarded as 
a rational function on $X$. We conclude 
$D = D' + \div(f)$ on $X$.
In other words, $D$ is linearly equivalent to 
$D' \in K$.
\end{proof}

As an immediate consequence, we obtain factoriality
of the ring of global sections
provided $K \to \Cl(X)$ is surjective,
see also~\cite{BeHa1}, \cite{ElKuWa} and~\cite{Ar2}.

\begin{theorem}
\label{thm:fact1}
Let $X$ be a normal prevariety, $K \subseteq \WDiv(X)$ 
a finitely generated subgroup and $\mathcal{S}$ the 
associated sheaf of divisorial algebras.
If the canonical map $K \to \Cl(X)$ is surjective,
then the algebra
$\Gamma(X,\mathcal{S})$ is a unique factorization 
domain.
\end{theorem}

\begin{proof}
According to Remark~\ref{rem:codim2divsheaf},
the algebra $\Gamma(X,\mathcal{S})$ equals 
$\Gamma(X_{\reg},\mathcal{S})$ and thus 
we may apply Theorem~\ref{thm:fact1smooth}.
\end{proof}

Divisibility and primality in 
the ring of global sections 
$\Gamma(X,\mathcal{S})$
can be characterized purely in terms 
of~$X$.

\begin{proposition}
\label{cor:divdivalg}
Let $X$ be a normal prevariety, 
$K \subseteq \WDiv(X)$ a finitely 
generated subgroup projecting onto 
$\Cl(X)$ and let $\mathcal{S}$ be the 
associated sheaf of divisorial algebras.
\begin{enumerate}
\item
An element $0 \ne f \in \Gamma(X,\mathcal{S}_D)$ divides 
an element $0 \ne g \in \Gamma(X,\mathcal{S}_E)$ 
if and only if
$\div_D(f) \le \div_E(g)$ holds.
\item
An element $0 \ne f \in \Gamma(X,\mathcal{S}_D)$ is prime 
if and only if the divisor $\div_D(f) \in \WDiv(X)$ 
is prime.
\end{enumerate}
\end{proposition}

\begin{proof}
We may assume that $X$ is smooth.
Then $\t{X} = \Spec_X \, \mathcal{S}$ 
exists, and Lemma~\ref{smoothpull}
reduces (i) and~(ii) to the corresponding 
statements on regular functions on $\t{X}$, 
which in turn are well known.
\end{proof}

\subsection{Geometry of the relative spectrum}
We collect basic geometric properties 
of the relative spectrum of a sheaf of 
divisorial algebras. 
We will use the following pullback
construction for Weil divisors.

\begin{remark}
\label{rem:weilpull}
Consider any 
morphism $\varphi \colon \t{X} \to X$ 
of normal prevarieties
such that the closure of 
$X \setminus \varphi(\t{X})$ is of 
codimension at least two in $X$.
Then we may define a pullback 
homomorphism for Weil divisors
$$ 
\varphi^* \colon \ \WDiv(X) \ \to \ \WDiv(\t{X})
$$
as follows: 
Given $D \in \WDiv(X)$,
consider its restriction
$D'$ to $X_{\reg}$, the 
usual pullback $\varphi^*(D')$
of Cartier divisors
on $\varphi^{-1}(X_{\reg})$
and define $\varphi^*(D)$ 
to be the Weil divisor obtained 
by closing the support of 
$\varphi^*(D')$. Note that we 
always have
\begin{eqnarray*}
\Supp(\varphi^*(D))
& \subseteq &
\varphi^{-1}(\Supp(D)).
\end{eqnarray*}
If for any closed $A \subseteq X$ 
of codimension at least two, 
$\varphi^{-1}(A) \subseteq \t{X}$
is as well of codimension at least two,
then $\varphi^*$ maps 
principal divisors to 
principal divisors, and we obtain 
a pullback homomorphism
$$ 
\varphi^* \colon \ \Cl(X) \ \to \ \Cl(\t{X}).
$$
\end{remark}

\begin{example}
Consider $X = V(\KK^4; \; T_1T_2 - T_3T_4)$ 
and $\t{X} = \KK^4$. Then we have a morphism
$$ 
p \colon \t{X} \ \to \ X, 
\qquad
z \ \mapsto \ 
(z_1z_2,z_3z_4, z_1z_3,z_2z_4).
$$
For the prime divisor $D = \KK \times 0 \times \KK \times 0$
on $X$, we have 
$$ 
\Supp(p^*(D))
\ = \ 
V(\t{X};Z_4)
\ \subsetneq \
V(\t{X};Z_4)
\cup 
V(\t{X};Z_2,Z_3)
\ = \ 
p^{-1}(\Supp(D)).
$$
In fact, $p \colon \t{X} \to X$ is the 
morphism determined by the sheaf of divisorial algebras
associated to $K = \ZZ D$.
\end{example}

\index{affine intersection!}
We say that a prevariety $X$ is of 
{\em affine intersection\/}  if
for any two affine open subsets
$U,U' \subseteq X$ the intersection $U \cap U'$ 
is again affine.
For example, every variety is of affine 
intersection.
Note that a prevariety $X$ is of affine 
intersection if it can be covered 
by open affine subsets 
$X_1, \ldots, X_s \subseteq X$ 
such that all intersections 
$X_i \cap X_j$ are affine.
Moreover, if $X$ is of affine intersection, 
then the complement of any affine open subset
$U \subsetneq X$ is of pure codimension one.

\begin{proposition}
\label{prop:divalgspec}
In the situation of~\ref{constr:divrelspec},
consider the pullback homomorphism
$p^* \colon \WDiv(X) \to \WDiv(\t{X})$
defined in~\ref{rem:weilpull}.
Then, for every 
$D \in K$ and every non-zero
$f \in \Gamma(X,\mathcal{S}_D)$,
we have 
\begin{eqnarray*}
\div(f) & = & p^*(\div_D(f)),
\end{eqnarray*}
where on the left hand side 
$f$ is a function on $\t{X}$, 
and on the right hand side a 
function on $X$.
If $X$ is of affine intersection and 
$X_{D,f}$ is affine, then we have moreover
\begin{eqnarray*}
\Supp(\div(f)) & = & p^{-1}(\Supp(\div_D(f))).
\end{eqnarray*}
\end{proposition}

\begin{proof}
By Lemma~\ref{smoothpull}, the first 
equation holds on $p^{-1}(X_{\reg})$.
By Proposition~\ref{prop:smallmap},
the complement $\t{X} \setminus p^{-1}(X_{\reg})$
is of codimension at least two
and thus the first equation holds 
on the whole $\t{X}$.
For the proof of the second one, consider
$$ 
X_{D,f} 
\ =  \
X \setminus \Supp(\div_D(f)),
\qquad\qquad
\t{X}_f 
\ =  \
\t{X} \setminus V(\t{X},f).
$$
Then we have to show that 
$p^{-1}(X_{D,f})$ equals $\t{X}_f$.
Since $f$ is invertible on $p^{-1}(X_{D,f})$,
we obtain $p^{-1}(X_{D,f}) \subseteq \t{X}_f$.
Moreover, Lemma~\ref{smoothpull} yields
\begin{eqnarray*}
p^{-1}(X_{D,f}) \ \cap \ p^{-1}(X_{\reg})
& = & 
\t{X}_f \cap p^{-1}(X_{\reg}).
\end{eqnarray*}
Thus the complement $\t{X}_f \setminus p^{-1}(X_{D,f})$
of the affine subset $p^{-1}(X_{D,f}) \subseteq \t{X}_f$
is of codimension at least two.
Since $p \colon \t{X} \to X$ is affine, 
the prevariety $\t{X}$ inherits the property 
to be of affine intersection from $X$
and hence $\t{X}_f \setminus p^{-1}(X_{D,f})$
must be empty.
\end{proof}

\begin{corollary}
\label{cor:homogzeroes}
Situation as in Construction~\ref{constr:divrelspec}.
Let $\t{x} \in \t{X}$ be a point such that 
$H \mal \t{x} \subseteq \t{X}$ is closed,
and let $0 \ne f \in \Gamma(X,\mathcal{S}_D)$.
Then we have 
\begin{eqnarray*}
f(\t{x}) \ = \ 0
& \iff &
p(\t{x}) \ \in \ \Supp(\div_D(f)).
\end{eqnarray*}
\end{corollary}

\begin{proof}
Remark~\ref{rem:weilpull} and 
Proposition~\ref{prop:divalgspec} 
show that $p(\Supp(\div(f)))$ 
is contained in $\Supp(\div_D(f))$.
Moreover, they coincide along 
the smooth locus of $X$ and
Theorem~\ref{prop:GclosGsep}
ensures that $p(\Supp(\div(f)))$ is 
closed. This gives
\begin{eqnarray*}
p(\Supp(\div(f)))
& = & 
\Supp(\div_D(f)).
\end{eqnarray*}  
Thus, $f(\t{x}) = 0$ implies
$p(\t{x}) \in \Supp(\div_D(f))$.
If $p(\t{x}) \in \Supp(\div_D(f))$
holds, then some 
$\t{x}' \in \Supp(\div(f))$ 
lies in the $p$-fiber of $\t{x}$.
Since $H \mal \t{x}$ is closed,
it is contained in the closure of 
$H \mal \t{x}'$, see 
Corollary~\ref{prop:goodquotfibers}.
This implies $\t{x} \in \Supp(\div(f))$.
\end{proof}


\begin{corollary}
\label{cor:quaff-1}
Situation as in Construction~\ref{constr:divrelspec}.
If $X$ is of affine intersection and 
covered by affine open subsets 
of the form $X_{D,f}$, where $D \in K$ and 
$f \in \Gamma(X,\mathcal{S}_D)$,
then $\t{X}$ is a quasiaffine variety.
\end{corollary}

\begin{proof}
According to Proposition~\ref{prop:divalgspec},
the prevariety $\t{X}$ is covered by open 
affine subsets of the form $\t{X}_f$ and thus 
is quasiaffine.
\end{proof}

\begin{corollary}
\label{cor:quaff2}
Situation as in Construction~\ref{constr:divrelspec}.
If $X$ is of affine intersection and 
$K \to \Cl(X)$ is surjective, 
then $\t{X}$ is a quasiaffine variety.
\end{corollary}

\begin{proof}
Cover $X$ by affine open sets $X_1, \ldots, X_r$.
Since $X$ is of affine intersection, every complement
$X \setminus X_i$ is of pure codimension one.
Since $K \to \Cl(X)$ is surjective, we obtain 
that $X \setminus X_i$ is the support of the
$D$-divisor of some $f \in \Gamma(X,\mathcal{S}_D)$.
The assertion thus follows from Corollary~\ref{cor:quaff-1}.
\end{proof}

\begin{proposition}
\label{prop:isogrp2divisors}
Situation as in Construction~\ref{constr:divrelspec}.
For $x \in X$, let $K_x^0 \subseteq K$
be the subgroup of divisors that are 
principal near $x$ and let 
$\t{x} \in p^{-1}(x)$ be a point with 
closed $H$-orbit.
Then the isotropy group $H_{\t{x}} \subseteq H$ 
is given by $H_{\t{x}} = \Spec \, \KK[K/K_x^0]$.
\end{proposition}

\begin{proof}
Replacing $X$ with a suitable affine 
neighbourhood of $x$, we may assume 
that $\t{X}$ is affine.
By Proposition~\ref{prop:isogrorbitlattice},
the isotropy group $H_{\t{x}}$ is 
$\Spec \, \KK[K/K_{\t{x}}]$
with the orbit group
$$
K_{\t{x}}
\ = \
\bangle{D \in K; \; f(\t{x}) \ne 0 
              \text{ for some } f \in \Gamma(X,\mathcal{S}_D)}
\ \subseteq \ 
K.
$$
Using Corollary~\ref{cor:homogzeroes}, 
we obtain that 
there exists an $f \in \Gamma(X,\mathcal{S}_D)$ 
with $f(\t{x}) \ne 0$ if and only if $D \in K_x^0$
holds. The assertion follows.
\end{proof}


\begin{corollary}
Situation as in Construction~\ref{constr:divrelspec}.
\begin{enumerate}
\item
If $X$ is locally factorial, then $H$ 
acts freely on $\t{X}$.
\item
If $X$ is $\QQ$-factorial, then $H$ 
acts with at most finite isotropy 
groups on $\t{X}$.
\end{enumerate}
\end{corollary}


\section{Cox sheaves and Cox rings}
\label{sec:coxrings}

\subsection{Free divisor class group}
As before, we work over an algebraically 
closed field $\KK$ of characteristic zero.
We introduce Cox sheaves and Cox rings 
for a prevariety with a free finitely generated 
divisor class group.
As an example, we compute 
in~\ref{exam:p1withdoubledpoint}
the Cox ring of a non-separated curve, 
the projective line with multipled points.

\begin{construction}
\label{constr:crfree}
\index{Cox sheaf!}%
\index{Cox ring!}%
Let $X$ be a normal prevariety
with free finitely generated divisor 
class group $\Cl(X)$.
Fix a subgroup $K \subseteq \WDiv(X)$
such that the canonical map 
$c \colon K \to \Cl(X)$ 
sending $D \in K$ to its class $[D] \in \Cl(X)$
is an isomorphism.
We define the {\em Cox sheaf\/} 
associated to $K$ to be
$$ 
\mathcal{R}
\ := \ 
\bigoplus_{[D] \in \Cl(X)} \mathcal{R}_{[D]},
\qquad\qquad
\mathcal{R}_{[D]}
\ := \ 
\mathcal{O}_X(D),
$$
where $D \in K$ represents $[D] \in \Cl(X)$
and the multiplication in $\mathcal{R}$ is 
defined by multiplying homogeneous sections
in the field of rational functions $\KK(X)$.
The sheaf $\mathcal{R}$ is a quasicoherent 
sheaf of normal integral $\mathcal{O}_X$-algebras 
and, up to isomorphy, it does not depend on 
the choice of the subgroup $K \subseteq \WDiv(X)$.
The {\em Cox ring\/} of $X$ is the algebra of global
sections
$$ 
\mathcal{R}(X)
\ := \ 
\bigoplus_{[D] \in \Cl(X)} \mathcal{R}_{[D]}(X),
\qquad\qquad
 \mathcal{R}_{[D]}(X)
\ := \
\Gamma(X,\mathcal{O}_X(D)).
$$
\end{construction}

\begin{proof}[Proof of Construction~\ref{constr:crfree}]
Given two subgroups $K,K' \subseteq \WDiv(X)$
projecting isomorphically onto $\Cl(X)$,
we have to show that the corresponding 
sheaves of divisorial algebras $\mathcal{R}$ 
and $\mathcal{R}'$ are isomorphic.
Choose a basis $D_1, \ldots, D_s$ for $K$
and define a homomorphism
$$
\eta \colon K \ \to \ \KK(X)^*,
\qquad
a_1D_1 + \ldots +  a_sD_s \ \mapsto \ f_1^{a_1} \cdots f_s^{a_s},
$$
where  $f_1, \ldots, f_s \in \KK(X)^*$ are 
such that the divisors $D_i - \div(f_i)$
form a basis of $K'$.
Then we obtain an isomorphism $(\psi,\t{\psi})$
of the sheaves of divisorial algebras 
$\mathcal{R}$ and $\mathcal{R}'$ by setting
$$
\begin{array}{rclcrcl}
\t{\psi} \colon K & \to & K', 
& \qquad &
D  & \mapsto &  - \div(\eta(D)) + D,
\\[1ex]
\psi \colon \mathcal{R} 
& \to & 
\mathcal{R}',
& \qquad &
\Gamma(U,\mathcal{R}_{[D]}) \ni f 
& \mapsto &  
\eta(D) f \in \Gamma(U,\mathcal{R}_{[\t{\psi}(D)]}).
\end{array}
$$
\end{proof}

\begin{example}
Let $X$ be the projective space $\PP_n$ and 
$D \subseteq \PP_n$ be a hyperplane. 
The class of $D$ generates $\Cl(\PP_n)$ freely. 
We take $K$ as the subgroup of $\WDiv(\PP_n)$ 
generated by $D$, and the Cox ring 
$\mathcal{R}\,(\PP_n)$ is the polynomial
ring $\KK[z_0, z_1, \ldots, z_n]$ 
with the standard grading. 
\end{example}

\begin{remark}
If $X \subseteq \PP_n$ 
is a closed normal subvariety whose divisor class
group is generated by a hyperplane section,
then $\mathcal{R}(X)$ coincides with 
$\Gamma(\b{X}, \mathcal{O})$,
where $\b{X} \subseteq \KK^{n+1}$ 
is the cone over $X$ if and only if $X$ 
is projectively normal.
\end{remark}

\begin{remark}
Let $s$ denote the rank of $\Cl(X)$.
Then Remark~\ref{RR} realizes the Cox ring 
$\mathcal{R}(X)$ as a graded subring 
of the Laurent polynomial ring:
\begin{eqnarray*}
\mathcal{R}(X)
& \subseteq & 
\KK(X)[T_1^{\pm 1}, \ldots, T_s^{\pm 1}].
\end{eqnarray*}
Using the fact that there are 
$f \in \mathcal{R}_{[D]}(X)$ with 
$X_{D,f}$ affine and Remark~\ref{rem:Sloc},
we see that this inclusion gives rises 
to an isomorphism of the quotient fields
\begin{eqnarray*}
\Quot({\mathcal R}(X)) 
& \cong & 
\KK(X)(T_1, \ldots, T_s).
\end{eqnarray*}
\end{remark}

\begin{proposition}
\label{prop:cralgfree}
Let $X$ be a normal prevariety with 
free finitely generated divisor class group.
\begin{enumerate}
\item
The Cox ring $\mathcal{R}(X)$ is a unique 
factorization domain.
\item
The units of the Cox ring 
are given by $\mathcal{R}(X)^* = \Gamma(X,\mathcal{O}^*)$.
\end{enumerate}
\end{proposition}

\begin{proof}
The first assertion is a direct consequence 
of Theorem~\ref{thm:fact1}.
To verify the second one, consider a unit 
$f \in \mathcal{R}(X)^*$.
Then $fg = 1 \in \mathcal{R}_0(X)$ holds
with some unit $g \in \mathcal{R}(X)^*$.
This can only happen, when $f$ and $g$ 
are homogeneous, say of degree $[D]$ and 
$-[D]$, and thus we obtain
$$
0 
\; = \; 
\div_0(1) 
\; = \; 
\div_{D}(f) + \div_{-D}(g) 
\; = \;
(\div(f) + D) + (\div(g) - D). 
$$
Since the divisors $(\div(f) + D)$ and $(\div(g) - D)$
are effective, we conclude that $D = -\div(f)$.
This means $[D] = 0$ and we obtain  
$f \in \Gamma(X,\mathcal{O}^*)$.
\end{proof}

\begin{example} 
\label{exam:p1withdoubledpoint}
Compare~\cite[Section~2]{HaSu}.
Take the projective line $\PP_1$, 
a~tuple $A = ( a_0, \ldots, a_r)$ 
of pairwise different points $a_i \in \PP_1$
and a tuple $\mathfrak{n}=(n_0, \ldots, n_r)$ 
of integers $n_i \in \ZZ_{\ge 1}$. 
We construct a non-separated smooth 
curve $\PP_1(A,\mathfrak{n})$ mapping birationally 
onto $\PP_1$ such that over each $a_i$ lie precisely $n_i$ 
points. Set
$$
X_{ij} 
\ := \ 
\PP_1 \setminus \bigcup_{k \ne i} a_k, 
\qquad 
0 \le i \le r,
\qquad
1\le j \le n_i.
$$
Gluing the $X_{ij}$ along the common open subset 
$\PP_1 \setminus \{ a_0, \ldots, a_r \}$
gives an irreducible smooth prevariety 
$\PP_1(A,\mathfrak{n})$ of dimension one. 
The inclusion maps $X_{ij} \to \PP_1 $ define  
a morphism $\pi \colon \PP_1(A, \mathfrak{n}) \to \PP_1$, 
which is locally an isomorphism. 
Writing $a_{ij}$ for the point in $\PP_1(A, \mathfrak{n})$ 
stemming from $a_i \in X_{ij}$, 
we obtain the fibre over any $a \in \PP_1$ as
\begin{eqnarray*}
\pi^{-1}(a)
& = &
\begin{cases}
\{a_{i1}, \ldots, a_{in_i}\} & a = a_i \text{ for some } 0 \le i \le r,
\\
\{a\}                        & a \ne a_i \text{ for all } 0 \le i \le r.
\end{cases}
\end{eqnarray*}

We compute the divisor class group of 
$\PP_1(A,\mathfrak{n})$.
Let $K'$ denote the group of Weil divisors on
$\PP_1(A,\mathfrak{n})$ generated by the prime 
divisors $a_{ij}$.
Clearly $K'$ maps onto the divisor class group.
Moreover, the group of principal divisors inside
$K'$ is
$$
K'_0
\ := \ 
K'  \cap \PDiv(\PP_1(A,\mathfrak{n}))
\ = \
\left\{ 
\sum_{{0 \le i \le r,} \atop {1 \le j \le n_i}} c_i a_{ij}; \; 
c_0 + \ldots + c_r = 0
\right\}.
$$
One directly checks that $K'$ 
is the direct sum of $K_0'$ and the 
subgroup $K \subseteq K'$ generated by 
$a_{01}, \ldots, a_{0n_0}$ and the 
$a_{i1}, \ldots, a_{in_i-1}$.
Consequently, the divisor class group of 
$\PP_1(A, \mathfrak{n})$ is given by
\begin{eqnarray*}
\Cl(\PP_1(A,\mathfrak{n}))
& = &
\bigoplus_{j=1}^{n_0} \ZZ \mal [a_{0j}] 
\ \oplus \
\bigoplus_{i=1}^r
\left( \bigoplus_{j=1}^{n_i-1} \ZZ \mal [a_{ij}] 
\right).
\end{eqnarray*}

We are ready to determine the Cox ring 
of the prevariety  
$\PP_1(A, \mathfrak{n})$. 
For every $0 \le i \le r$,
define a monomial
$$
T_i 
\ := \ 
T_{i1} \cdots T_{in_i}
\ \in \ 
\KK[T_{ij}; \; 0 \le i \le r, \; 1 \le j \le n_i].
$$
Moreover, for every $a_i \in \PP_1$ 
fix a presentation 
$a_i = [b_i,c_i]$ with $b_i,c_i \in \KK$
and for every $0 \le i \le r-2$ set 
$k = j+1 = i+2$ 
and define a trinomial
\begin{eqnarray*}
g_i 
& := &
(b_jc_k - b_kc_j)T_i
\ + \ 
(b_kc_i - b_ic_k)T_j
\ + \ 
(b_ic_j - b_jc_i)T_k.
\end{eqnarray*}
We claim that for $r \le 1$ the Cox ring 
$\mathcal{R}\,(\PP_1(A, \mathfrak{n}))$
is isomorphic to the polynomial ring $\KK[T_{ij}]$, 
and for $r \ge 2$ it has a presentation
\begin{eqnarray*}
\mathcal{R}(\PP_1(A,\mathfrak{n}))
&  \cong &
\KK[T_{ij}; \; 0 \le i \le r, \; 1 \le j \le n_i]
\  / \ 
\bangle{g_i; \; 0 \le i \le r-2},
\end{eqnarray*}
where, in both cases, the grading is given by 
$\deg(T_{ij}) = [a_{ij}]$. 
Note that all relations are 
homogeneous of degree
$$ 
\deg(g_i)
\ = \ 
[a_{i1} + \ldots + a_{in_i}]
\ = \ 
[a_{01} + \ldots + a_{0n_0}].
$$

Let us verify this claim.
Set for short $X := \PP_1(A,\mathfrak{n})$ and $Y := \PP_1$.
Let $K \subseteq \WDiv(X)$ be the subgroup
generated by all $a_{ij} \in X$ different from 
$a_{1n_1}, \ldots, a_{rn_r}$,
and let $L \subseteq \WDiv(Y)$ be the subgroup
generated by $a_0 \in Y$.
Then we may view the Cox rings $\mathcal{R}(X)$ 
and $\mathcal{R}(Y)$ as the rings of 
global sections of the sheaves of divisorial algebras
$\mathcal{S}_X$ and $\mathcal{S}_Y$ associated 
to $K$ and $L$.
The canonical morphism $\pi \colon X \to Y$ 
gives rise to injective pullback homomorphisms
$$
\pi^* \colon L \ \to \ K,
\qquad
\qquad
\pi^* \colon \Gamma(Y,\mathcal{S}_Y) \ \to \ \Gamma(X,\mathcal{S}_X).
$$
For any divisor $a_{ij} \in K$, 
let $T_{ij} \in \Gamma(X,\mathcal{S}_X)$ 
denote its canonical section, 
i.e., the rational function $1 \in \Gamma(X,\mathcal{S}_{X,a_{ij}})$.
Moreover, let $[z,w]$ be the homogeneous coordinates
on $\PP_1$ and consider the sections
$$
S_i 
\ := \ 
\frac{b_iw-c_iz}{b_0w-c_0z}
\ \in \
\Gamma(Y,\mathcal{S}_{Y,a_0}),
\qquad
0 \le i \le r.
$$
Finally, set 
$d_{in_i} :=  a_{01}+ \ldots + a_{0n_0} - a_{i1} -  \ldots - a_{in_i-1} \in K$ 
and define homogeneous sections
$$
T_{in_i} 
\ := \ 
\pi^* S_i (T_{i1} \cdots T_{in_i-1})^{-1}
\ \in \
\Gamma(X,\mathcal{S}_{X,d_{in_i}}),
\qquad
1 \le i \le r.
$$ 
We show that the sections $T_{ij}$, 
where $0 \le i \le r$
and $1 \le j \le n_i$, generate the 
Cox ring $\mathcal{R}(X)$.
Note that we have 
$$
\div_{a_{ij}} (T_{ij}) \ = \ a_{ij}, 
\qquad
\div_{d_{in_i}} (T_{in_i}) \ = \ a_{in_i}.
$$  
Consider $D \in K$ and $h \in \Gamma(X,\mathcal{S}_D)$.
If there occurs an $a_{ij}$ in $\div_D(h)$, 
then we may divide $h$ in $\Gamma(X,\mathcal{S})$ 
by the corresponding $T_{ij}$,
use Proposition~\ref{cor:divdivalg}~(i).
Doing this as long as possible, we arrive 
at some $h' \in \Gamma(X,\mathcal{S}_{D'})$
such that $\div_{D'}(h')$
has no components $a_{ij}$.
But then $D'$ is a pullback divisor and hence
$h'$ is contained in 
$$
\pi^*(\Gamma(Y,\mathcal{S}_Y)) 
\ = \
\KK[\pi^*S_0,\pi^*S_1]
\ = \ 
\KK[T_{01} \cdots T_{0n_0}, T_{11} \cdots T_{1n_1}].
$$
Finally, we have to determine the relations
among the sections 
$T_{ij} \in \Gamma(X,\mathcal{S}_X)$.
For this, we first note that 
among the $S_i \in \Gamma(Y,\mathcal{S}_Y)$
we have the relations
$$
(b_jc_k - b_kc_j)S_i
\ + \ 
(b_kc_i - b_ic_k)S_j
\ + \ 
(b_ic_j - b_jc_i)S_k
\ = \ 
0,
\quad
j \ = \ i+1, 
\quad
k \ = \ i +2.
$$
Given any nontrivial homogeneous relation
$F = \alpha_1F_1 + \ldots + \alpha_lF_l = 0$
with $\alpha_i \in \KK$ and 
pairwise different monomials $F_i$ 
in the $T_{ij}$,
we achieve by subtracting suitable multiples 
of pullbacks of the above  relations 
a homogeneous relation
$$
F' 
\ = \ 
\alpha_1'F''_1 \pi^*S_0^{k_1}\pi^*S_1^{l_1} 
+ \ldots  + 
\alpha_m'F''_m \pi^*S_0^{k_m}\pi^*S_1^{l_m} 
\ = \ 
0
$$
with pairwise different monomials 
$F_j''$, none of which has any factor $\pi^*S_i$.
We show that $F'$ must be trivial.
Consider the multiplicative group $M$ of 
Laurent monomials in the $T_{ij}$ and 
the degree map
$$ 
M \ \to \ K,
\qquad
T_{ij}
\ \mapsto \ 
\deg(T_{ij})
\ = \ 
\begin{cases}
a_{ij},
& 
i = 0 \text{ or } j \le  n_i-1,
\\
d_{in_i}
& 
i \ge 1 \text{ and } j = n_i.
\end{cases}
$$
The kernel of this degree map is generated by 
the Laurent monomials $\pi^*S_0/\pi^*S_i$, where 
$1 \le i \le r$.
The monomials of $F'$ are all of the 
same $K$-degree and thus any two of them 
differ by a product of (integral) powers of the 
$\pi^*S_i$. 
It follows that all the $F''_j$ coincide.
Thus, we obtain the relation 
\begin{eqnarray*}
\alpha_1'\pi^*S_0^{k_1}\pi^*S_1^{l_1} 
\ + \ \ldots \ + \ 
\alpha_m'\pi^*S_0^{k_m}\pi^*S_1^{l_m} 
& = &
0.
\end{eqnarray*}
This relation descends to a relation 
in $\Gamma(Y,\mathcal{S}_Y)$, which is 
the polynomial ring $\KK[S_0,S_1]$.
Consequently, we obtain 
$\alpha_1' = \ldots = \alpha_m' = 0$.
\end{example}

\subsection{Torsion in the class group}
Again we work over an algebraically closed field
$\KK$ of characteristic zero.
We extend the definition of Cox sheaf and 
Cox ring to normal prevarieties~$X$ having
a finitely generated divisor class group~$\Cl(X)$
with torsion.
The idea is to take a subgroup 
$K \subseteq \WDiv(X)$ 
projecting onto~$\Cl(X)$,
to consider its associated sheaf of 
divisorial algebras $\mathcal{S}$ 
and to identify in a systematic
manner  homogeneous components $\mathcal{S}_D$ 
and $\mathcal{S}_{D'}$, 
whenever $D$ and $D'$ are
linearly equivalent.

\begin{construction}
\label{constr:crtorsion}
\index{Cox ring!}\index{Cox sheaf!}%
Let $X$ be a normal prevariety with 
$\Gamma(X,\mathcal{O}^*) = \KK^*$
and finitely generated divisor class 
group $\Cl(X)$.
Fix a subgroup 
$K \subseteq \WDiv(X)$ 
such that the map 
$c \colon K \to \Cl(X)$
sending $D \in K$ to its class 
$[D] \in \Cl(X)$
is surjective.
Let $K^0 \subseteq K$ be the kernel 
of $c$,
and let $\chi \colon K^0 \to \KK(X)^*$
be a character, i.e.~a group homomorphism, 
with
$$
\div(\chi(E))
\ = \ 
E,
\qquad
\text{for all } 
E \in K^0.
$$
Let $\mathcal{S}$ be the sheaf of divisorial algebras 
associated to $K$ and denote by $\mathcal{I}$ 
the sheaf of ideals of $\mathcal{S}$
locally generated by the sections
$1 - \chi(E)$, where $1$ is homogeneous
of degree zero, $E$ runs through $K^0$ and 
$\chi(E)$ is homogeneous of degree $-E$.
The {\em Cox sheaf\/} associated to 
$K$ and $\chi$ is the quotient sheaf 
$\mathcal{R} := \mathcal{S}/\mathcal{I}$
together with the $\Cl(X)$-grading
$$ 
\mathcal{R}
\ = \
\bigoplus_{[D] \in \Cl(X)}  \mathcal{R}_{[D]},
\qquad\qquad
\mathcal{R}_{[D]} 
\ := \ 
\pi \left( \bigoplus_{D' \in c^{-1}([D])} \mathcal{S}_{D'} \right).
$$
where $\pi \colon \mathcal{S} \to \mathcal{R}$ 
denotes the projection.
The Cox sheaf $\mathcal{R}$ is a quasicoherent sheaf 
of $\Cl(X)$-graded $\mathcal{O}_X$-algebras.
The {\em Cox ring\/} is the ring of global 
sections
$$ 
\mathcal{R}(X)
\ := \
\bigoplus_{[D] \in \Cl(X)}  \mathcal{R}_{[D]}(X),
\qquad\qquad
\mathcal{R}_{[D]}(X)
\ := \ 
\Gamma(X,\mathcal{R}_{[D]}).
$$
For any open set $U \subseteq X$, the canonical 
homomorphism 
$\Gamma(U,\mathcal{S})/\Gamma(U,\mathcal{I}) 
\to \Gamma(U,\mathcal{R})$
is an isomorphism.
In particular, we have 
\begin{eqnarray*}
\mathcal{R}(X)
& \cong & 
\Gamma(X,\mathcal{S})/\Gamma(X,\mathcal{I}).
\end{eqnarray*}
\end{construction}

All the claims made in this construction 
will be verified as separate Lemmas 
in the next subsection. 
The assumption $\Gamma(X,\mathcal{O}^*) = \KK^*$ is 
crucial for the following uniqueness statement
on Cox sheaves and rings.

\begin{proposition}
\label{prop:coxsheafunique}
Let $X$ be a normal prevariety with 
$\Gamma(X,\mathcal{O}^*) = \KK^*$
and finitely generated divisor class 
group $\Cl(X)$.
If $K,\chi$ and $K',\chi'$ are data
as in Construction~\ref{constr:crtorsion},
then there is a graded isomorphism 
of the associated Cox sheaves.
\end{proposition}

Also this will be proven in the next subsection.
The construction of Cox sheaves 
(and thus also Cox rings) of a prevariety
$X$ can be made canonical by fixing a suitable 
point $x \in X$.

\begin{construction}
\index{Cox ring! canonical}\index{Cox sheaf! canonical}%
Let $X$ be a normal prevariety with 
$\Gamma(X,\mathcal{O}^*) = \KK^*$
and finitely generated divisor class group $\Cl(X)$. 
Fix a point $x \in X$ with factorial local 
ring $\mathcal{O}_{X,x}$.
For the subgroup 
\begin{eqnarray*}
K^{x} 
& := & 
\{D \in \WDiv(X); \ x \not\in \Supp(D)\}
\end{eqnarray*}
let $\mathcal{S}^{x}$ be the associated
sheaf of divisorial algebras and
let $K^{x,0} \subseteq K^{x}$ denote 
the subgroup
consisting of principal divisors.
Then, for each $E \in K^{x,0}$, 
there is a unique section
$f_E \in \Gamma(X,\mathcal{S}_{-E})$,
which is defined near $x$ and satisfies
$$ 
\div(f_E) \ = \ E,
\qquad\qquad
f_E(x) \ = \ 1.
$$ 
The map $\chi^x \colon K^x \to \KK(X)^*$ 
sending $E$ to $f_E$ is a character
as in Construction~\ref{constr:crtorsion}.
We call the Cox sheaf $\mathcal{R}^x$ 
associated to $K^x$ and $\chi^x$ the 
{\em canonical Cox sheaf of the pointed
space $(X,x)$.}
\end{construction}


\begin{example}[An affine surface with torsion in the divisor class group]
Consider the two-dimensional affine quadric
$$
X 
\ := \ 
V(\KK^3 ; \; T_1T_2 - T_3^2)
\ \subseteq \
\KK^3.
$$
We have the functions $f_i := T_{i \vert X}$
on $X$ and with the prime divisors 
$D_1 := V(X ; f_1)$ and $D_2 := V(X ; f_2)$
on $X$, we have 
$$ 
\div(f_1) \ = \ 2D_1,
\qquad
\div(f_2) \ = \ 2D_2,
\qquad
\div(f_3) \ = \ D_1+D_2.
$$
The divisor class group $\Cl(X)$ is of order two; 
it is generated by $[D_1]$.
For $K := \ZZ D_1$, let $\mathcal{S}$ 
denote the associated sheaf of divisorial algebras.
Consider the sections 
$$
g_1 := 1 \ \in \ \Gamma(X,\mathcal{S}_{D_1}),
\qquad 
g_2 := f_3f_1^{-1} \ \in \ \Gamma(X,\mathcal{S}_{D_1}),
$$
$$
g_3 := f_1^{-1} \ \in \ \Gamma(X,\mathcal{S}_{2D_1}),
\qquad
g_4 := f_1 \ \in \ \Gamma(X,\mathcal{S}_{-2D_1}).
$$
Then $g_1,g_2$ generate $\Gamma(X,\mathcal{S}_{D_1})$
as a $\Gamma(X,\mathcal{S}_{0})$-module, 
and $g_3,g_4$ are inverse to each other. 
Moreover, we have 
$$ 
f_1 \ = \ g_1^2g_4,
\qquad
f_2 \ = \ g_2^2g_4,
\qquad
f_3 \ = \ g_1g_2g_4.
$$
Thus, $g_1,g_2,g_3$ and $g_4$ generate the 
$\KK$-algebra
$\Gamma(X,\mathcal{S})$. 
Setting $\deg(Z_i) := \deg(g_i)$, we obtain a 
$K$-graded epimorphism 
$$
\KK[Z_1,Z_2,Z_3^{\pm 1}] 
\ \to \ 
\Gamma(X,\mathcal{S}),
\qquad
Z_1 \mapsto g_1, \
Z_2 \mapsto g_2, \
Z_3 \mapsto g_3,
$$
which, by dimension reasons, is even an 
isomorphism.
The kernel of the projection 
$K \to \Cl(X)$ is $K^0 = 2\ZZ D_1$ 
and a character as in 
Construction~\ref{constr:crtorsion} 
is 
$$ 
\chi \colon K^0 \ \to \ \KK(X)^*,
\qquad
2nD_1 \ \mapsto \ f_1^{n}.
$$
The ideal $\mathcal{I}$ is generated by 
$1 - f_1$, where 
$f_1 \in \Gamma(X,\mathcal{S}_{-2D_1})$,
see Remark~\ref{rem:idfingen} below.
Consequently, the Cox ring of $X$ 
is given as
$$
\mathcal{R}(X) 
\ \cong \ 
\Gamma(X,\mathcal{S})
/ 
\Gamma(X,\mathcal{I})
\ \cong \ 
\KK[Z_1,Z_2,Z_3^{\pm 1}]/ \bangle{1-Z_3^{-1}}
\ \cong \
\KK[Z_1,Z_2],
$$ 
where the $\Cl(X)$-grading on the polynomial 
ring $\KK[Z_1,Z_2]$ is given by 
$\deg(Z_1) = \deg(Z_2) = [D_1]$.
\end{example}

\subsection{Well-definedness}
Here we prove the claims made in Construction~\ref{constr:crtorsion}
and Proposition~\ref{prop:coxsheafunique}.
In particular, we show that, up to isomorphy,
Cox sheaf and Cox ring do not depend on the choices 
made in their construction.
 
\begin{lemma}
\label{lem:secpres}
Situation as in Construction~\ref{constr:crtorsion}.
Consider the $\Cl(X)$-grading of the 
sheaf~$\mathcal{S}$ defined by
$$ 
\mathcal{S}
\ = \ 
\bigoplus_{[D] \in \Cl(X)} \mathcal{S}_{[D]} ,
\qquad
\qquad
\mathcal{S}_{[D]} 
\ := \ 
\bigoplus_{D' \in c^{-1}([D])} \mathcal{S}_{D'}.
$$
Given $f \in \Gamma(U,\mathcal{I})$ and 
$D \in K$, 
the $\Cl(X)$-homogeneous component 
$f_{[D]} \in \Gamma(U,\mathcal{S}_{[D]})$
of $f$ has a unique representation
$$ 
f_{[D]}
\ = \ 
\sum_{E \in K^0} (1 - \chi(E))  f_{E},
\qquad
\text{where }
f_E \in \Gamma(U,\mathcal{S}_D)
\text{ and } 
\chi(E) \in \Gamma(U,\mathcal{S}_{-E}).
$$
In particular, the sheaf $\mathcal{I}$ of ideals 
is $\Cl(X)$-homogeneous. 
Moreover, if $f \in \Gamma(U,\mathcal{I})$ 
is $K$-homogeneous, then it is the zero section.
\end{lemma}

\begin{proof}
To obtain uniqueness of the representation 
of $f_{[D]}$, 
observe that for every $0 \ne E \in K^0$, 
the product $-\chi(E)  f_{E}$ is the 
$K$-homogeneous component of degree $D-E$
of $f_{[D]}$.
We show existence.
By definition of the sheaf of ideals 
$\mathcal{I}$,
every germ $f_x \in \mathcal{I}_x$
can on a suitable neighbourhood $U_x$ 
be represented by a section 
$$
g 
\ = \ 
\sum_{E \in K^0} (1 - \chi(E)) g_{E},
\qquad
\text{where }
g_E \in \Gamma(U_x,\mathcal{S}).
$$
Collecting the $\Cl(X)$-homogeneous parts
on the right hand side represents 
the $\Cl(X)$-homogeneous part 
$h \in \Gamma(U_x,\mathcal{S}_{[D]})$ 
of degree $[D]$ of 
$g \in \Gamma(U_x,\mathcal{S})$
as follows:
$$
h 
\ = \ 
\sum_{E \in K^0} (1 - \chi(E)) h_{E},
\qquad
\text{where }
h_E \in \Gamma(U_x,\mathcal{S}_{[D]}).
$$
Note that we have $h \in \Gamma(U_x,\mathcal{I})$
and $h$ represents $f_{[D],x}$.
Now, developping each 
$h_E \in \Gamma(U_x,\mathcal{S}_{[D]})$
according to the $K$-grading gives 
representations
$$
h_E
\ = \ 
\sum_{D' \in D+K^0} h_{E,D'},
\qquad
\text{where }
h_{E,D'} \in \Gamma(U_x,\mathcal{S}_{D'}).
$$
The section
$h_{E,D'}' := \chi(D'-D)h_{E,D'}$ 
is $K$-homogeneous of degree $D$,
and we have the identity
$$ 
(1-\chi(E)) h_{E,D'}
\ = \ 
(1 - \chi(E+D-D')) h_{E,D'}'
-
(1 - \chi(D-D')) h_{E,D'}'.
$$
Plugging this into the representation of $h$
establishes the desired representation of $f_{[D]}$ 
locally. 
By uniqueness, we may glue the local representations.
\end{proof}

\begin{remark}
\label{rem:idfingen}
Situation as in Construction~\ref{constr:crtorsion}.
Then, for any two divisors 
$E, E' \in K^0$, one has the identities
\begin{eqnarray*}
1 \ - \ \chi(E + E') 
& = &
(1 \  - \ \chi(E)) \ + \ (1 \, - \, \chi(E')) \chi(E),
\\
1 \ - \ \chi( - E) 
&  = &
(1 \ - \ \chi(E))( - \chi(-E)).
\end{eqnarray*}
Together with Lemma~\ref{lem:secpres},
this implies that for any basis
$E_1, \ldots, E_s$ of $K^0$ and 
any open $U \subseteq X$,
the ideal $\Gamma(U, \mathcal{I})$ is generated
by $1 - \chi(E_i)$, where $ 1 \le i \le s$.
\end{remark}

\begin{lemma}
\label{lem:clhom2Khom}
Situation as in Construction~\ref{constr:crtorsion}.
If $f \in \Gamma(U,\mathcal{S})$ is 
$\Cl(X)$-homogeneous of degree $[D]$
for some $D \in K$,
then there is a $K$-homogeneous 
$f' \in \Gamma(U,\mathcal{S})$ of 
degree $D$ with 
$f-f' \in \Gamma(U,\mathcal{I})$.
\end{lemma}

\begin{proof}
Writing the $\Cl(X)$-homogeneous $f$ as a sum 
of $K$-homogeneous functions $f_{D'}$, 
we obtain the assertion by means of the 
following trick:
$$ 
f
\ = \ 
\sum_{D' \in D+K^0} f_{D'}
\ = \ 
\sum_{D' \in D+K^0} \chi(D' -D) f_{D'} 
+
\sum_{D' \in D+K^0} (1-\chi(D' -D)) f_{D'}.
$$
\end{proof}

\begin{lemma}
\label{lem:gradedcomps}
Situation as in Construction~\ref{constr:crtorsion}.
Then, for every $D \in K$, we have an 
isomorphism of sheaves
$\pi_{\vert \mathcal{S}_D} \colon \mathcal{S}_D \to \mathcal{R}_{[D]}$.
\end{lemma}

\begin{proof}
Lemma~\ref{lem:secpres} shows that
the homomorphism
$\pi_{\vert \mathcal{S}_D}$
is stalkwise injective and 
from Lemma~\ref{lem:clhom2Khom}
we infer that it is stalkwise 
surjective.
\end{proof}

\begin{lemma}
\label{lem:globfact}
Situation as in Construction~\ref{constr:crtorsion}.
Then, for every open subset $U \subseteq X$, 
we have a canonical isomorphism
\begin{eqnarray*}
\Gamma(U,\mathcal{S}) 
 / 
\Gamma(U,\mathcal{I})
& \cong & 
\Gamma(U,\mathcal{S}/\mathcal{I}).
\end{eqnarray*}
\end{lemma}

\begin{proof}
The canonical map
$\psi \colon \Gamma(U,\mathcal{S})/\Gamma(U,\mathcal{I}) 
\to 
\Gamma(U,\mathcal{S}/\mathcal{I})$ 
is injective.
In order to see that it is as well surjective, 
let $h \in \Gamma(U,\mathcal{S}/\mathcal{I})$
be given.
Then there are a covering of $U$ by open subsets 
$U_i$ and sections $g_i \in \Gamma(U_i,\mathcal{S})$ 
such that $h_{\vert U_i} = \psi(g_i)$ holds and 
$g_j - g_i$ belongs to $\Gamma(U_i \cap U_j,\mathcal{I})$.
Consider the $\Cl(X)$-homogeneous parts 
$g_{i,[D]} \in  \Gamma(U_i,\mathcal{S}_{[D]})$
of $g_i$.
By Lemma~\ref{lem:secpres}, the ideal sheaf 
$\mathcal{I}$ is homogeneous and thus also 
$g_{j,[D]} - g_{i,[D]}$ belongs to 
$\Gamma(U_i\cap U_j,\mathcal{I})$.
Moreover, Lemma~\ref{lem:clhom2Khom} provides
$K$-homogeneous $f_{i,D}$ with 
$f_{i,D} - g_{i,[D]}$ in $\Gamma(U_i,\mathcal{I})$.
The differences $f_{j,D}-f_{i,D}$ lie in 
$\Gamma(U_i \cap U_j,\mathcal{I})$ and hence,
by Lemma~\ref{lem:secpres}, vanish.
Thus, the $f_{i,D}$ fit together to 
$K$-homogeneous sections
$f_D \in \Gamma(U,\mathcal{S})$.
By construction, $f = \sum f_D$ satisfies 
$\psi(f) = h$. 
\end{proof}


\begin{proof}[Proof of Proposition~\ref{prop:coxsheafunique}]
In a first step, we reduce to Cox sheaves arising 
from finitely generated subgroups of $\WDiv(X)$.
So, let $K \subseteq \WDiv(X)$ and 
$\chi \colon K^0 \to \KK(X)^*$ be any data
as in~\ref{constr:crtorsion}.
Choose a finitely generated subgroup 
$K_1 \subseteq K$ projecting onto $\Cl(X)$.
Restricting $\chi$ gives a character 
$\chi_1 \colon K_1^0 \to \KK(X)^*$.
The inclusion $K_1 \to K$ defines an
injection $\mathcal{S}_1 \to \mathcal{S}$
sending the ideal $\mathcal{I}_1$ defined 
by $\chi_1$ to the ideal $\mathcal{I}$ 
defined by $\chi$.
This gives a $\Cl(X)$-graded injection 
$\mathcal{R}_1 \to \mathcal{R}$ of the 
Cox sheaves associated to $K_1,\chi_1$ 
and $K,\chi$ respectively.
Lemma~\ref{lem:clhom2Khom} shows that 
every $\Cl(X)$-homogeneous section of 
$\mathcal{R}$ can be represented by a 
$K_1$-homogeneous section of $\mathcal{S}$.
Thus, $\mathcal{R}_1 \to \mathcal{R}$ is
also surjective.

Next
we show that for a fixed finitely generated 
$K \subseteq \WDiv(X)$,
any two characters $\chi,\chi' \colon K^0 \to \KK(X)^*$ 
as in~\ref{constr:crtorsion}
give rise to isomorphic Cox sheaves 
$\mathcal{R}'$ and $\mathcal{R}$.
For this note that the product $\chi^{-1} \chi'$ 
sends $K^0$ to $\Gamma(X,\mathcal{O}^*)$.
Using $\Gamma(X,\mathcal{O}^*) = \KK^*$,
we may extend $\chi^{-1} \chi'$ to a homomorphism 
$\vartheta \colon K \to \Gamma(X,\mathcal{O}^*)$
and obtain a graded automorphism $(\alpha,\id)$ 
of $\mathcal{S}$ by 
$$ 
\alpha_D \colon \mathcal{S}_D \ \to \ \mathcal{S}_D,
\qquad
f \ \mapsto \ \vartheta(D) f.
$$
By construction, this automorphism sends 
the ideal $\mathcal{I}'$ to the ideal 
$\mathcal{I}$ and induces a graded
isomorphism from $\mathcal{S}/\mathcal{I}'$
onto $\mathcal{S}/\mathcal{I}$.

Now consider two finitely generated subgroups
$K,K' \subseteq \WDiv(X)$ both projecting 
onto $\Cl(X)$.
Then we find a homomorphism $\t{\alpha} \colon K \to K'$ 
such that the following diagram is commutative
$$ 
\xymatrix{
K 
\ar[rr]^{\t{\alpha}}
\ar[dr]
&& 
K'
\ar[dl]
\\
&
{\Cl(X)}
&
}
$$
This homomorphism $\t{\alpha} \colon K \to K'$ 
must be of the form $\t{\alpha}(D) = D - \div(\eta(D))$
with a homomorphism $\eta \colon K \to \KK(X)^*$.
Choose a character 
$\chi' \colon K'^0 \to \KK(X)^*$ 
as in~\ref{constr:crtorsion}.
Then, for $D \in K^0$, we have 
$$ 
D - \div(\eta(D))
\ = \ 
\t{\alpha}(D)
\ = \ 
\div(\chi'(\t{\alpha}(D))).
$$
Thus, $D$ equals the divisor of the function 
$\chi(D) := \chi'(\t{\alpha}(D)) \eta(D)$.
This defines a character $\chi \colon K^0 \to \KK(X)^*$.
Altogether, we obtain a morphism $(\alpha,\t{\alpha})$ 
of the sheaves of divisorial algebras 
$\mathcal{S}$ and $\mathcal{S}'$ 
associated to $K$ and $K'$ by 
$$ 
\alpha_D \colon \mathcal{S}_D \ \to \ \mathcal{S}_{\t{\alpha}(D)}',
\qquad
f \ \mapsto \ \eta(D) f.
$$
By construction, it sends the ideal $\mathcal{I}$ 
defined by $\chi$ to the ideal $\mathcal{I}'$ defined 
by  $\chi'$.
Using Lemma~\ref{lem:gradedcomps},
we see that the induced homomorphism 
$\mathcal{R} \to \mathcal{R}'$ is an isomorphism 
on the homogeneous components and thus it is an 
isomorphism.
\end{proof}

\subsection{Examples}
For a normal prevariety $X$ with 
a free finitely generated divisor 
class group, 
we obtained in Proposition~\ref{prop:cralgfree}
that the Cox ring is a 
unique factorization domain 
having $\Gamma(X,\mathcal{O}^*)$ 
as its units.
Here we provide two examples showing 
that these statements need not hold any 
more if there is torsion in the divisor 
class group.
As usual, $\KK$ is an algebraically 
closed field of characteristic zero.

\begin{example}
[An affine surface with non-factorial Cox ring]
\label{exam:sl2n}
Consider the smooth affine surface 
\begin{eqnarray*}
Z
& := & 
V(\KK^3; T_1^2 - T_2T_3 - 1).
\end{eqnarray*}
We claim that $\Gamma(Z,\mathcal{O}^*) = \KK^*$
and $\Cl(Z) \cong \ZZ$ hold.
To see this, consider $f_i := T_{i \vert Z}$ 
and the prime divisors
$$
D_+ 
\ := \ 
V(Z;f_1-1, f_2)
\ = \ 
\{1\} \times \{0\} \times \KK,
$$
$$
D_- 
\ := \ 
V(Z;f_1+1, f_2)
\ = \ 
\{-1\} \times \{0\} \times \KK.
$$
Then we have $\div(f_2) = D_+ + D_-$.
In particular, $D_+$ 
is linearly equivalent to $-D_-$. 
Moreover, we have 
$$
Z \setminus \Supp(\div(f_2)) 
\ = \ 
Z_{f_2}
\ \cong \
\KK^* \times \KK.
$$
This gives $\Gamma(Z, \Of)^* = \KK^*$,
and shows that $\Cl(Z)$ is 
generated by the class~$[D_+]$.
Now suppose that $n[D_+] = 0$ 
holds for some $n > 0$.
Then we have $nD_+ = \div(f)$ 
with $f \in \Gamma(Z, \Of)$
and $f_2^n = fh$ holds with some 
$h \in \Gamma(Z, \Of)$ satisfying 
$\div(h) = nD_-$.
Look at the $\ZZ$-grading 
of $\Gamma(Z,\mathcal{O})$ given by
$$
\deg(f_1) = 0, \qquad \deg(f_2) = 1, \qquad \deg(f_3) = -1.
$$
Any element of positive degree 
is a multiple of $f_2$. 
It follows that in the decomposition 
$f_2^n = fh$ one of the factors 
$f$ or $h$ must be a multiple of $f_2$, 
a contradiction. 
This shows that $\Cl(Z)$ is freely 
generated by $[D_+]$.

Now consider the involution $Z \to Z$
sending $z$ to $-z$ and let 
$\pi \colon Z \to X$ denote the 
quotient of the corresponding free
$\ZZ/2\ZZ$-action. 
We claim that $\Cl(X)$ is isomorphic 
to $\ZZ/2\ZZ$ and is generated by the 
class of $D := \pi(D_+)$. 
Indeed, the subset  
$$
X \setminus \Supp(D) 
\ = \ 
\pi(Z_{f_2})
\ \cong \ 
\KK^* \times \KK
$$
is factorial and $2D$ equals $\div(f_2^2)$. 
Moreover, the divisor $D$ is not principal, because 
$\pi^*(D) = D_+ + D_-$ is not the divisor 
of a $\ZZ/2\ZZ$-invariant function on $Z$. 
This verifies our claim. 
Moreover, we have 
$\Gamma(X,\mathcal{O}^*) = \KK^*$.

In order to determine the Cox ring of $X$,
take $K = \ZZ D \subseteq \WDiv(X)$,
and let~$\mathcal{S}$ denote the associated 
sheaf of divisorial algebras.
Then, as $\Gamma(X, \mathcal{S}_0)$-modules,
$\Gamma(X, \mathcal{S}_D)$ and $\Gamma(X, \mathcal{S}_{-D})$ 
are generated by the sections
$$
a_1 := 1, 
\
a_2 := f_1f_2^{-1}, 
\ 
a_3 := f_2^{-1}f_3
\
\in \Gamma(X,\mathcal{S}_D),
$$
$$
b_1 := f_1f_2, 
\
b_2 := f_2^2, 
\ 
b_3 := f_2f_3
\
\in \Gamma(X,\mathcal{S}_{-D}).
$$
Thus, using the fact that $f_2^{\pm 2}$ 
define invertible 
elements of degree $\mp 2D$, 
we see that $a_1,a_2,a_3,b_1,b_2,b_3$ 
generate the algebra $\Gamma(X, \mathcal{S})$.
Now, take the character $\chi \colon K^0 \to \KK(X)^*$
sending $2nD$ to $f_2^{2n}$. 
Then, by Remark~\ref{rem:idfingen},
the associated ideal $\Gamma(X, \mathcal{I})$ 
is generated by $1-f_2^2$. 
The generators of the factor algebra 
$\Gamma(X, \mathcal{S}) \, / \, \Gamma(X, \mathcal{I})$ are
$$
Z_1 = a_2 + \mathcal{I} = b_1 + \mathcal{I}, \qquad  
Z_2 = a_1 + \mathcal{I} = b_2 + \mathcal{I}, \qquad
Z_3 = a_3 + \mathcal{I} = b_3 + \mathcal{I}. 
$$
The defining relation is $Z_1^2-Z_2Z_3=1$.
Thus the Cox ring $\mathcal{R}(X)$
is isomorphic to $\Gamma(Z,\Of)$.
In particular, it is not a factorial ring.
\end{example}

\begin{example}
[A surface with only constant invertible functions
but non-constant invertible elements in the Cox ring] 
\label{ex:nonconstunits}
Consider the affine surface 
\begin{eqnarray*}
X 
& := & 
V( \KK^3; \; T_1T_2T_3 - T_1^2 - T_2^2 - T_3^2 + 4).
\end{eqnarray*}
This is the quotient space of the torus 
$\TT^2 := (\KK^*)^2$ with respect to
the $\ZZ/2\ZZ$-action defined by the involution 
$t \mapsto t^{-1}$;  
the quotient map is explicitly 
given as
$$ 
\pi \colon \TT^2 \  \to \ X,
\qquad
t 
\ \mapsto \ 
(t_1 + t_1^{-1},  t_2 + t_2^{-1},  t_1t_2 + t_1^{-1}t_2^{-1}).
$$
Since every $\ZZ/2\ZZ$-invariant invertible function
on $\TT^2$ is constant, we have $\Gamma(X, \Of^*) = \KK^*$. 
Moreover, using \cite[Proposition~5.1]{KKV}, 
one verifies  
$$
\Cl(X) \ \cong \ \ZZ/2\ZZ \oplus \ZZ/2\ZZ \oplus \ZZ/2\ZZ,
\qquad
\Pic(X) \ = \ 0.
$$ 

Let us see that the Cox ring $\mathcal{R}(X)$ 
has non-constant invertible elements.
Set $f_i := T_{i \vert X}$ and consider the divisors
$$
D_{\pm}
\ := \
V(X ; f_1 \pm 2, \ f_2 \pm f_3),
\qquad
D \ := \ D_+ + D_-.
$$
Then, using the relations
$(f_1 \pm 2)(f_2f_3 - f_1 \pm 2)  =  (f_2 \pm f_3)^2$,
one verifies $\div(f_1 \pm 2) = 2D_{\pm}$.
Consequently, we obtain
\begin{eqnarray*}
2D
& = & 
\div(f_1^2-4).
\end{eqnarray*}
Moreover, $D$ is not principal, because 
otherwise $f_1^2-4$ 
must be a square and hence also 
$\pi^*(f_1^2-4)$ is a square,
which is impossible due to
$$
\pi^*(f_1^2-4) 
\ = \ 
t_1^2 + t_1^{-2}-4
\ = \ 
(t_1 + t_1^{-1}+2)
(t_1 + t_1^{-1}-2).
$$

Now choose Weil divisors $D_i$ on $X$
such that $D,D_2,D_3$ form a basis
for a group $K \subseteq \WDiv(X)$ projecting onto 
$\Cl(X)$,
and let $\mathcal{S}$ be the associated sheaf of divisorial 
algebras.
As usual, let $K^0 \subseteq K$ be the subgroup 
consisting of principal divisors and 
fix a character $\chi \colon K^0 \to \KK(X)^*$ with 
$\chi(2D)  = f_1^2-4$.
By Remark~\ref{rem:idfingen},
the associated ideal $\Gamma(X,\mathcal{I})$
in $\Gamma(X,\mathcal{S})$
is generated by 
$$ 
1 - \chi(2D), \quad 1 - \chi(2D_2), \quad 1 - \chi(2D_3),
$$
where $\chi(2D) = f_1^2-4$ lives in $\Gamma(X, \mathcal{S}_{-2D})$.
Now consider $f_1 \in \Gamma(X,\mathcal{S}_0)$ 
and the canonical section $1_D \in \Gamma(X,\mathcal{S}_D)$.
Then we have 
$$ 
(f_1+1_D)(f_1-1_D)
\ = \ 
f_1^2 - 1_D^2
\ = \ 
4 - 1_D^2 \cdot (1-\chi(2D))
\ \in \ 
\KK^* + \Gamma(X,\mathcal{I}).
$$
Consequently, the section $f_1+1_D  \in \Gamma(X,\mathcal{S})$
defines a unit in $\Gamma(X,\mathcal{R})$.
Note that $f_1+1_D$ is not $\Cl(X)$-homogeneous.
\end{example}


\section{Algebraic properties of the Cox ring}
\label{sec:cralgprops}

\subsection{Integrity and Normality}
As before, we work over an algebraically 
closed field $\KK$ of characteristic zero.
The following statement ensures in particular 
that the Cox ring is always  
a normal integral ring.

\begin{theorem}
\label{prop:crnormal}
Let $X$ be a normal prevariety with 
only constant invertible functions,
finitely generated divisor class 
group, and Cox sheaf $\mathcal{R}$.
Then, for every open $U \subseteq X$,
the ring $\Gamma(U,\mathcal{R})$ 
is integral and normal.
\end{theorem}

The proof is based on the geometric 
construction~\ref{constr:univtorssmooth}
which is also used later and therefore 
occurs separately.
We begin with two preparing observations.

\begin{lemma}
\label{rem:codim2divsheaf2}
Situation as in Construction~\ref{constr:crtorsion}.
For any two
open subsets $V \subseteq U \subseteq X$ 
such that $U \setminus V$ is of 
codimension at least two in $U$, one
has the restriction isomorphism 
$$
\Gamma(U, \mathcal{R}) 
\ \to \
\Gamma(V, \mathcal{R}).
$$
In particular, the algebra $\Gamma(U, \mathcal{R})$
equals the algebra $\Gamma(U_\reg, \mathcal{R})$,
where $U_\reg \subseteq U$ denotes the set
of smooth points.
\end{lemma}

\begin{proof}
According to Remark~\ref{rem:codim2divsheaf},
the restriction  
$\Gamma(U, \mathcal{S}) \to \Gamma(V, \mathcal{S})$
is an isomorphism.
Lemma~\ref{lem:secpres} ensures that 
$\Gamma(U, \mathcal{I})$ is mapped isomorphically 
onto $\Gamma(V, \mathcal{I})$ under this
isomorphism.
By Lemma~\ref{lem:globfact}, we have 
$\Gamma(U, \mathcal{R}) = 
\Gamma(U, \mathcal{S})/\Gamma(U, \mathcal{I})$
and 
$\Gamma(V, \mathcal{R}) = 
\Gamma(V, \mathcal{S})/\Gamma(V, \mathcal{I})$,
which gives the assertion.
\end{proof}

\begin{lemma}
\label{lem:Iradical}
Situation as in Construction~\ref{constr:crtorsion}.
Then for every open $U \subseteq X$,
the ideal $\Gamma(U,\mathcal{I}) \subseteq \Gamma(U,\mathcal{S})$ 
is radical.
\end{lemma}

\begin{proof}
By Lemma~\ref{lem:globfact}, 
the ideal $\Gamma(U,\mathcal{I})$
is radical if and only if the algebra
$\Gamma(U,\mathcal{R})$ has no nilpotent
elements.
Proposition~\ref{prop:coxsheafunique}
thus allows us to assume that $\mathcal{S}$ 
arises from a finitely generated group $K$.
Moreover, by Remark~\ref{rem:codim2divsheaf},
we may assume that $X$ is smooth
and it suffices to verify the assertion for 
affine $U \subseteq X$.
We consider $\t{U} = \Spec \, \Gamma(U,\mathcal{S})$
and the zero set $\rq{U} \subseteq \t{U}$ of 
$\Gamma(U,\mathcal{I})$.
Note that $\rq{U}$ is invariant under the 
action of the quasitorus 
$H_X = \Spec \, \KK[\Cl(X)]$ 
on $\t{U}$ given by the $\Cl(X)$-grading.

Now, let $f \in \Gamma(U,\mathcal{S})$ with 
$f^n \in \Gamma(U,\mathcal{I})$ for some $n > 0$.
Then $f$ and thus also every 
$\Cl(X)$-homogeneous 
component $f_{[D]}$ of $f$ 
vanishes along $\rq{U}$.
Consequently,
$f_{[D]}^{m} \in  \Gamma(U,\mathcal{I})$
holds for some $m > 0$.
By Lemma~\ref{lem:clhom2Khom},
we may write $f_{[D]} = f_D + g$ 
with $f_D \in \Gamma(U,\mathcal{S}_D)$
and $g \in \Gamma(U,\mathcal{I})$. 
We obtain $f_D^{m} \in  \Gamma(U,\mathcal{I})$. 
By Lemma~\ref{lem:secpres},
this implies  $f_D^{m}= 0$ and 
thus $f_D = 0$, 
which in turn gives  
$f_{[D]} \in  \Gamma(U,\mathcal{I})$
and hence 
$f \in  \Gamma(U,\mathcal{I})$.
\end{proof}

\begin{construction}
\label{constr:univtorssmooth}
Situation as in Construction~\ref{constr:crtorsion}.
Assume that $K \subseteq \WDiv(X)$ is 
finitely generated and $X$ is smooth.
Consider $\t{X} := \Spec_X \, \mathcal{S}$
with the action of the torus 
$H := \Spec \, \KK[K]$ 
and the geometric 
quotient $p \colon \t{X} \to X$
as in Construction~\ref{constr:divrelspec}.
Then,
with $\rq{X} := V(\mathcal{I})$ and 
$H_X := \Spec \; \KK[\Cl(X)]$,
we have a commutative diagram
$$ 
\xymatrix{
{\rq{X}}
\ar[rr]
\ar[d]^{/H_X}_{q_X}
&&
{\t{X}}
\ar[d]^{/H}_{p}
\\
X
\ar@{=}[rr]
&&
X
}
$$
The prevariety $\rq{X}$ is smooth,
and, if $X$ is of affine intersection,
then it is quasiaffine.
The quasitorus $H_X \subseteq H$ acts freely 
on $\rq{X}$ and $q_X \colon \rq{X} \to X$
is a geometric quotient for this action;
in particular, it is an \'etale $H_X$-principal 
bundle.
Moreover, we have a canonical isomorphism 
of sheaves
\begin{eqnarray*}
\mathcal{R}
& \cong & 
(q_X)_*(\mathcal{O}_{\rq{X}}).
\end{eqnarray*}
\end{construction}

\begin{proof}
With the restriction $q_X \colon \rq{X} \to X$ 
of $p \colon \t{X} \to X$ we obviously obtain 
a commutative diagram as above.
Moreover, Lemma~\ref{lem:Iradical}
gives us $\mathcal{R} \cong q_*(\mathcal{O}_{\rq{X}})$.
Since the ideal $\mathcal{I}$ is $\Cl(X)$-homogeneous,
the quasitorus $H_X \subseteq H$ leaves 
$\rq{X}$ invariant.
Moreover, we see that $q_X \colon \rq{X} \to X$ 
is a good quotient for this action,
because we have the canonical isomorphisms
$$ 
(q_X)_*(\mathcal{O}_{\rq{X}})_0
\ \cong \ 
\mathcal{R}_0
\ \cong \ 
\mathcal{O}_X
\ \cong \ 
\mathcal{S}_0
\ \cong \ 
p_*(\mathcal{O}_{\t{X}})_0.
$$ 
Freeness of the $H_X$-action on $\rq{X}$ is due to
the fact that $H_X$ acts as a subgroup of
the freely acting~$H$, 
see Remark~\ref{rem:principalbundle}.
As a consequence, we see
that $q_X \colon \rq{X} \to X$ 
is a geometric quotient. 
Luna's Slice Theorem \cite{Lu} gives
commutative diagrams
$$ 
\xymatrix{
H_X \times S
\ar[r]
\ar[d]_{\pr_S}
&
q_X^{-1}(U)
\ar@{}[r]|{\subseteq}
\ar[d]_{q_X}
&
{\rq{X}}
\ar[d]^{q_X}
\\
S
\ar[r]
&
U
\ar@{}[r]|{\subseteq}
&
X
}
$$
where $U \subseteq X$ are open sets covering $X$ 
and the horizontal arrows are \'etale morphisms.
By~\cite[Proposition~I.3.17]{Mi}, \'etale morphisms 
preserve smoothness and thus $\rq{X}$ 
inherits smoothness from $X$.
If $X$ is of affine intersection,
then $\t{X}$ is quasiaffine, 
see~Corollary~\ref{cor:quaff2},
and thus $\rq{X}$ is quasiaffine.
\end{proof}

\begin{lemma}
\label{lem:polirred}
Let $\LL$ be a field of characteristic zero
containing all roots of unity,
and assume that $a \in \LL$ 
is not a proper power.
Then, for any $n \in \ZZ_{\ge 1}$, 
the polynomial $1-at^n$ 
is irreducible in $\LL[t,t^{-1}]$.
\end{lemma}

\begin{proof}
Over the algebraic closure of $\LL$
we have $1 - at^n = (1 -  a_1t) \cdots ( 1 - a_nt)$, 
where $a_i^n = a$ and any two $a_i$ differ by 
a $n$-th root of unity.
If $1-at^n$ would split over $\LL$
non-trivially into
$h_1(t) h_2(t)$, 
then $a_1^k$ must be
contained in $\LL$ for some $k<n$.
But then also $a_1^d$ lies in $\LL$
for the greatest
common divisor $d$ of $n$ and $k$.
Thus $a$ is a proper power, a contradiction.
\end{proof}

\begin{proof}[Proof of Theorem~\ref{prop:crnormal}]
According to Proposition~\ref{prop:coxsheafunique}
and Lemma~\ref{rem:codim2divsheaf2},
we may assume that we are in the setting 
of Construction~\ref{constr:univtorssmooth},
where it suffices to prove that $\rq{X}$
is irreducible.
Since $q_X \colon \rq{X} \to X$ is surjective,
some irreducible component 
$\rq{X}_1 \subseteq \rq{X}$ 
dominates $X$.
We verify that $\rq{X}_1$ equals $\rq{X}$
by checking that $q_X^{-1}(U)$ 
is irreducible
for suitable open neighbourhoods $U \subseteq X$
covering $X$.

Let $D_1, \ldots, D_s$ be a basis of $K$ 
such that $n_1D_1, \ldots, n_kD_k$, 
where $1 \le k \le s$, is a basis of $K^0$.
Enlarging $K$, if necessary, we may assume 
that the $D_i$ are primitive, i.e., 
no proper multiples.
We take subsets $U \subseteq X$ such that 
on $U$ every $D_i$ is principal, 
say $D_i = \div(f_i)$.
Then, with $\deg(T_i) := D_i$,
Remark~\ref{rem:locprinc2split} provides
a $K$-graded isomorphism 
$$ 
\Gamma(U,\mathcal{O}) \otimes_{\KK} \KK[T_1^{\pm 1}, \ldots, T_s^{\pm 1}]
\ \to \ 
\Gamma(U,\mathcal{S}),
\qquad
g \otimes T_1^{\nu_1} \cdots T_s^{\nu_s} 
\ \mapsto \ 
g  f_1^{-\nu_1} \cdots f_s^{-\nu_s}.
$$
In particular, this identifies $p^{-1}(U)$ with  
$U \times \TT^s$, where $\TT^s := (\KK^*)^s$.
According to Remark~\ref{rem:idfingen}, the ideal 
$\Gamma(U, \mathcal{I})$
is generated by $1 - \chi(n_iD_i)$, 
where $1 \le i \le k$. 
Thus $q_X^{-1}(U)$ is given in $U \times \TT^s$ by 
the equations
$$
1-\chi(n_iD_i)f_i^{n_i}T_i^{n_i}  \ = \ 0,
\qquad
1 \le i \le k.
$$ 
To obtain irreducibility of  $q_X^{-1}(U)$,
it suffices to show that each 
$1-\chi(n_iD_i)f_i^{n_i}T_i^{n_i}$ is irreducible 
in $\KK(X)[T_i^{\pm 1}]$.
With respect to the variable $S_i := f_iT_i$, 
this means to verify irreducibility of 
\begin{eqnarray*}
1 \ - \ \chi(n_iD_i)S_i^{n_i}
& \in & 
\KK(X)[S_i^{\pm 1}].
\end{eqnarray*}
In view of Lemma~\ref{lem:polirred},
we have to show that $\chi(n_iD_i)$ is not a 
proper power in $\KK(X)$. 
Assume the contrary. 
Then we obtain $n_iD_i = k_i \div(h_i)$ with some 
$h_i \in \KK(X)$.
Since $D_i$ is primitive, $k_i$ divides $n_i$
and thus, $n_i/k_i D_i$ is principal. 
A contradiction to the choice of $n_i$.

The fact that each ring $\Gamma(U,\mathcal{R})$
is normal follows directly from the fact that
it is the ring of functions of an open subset
of the smooth prevariety $\rq{X}$.
\end{proof}

\subsection{Localization and units}
We treat localization by homogeneous elements
and consider the units of the Cox ring 
$\mathcal{R}(X)$ of a normal prevariety $X$ 
defined over an algebraically closed field $\KK$
of characteristic zero.
The main tool is the divisor
of a homogeneous element of $\mathcal{R}(X)$,
which we first define precisely.

\index{$[D]$-divisor!}%
In the setting of Construction~\ref{constr:crtorsion},
consider a divisor $D \in K$ and a non-zero element
$f \in \mathcal{R}_{[D]}(X)$.
According to Lemma~\ref{lem:clhom2Khom}, 
there is a (unique) element
$\t{f} \in \Gamma(X, \mathcal{S}_D)$ 
with $\pi(\t{f}) = f$, where 
$\pi \colon \mathcal{S} \to \mathcal{R}$ 
denotes the projection.
We define the {\em $[D]$-divisor\/} of  $f$ 
to be the effective Weil divisor
$$
\div_{[D]}(f) 
\ := \
\div_D(\t{f})
\ = \ 
\div(\t{f}) + D
\ \in \ 
\WDiv(X).
$$ 

\begin{lemma}
\label{lem:dddivprops}
The $[D]$-divisor depends neither  
on representative $D \in K$
nor on the choices made 
in~\ref{constr:crtorsion}.
Moreover, the following holds.
\begin{enumerate}
\item
For every effective $E \in \WDiv(X)$ 
there are $[D] \in \Cl(X)$ and 
$f \in \mathcal{R}_{[D]}(X)$ with
$E = \div_{[D]}(f)$. 
\item 
Let $[D] \in \Cl(X)$ and $0 \ne f \in \mathcal{R}_{[D]}(X)$.
Then $\div_{[D]}(f) = 0$
implies  
$[D] = 0$ in $\Cl(X)$.
\item
For any two non-zero homogeneous elements
$f \in \mathcal{R}_{[D_1]}(X)$
and $g \in \mathcal{R}_{[D_2]}(X)$, we have
\begin{eqnarray*}
\div_{[D_1] + [D_2]}(fg) 
& = & 
\div_{[D_1]}(f) \ + \ \div_{[D_2]}(g).
\end{eqnarray*}
\end{enumerate}
\end{lemma}

\begin{proof}
Let $f \in \mathcal{R}_{[D]}(X)$, 
consider any two isomorphisms
$\varphi_i \colon \mathcal{O}_X(D_i) \to \mathcal{R}_{[D]}$
and let $\t{f}_i$ be the sections with 
$\varphi_i(\t{f}_i) = f$.
Then $\varphi_2^{-1} \circ \varphi_1$ is 
multiplication
with some $h \in \KK(X)^*$ satisfying
$\div(h) = D_1-D_2$.
Well-definedness of the $[D]$-divisor
thus follows from 
$$ 
\div_{D_1}(\t{f}_1)
\ = \ 
\div(h\t{f}_1) + D_2
\ = \ 
\div_{D_2}(\t{f}_2).
$$
If  $\div_{[D]}(f) = 0$ holds as in~(ii), 
then, for a representative 
$\t{f} \in \Gamma(X,\mathcal{O}_X(D))$ 
of $f \in \mathcal{R}_{[D]}(X)$, 
we have $\div_D(\t{f}) = 0$ and 
hence $D$ is principal.
Observations~(i) and ~(iii) are obvious.
\end{proof}

\index{$[D]$-localization!}%
For every non-zero homogeneous element 
$f \in \mathcal{R}_{[D]}(X)$, 
we define the {\em $[D]$-localization\/} 
of $X$ by $f$ to be the  open subset 
$$ 
X_{[D],f}
\ := \ 
X \setminus \Supp(\div_{[D]}(f))
\ \subseteq \ 
X.
$$

\begin{proposition}
\label{prop:crloc}
Let $X$ be a normal prevariety with 
only constant invertible functions,
finitely generated divisor class group
and Cox ring $\mathcal{R}(X)$.
Then, for every non-zero homogeneous 
$f \in \mathcal{R}_{[D]}(X)$, 
we have a canonical isomorphism
\begin{eqnarray*}
\Gamma(X_{[D],f},\mathcal{R})
& \cong & 
\Gamma(X,\mathcal{R})_f.
\end{eqnarray*}
\end{proposition}

\begin{proof}
Let the divisor $D \in K$ represent $[D] \in \Cl(X)$ 
and consider the section
$\t{f} \in \Gamma(X, \mathcal{S}_D)$ 
with $\pi(\t{f}) = f$.
According to Remark~\ref{rem:Sloc}, 
we have 
\begin{eqnarray*}
\Gamma(X_{D,\t{f}},\mathcal{S})
& \cong & 
\Gamma(X,\mathcal{S})_{\t{f}}.
\end{eqnarray*}
The assertion thus follows from Lemma~\ref{lem:globfact}
and the fact that localization is compatible with 
passing to the factor ring.
\end{proof}

We turn to the units of the Cox ring 
$\mathcal{R}(X)$;
the following result says in particular, 
that for a complete normal variety $X$ 
they are all constant.

\begin{proposition}
\label{prop:inv2const2}
Let $X$ be a normal prevariety with 
only constant invertible functions,
finitely generated divisor class 
group and Cox ring $\mathcal{R}(X)$.
\begin{enumerate}
\item
Every homogeneous invertible element of 
$\mathcal{R}(X)$ is constant.
\item
If $\Gamma(X,\mathcal{O}) = \KK$ holds, then 
every invertible element of $\mathcal{R}(X)$
is constant.
\end{enumerate}
\end{proposition}

\begin{proof}
For~(i), let $f \in \mathcal{R}(X)^*$ be homogeneous of
degree $[D]$.
Then its inverse $g \in \mathcal{R}(X)^*$ is 
homogeneous of degree $-[D]$, and
$fg = 1$ lies in $\mathcal{R}(X)^*_0 = \KK^*$.
By Lemma~\ref{lem:dddivprops}~(iii), we have 
$$
0 
\ = \ 
\div_{0}(fg)
\ = \ 
\div_{[D]}(f) + \div_{[-D]}(g).
$$
Since the divisors $\div_{[D]}(f)$ 
and $\div_{[-D]}(g)$ 
are effective, they both vanish.
Thus, Lemma~\ref{lem:dddivprops}~(ii)
yields $[D]=0$.
This implies 
$f \in \mathcal{R}(X)^*_0 = \KK^*$ as wanted.

For~(ii), we have to show that any invertible 
$f \in \mathcal{R}(X)$ is of degree zero.
Choose a decomposition $\Cl(X) = K_0 \oplus K_t$ 
into a free part and the torsion part, and 
consider the coarsified grading 
$$ 
\mathcal{R}(X)
\ = \ 
\bigoplus_{w \in K_0} R_w,
\qquad\qquad
R_w
\ := \ 
\bigoplus_{u \in K_t} \ \mathcal{R}(X)_{w+u}.
$$
Then, as any invertible element of the $K_0$-graded
integral ring $\mathcal{R}(X)$, also $f$ 
is necessarily $K_0$-homogeneous of some degree 
$w \in K_0$.
Decomposing $f$ and $f^{-1}$ into $\Cl(X)$-homogeneous 
parts we get representations
$$ 
f 
\ = \ 
\sum_{u \in K_t} f_{w+u},
\qquad
\qquad
f^{-1} 
\ = \ 
\sum_{u \in K_t} f^{-1}_{-w+u}.
$$
Because of $ff^{-1} = 1$, we have $f_{w+v}f^{-1}_{-w-v} \ne 0$
for at least one $v \in K_t$.
Since $\Gamma(X,\mathcal{O}) = \KK$ holds,
$f_{w+v}f^{-1}_{-w-v}$ must be a non-zero constant.
Using Lemma~\ref{lem:dddivprops} we conclude
$w+v = 0$ as before.
In particular, $w = 0$ holds and thus
each $f_{w+u}$ has a torsion degree.
For a suitable power $f_{w+u}^n$ we have
$n\div_{w+u}(f_{w+u}) = 0$, 
which implies $f_{w+u} = 0$ 
for any $u \ne 0$.
\end{proof}

\begin{remark}
The affine surface $X$ treated in Example~\ref{ex:nonconstunits}
shows that requiring $\Gamma(X,\mathcal{O}^*)=\KK^*$ is 
in general not enough in order to ensure that all 
units of the Cox ring are constant.
\end{remark}

\subsection{Divisibility properties}
For normal prevarieties $X$ with a free finitely 
generated divisor class group, we saw that the
Cox ring admits unique factorization.
If we have torsion in the divisor class group this
does not need to hold any more.
However, restricting to homogeneous elements 
leads to a framework for a reasonable divisibility
theory; the precise notions are the following.

\begin{definition}
\label{def:factgrad}
\index{$K$-prime!}%
\index{$K$-height!}%
\index{factorially graded!}%
Consider an abelian group $K$ and 
a $K$-graded integral $\KK$-algebra
$R = \bigoplus_{w \in K} R_w$.
\begin{enumerate}
\item
A non-zero non-unit $f \in R$ is {\em $K$-prime\/} if
it is homogeneous and $f | gh$ with homogeneous 
$g,h \in R$ implies $f | g$ or $f | h$.
\item
We say that $R$ is {\em factorially graded\/} 
if every homogeneous  non-zero non-unit $f \in R$
is a product of $K$-primes.
\item
An ideal $\mathfrak{a} \lhd R$ is {\em $K$-prime\/} if
it is homogeneous and for any two homogeneous $f,g \in R$  
with $fg \in \mathfrak{a}$ one has either $f \in \mathfrak{a}$
or $g \in \mathfrak{a}$. 
\item
A $K$-prime ideal $\mathfrak{a} \lhd R$ has {\em $K$-height\/}
$d$ if $d$ is maximal admitting a chain 
$\mathfrak{a}_0 
\subset \mathfrak{a}_1 \subset \ldots \subset 
\mathfrak{a}_d=\mathfrak{a}$
of $K$-prime ideals. 
\end{enumerate}
\end{definition}

Let us look at these concepts also from the 
geometric point of view.
Consider a prevariety $Y$ with an action of an 
algebraic group $H$.
Then $H$ acts also on the group $\WDiv(Y)$ of Weil 
divisors via 
\begin{eqnarray*}
h \cdot \sum a_D D
& := & 
\sum a_D (h \mal D).
\end{eqnarray*}
\index{$H$-prime divisor!}%
\index{divisor!$H$-prime}%
\index{$H$-factorial!}%
By an {\em $H$-prime divisor} we mean a non-zero 
sum $\sum a_D D$ with prime divisors $D$ 
such that $a_D \in \{0,1\}$ always holds
and the $D$ with $a_D=1$ are transitively permuted by~$H$. 
Note that every $H$-invariant divisor 
is a unique sum of $H$-prime divisors. 
We say that $Y$ is {\em $H$-factorial}
if every $H$-invariant Weil divisor on $Y$ 
is principal.

\begin{proposition}
\label{hfactorial}
Let $H = \Spec \, \KK[K]$ be a quasitorus
and $W$ an irreducible normal 
quasiaffine $H$-variety.
Consider the $K$-graded algebra 
$R := \Gamma(W, \mathcal{O})$
and assume $R^* = \KK^*$.
Then the following statements are equivalent.
\begin{enumerate}
\item 
Every $K$-prime ideal of $K$-height one in 
$R$ is principal. 
\item
The variety $W$ is $H$-factorial.
\item
The ring $R$ is factorially graded.
\end{enumerate}
Moreover, if one of these statements holds, 
then a homogeneous non-zero non-unit $f \in R$
is $K$-prime if and only if the divisor $\div(f)$ is 
$H$-prime, and every $H$-prime divisor
is of the form $\div(f)$ with a $K$-prime $f \in R$.
\end{proposition}

%
%
%

\begin{proof}
Assume that~(i) holds and let $D$ be an 
$H$-invariant Weil divisor on $W$. 
Write $D=D_1 + \ldots + D_r$ with
$H$-prime divisors $D_i$. 
Then the vanishing ideal $\mathfrak{a}_i$
of $D_i$ is of $K$-height one,
and~(i) guarantees that it is principal,
say $\mathfrak{a}_i  = \bangle{f_i}$.
Thus $D_i = \div(f_i)$ and $D = \div(f_1 \cdots f_r)$
hold, which proves~(ii).

Assume that (ii) holds. Given a homogeneous
element $0 \ne f \in R \setminus R^*$,
write $\div(f) = D_1 + \ldots +D_r$ with 
$H$-prime divisors $D_i$. 
Then $D_i = \div(f_i)$ holds,
where, because of $R^* = \KK^*$,
the elements $f_i$ 
are homogeneous.
One verifies directly 
that the $f_i$ are $K$-prime. 
Thus we have 
$f=\alpha f_1 \cdots f_r$ with $\alpha \in \KK^*$
as required in~(iii). 

If~(iii) holds and $\mathfrak{a}$ is a 
$K$-prime ideal of $K$-height one,
then we take any homogeneous 
$0 \ne f \in \mathfrak{a}$ 
and find a $K$-prime factor $f_1$ 
of $f$ with $f_1 \in \mathfrak{a}$. 
This gives inclusions
$0 \subsetneq \bangle{f_1} \subseteq \mathfrak{a}$ 
of $K$-prime ideals,
which implies $\mathfrak{a}=\bangle{f_1}$.
\end{proof}

\begin{corollary}
\label{facthfact}
Under the assumptions of Proposition~\ref{hfactorial},
factoriality of the algebra $R$ implies that it is 
factorially graded.
\end{corollary}

We are ready to study the divisibility theory 
of the Cox ring. Here comes the main result;
it applies in particular to complete varieties,
see Corollary~\ref{corgradfactcomplete}.

\begin{theorem}
\label{thm:crfactgrad}
Let $X$ be an irreducible normal prevariety 
of affine intersection with 
only constant invertible functions and
finitely generated divisor class group.
If the Cox ring $\mathcal{R}(X)$ satisfies 
$\mathcal{R}(X)^* = \KK^*$, 
then it is factorially graded.
\end{theorem}

\begin{lemma}
\label{lem:smoothtorsionpull}
In the situation of Construction~\ref{constr:univtorssmooth},
every non-zero element 
$f \in \Gamma(X,\mathcal{R}_{[D]})$ 
satisfies
\begin{eqnarray*}
\div(f)
& = & 
q_X^*(\div_{[D]}(f)), 
\end{eqnarray*}
where on the left hand side $f$ 
is a regular function on $\rq{X}$ 
and on the right hand side~$f$ is 
an element on $\mathcal{R}(X)$.
\end{lemma}

\begin{proof}
In the notation of~\ref{constr:univtorssmooth},
let $D \in K$ represent $[D] \in \Cl(X)$,
and let $\t{f} \in \Gamma(X,\mathcal{S}_D)$ 
project to 
$f \in \Gamma(X,\mathcal{R}_{[D]})$.
The commutative diagram 
of~\ref{constr:univtorssmooth} yields
$$ 
\div(f) 
\ = \ 
\imath^*(\div(\t{f}))
\ = \ 
\imath^*(p^*(\div_{D}(\t{f})))
\ = \ 
q_X^*(\div_{[D]}(f)),
$$
where $\imath \colon \rq{X} \to \t{X}$ denotes the inclusion
and the equality 
$\div(\t{f})= p^*(\div_{D}(\t{f}))$ 
was established
in Lemma~\ref{smoothpull}.
\end{proof}

\begin{lemma}
\label{prop:torshfact}
In the situation of Construction~\ref{constr:univtorssmooth},
the prevariety $\rq{X}$ is irreducible, smooth and $H$-factorial.
\end{lemma}

\begin{proof}
As remarked in Construction~\ref{constr:univtorssmooth},
the prevariety $\rq{X}$ is smooth
and due to Proposition~\ref{prop:crnormal},
it is irreducible.
Let $\rq{D}$ be an invariant Weil divisor on~$\rq{X}$.
Using, for example, the fact 
that $q_X \colon \rq{X} \to X$ is an 
\'{e}tale principal bundle, 
we see that $\rq{D} = q_X^*(D)$ holds
with a Weil divisor $D$ on $X$.
Thus, we have to show that all pullback
divisors $q_X^*(D)$ are principal.
For this, it suffices to 
consider effective divisors
$D$ on $X$, and these are treated by 
Lemmas~\ref{lem:dddivprops} 
and~\ref{lem:smoothtorsionpull}.
\end{proof}

\begin{proof}[Proof of Theorem~\ref{thm:crfactgrad}]
According to Lemma~\ref{rem:codim2divsheaf2},
we may assume that $X$ is smooth.
Then $\mathcal{R}(X)$ is the algebra 
of regular functions of the  quasiaffine 
variety $\rq{X}$ constructed 
in~\ref{constr:univtorssmooth}.
Lemma~\ref{prop:torshfact} tells us 
that $\rq{X}$ is irreducible, smooth
and $H$-factorial.
Thus, Proposition~\ref{hfactorial} gives 
the assertion.
\end{proof}

\begin{corollary}
Let $X$ be a normal prevariety of affine intersection
with $\Gamma(X,\mathcal{O}) = \KK$ 
and finitely generated divisor class group.
Then the Cox ring $\mathcal{R}(X)$ 
is factorially graded.
\end{corollary}

\begin{proof}
According to Proposition~\ref{prop:inv2const2},
the assumption $\Gamma(X,\mathcal{O}) = \KK$
ensures $\mathcal{R}(X)^* = \KK^*$.
Thus Theorem~\ref{thm:crfactgrad} applies.
\end{proof}

\begin{corollary}
\label{corgradfactcomplete}
Let $X$ be a complete normal variety 
with finitely generated divisor class group.
Then the Cox ring $\mathcal{R}(X)$ 
is factorially graded.
\end{corollary}

As in the torsion free case,
see Proposition~\ref{cor:divdivalg}, 
divisibility and primality of homogeneous 
elements in the Cox ring $\mathcal{R}(X)$ 
can be characterized in terms of $X$.

\begin{proposition}
Let $X$ be a normal prevariety of affine intersection
with only constant invertible functions 
and finitely generated divisor class group.
Suppose that the Cox ring $\mathcal{R}(X)$ satisfies 
$\mathcal{R}(X)^* = \KK^*$.
\begin{enumerate}
\item
An element $0 \ne f \in \Gamma(X,\mathcal{R}_{[D]})$ 
divides  $0 \ne g \in \Gamma(X,\mathcal{R}_{[E]})$ 
if and only if $\div_{[D]}(f) \le \div_{[E]}(g)$
holds.
\item 
An element $0 \ne f \in  \Gamma(X,\mathcal{R}_{[D]})$ 
is $\Cl(X)$-prime
if and only if the divisor $\div_{[D]}(f) \in \WDiv(X)$ 
is prime. 
\end{enumerate}
\end{proposition}

\begin{proof}
According to Lemma~\ref{rem:codim2divsheaf2},
we may assume that $X$ is smooth.
Then Construction~\ref{constr:univtorssmooth}
presents $X$ as the geometric quotient 
of the smooth quasiaffine $H_X$-variety
$\rq{X}$, which has $\mathcal{R}(X)$ 
as its algebra of regular functions.
The first statement follows immediately 
from Lemma~\ref{lem:smoothtorsionpull}
and, for the second one, we additionally 
use Proposition~\ref{hfactorial}.
\end{proof}

\begin{remark}
Let $X$ be a prevariety of affine intersection
with only constant invertible functions,
finitely generated divisor class group 
and a Cox ring $\mathcal{R}(X)$ with only 
constant invertible elements.
Then the assignement 
$f \mapsto \div_{[D]}(f)$ induces 
an isomorphism from the 
multiplicative semigroup of homogeneous 
elements of $\mathcal{R}(X)$ modulo units 
onto the semigroup $\WDiv^+(X)$ of 
effective Weil divisors on $X$.
The fact that $\mathcal{R}(X)$ is 
factorially graded
reflects the fact that every effective 
Weil divisor 
is a unique non-negative linear combination 
of prime divisors.
\end{remark}

\begin{remark}
For the affine surface $X$ 
considered in Example~\ref{exam:sl2n},
the Cox ring $\mathcal{R}(X)$ 
is factorially $\ZZ/2\ZZ$-graded but 
it is not a factorial ring.
\end{remark}


\section{Geometric realization of the Cox sheaf}
\label{sec:univtorsors}

\subsection{Characteristic spaces}
We study the geometric realization of a
Cox sheaf, its relative spectrum,
which we call a {\em characteristic space\/}.
For locally factorial varieties, 
e.g. smooth ones, this concept coincides 
with the universal torsor
introduced by Colliot-Th\'el\`ene and Sansuc
in~\cite{CTSaA}, see also~\cite{CTSa}
and~\cite{Sk}.
As soon as we have non-factorial singularities, 
the characteristic 
space is not any more a torsor, i.e. a
principal bundle, as we will see later.
As before, we work with normal prevarieties
defined over an algebraically closed field of 
characteristic zero.
First we provide two statements on local
finite generation of Cox sheaves.

\begin{proposition}
\label{prop:crfg2csfg}
Let $X$ be a normal prevariety of 
affine intersection with only 
constant invertible functions
and finitely generated divisor class group.
If the Cox ring $\mathcal{R}(X)$ is finitely 
generated, then the Cox sheaf 
$\mathcal{R}$ is locally of finite type.
\end{proposition}

\begin{proof}
The assumption that $X$ is of affine intersection
guarantees that it is covered by open 
affine subsets $X_{[D],f}$, 
where $[D] \in \Cl(X)$ 
and $f \in \mathcal{R}_{[D]}(X)$.
By Proposition~\ref{prop:crloc}, 
we have 
$\Gamma(X_{[D],f},\mathcal{R}) = \mathcal{R}(X)_f$,
which gives the assertion.
\end{proof}

\begin{proposition}
\label{prop:qfact2csfg}
Let $X$ be a normal prevariety with only 
constant invertible functions
and finitely generated divisor class group.
If $X$ is $\QQ$-factorial, then any Cox sheaf
$\mathcal{R}$ is locally of finite type.
\end{proposition}

\begin{proof}
By definition, the Cox sheaf $\mathcal{R}$ 
is the quotient of a sheaf of divisorial 
algebras~$\mathcal{S}$ by some ideal 
sheaf $\mathcal{I}$.
According to Proposition~\ref{prop:coxsheafunique},
we may assume that $\mathcal{S}$ arises from 
a finitely generated subgroup $K \subseteq \WDiv(X)$.
Proposition~\ref{prop:qfact2locfin} then
tells us that $\mathcal{S}$ is locally of 
finite type,
and Lemma~\ref{lem:globfact} ensures 
that the quotient $\mathcal{R} = \mathcal{S}/\mathcal{I}$
can be taken at the level of sections.
\end{proof}

We turn to the relative spectrum of a Cox sheaf.
The following generalizes 
Construction~\ref{constr:univtorssmooth}, 
where the smooth case is considered.

\begin{construction}
\label{constr:univtors}
\index{space!characteristic}%
\index{characteristic space!}%
\index{characteristic quasitorus!}%
\index{quasitorus! characteristic}%
Let $X$ be a normal prevariety with 
$\Gamma(X,\mathcal{O}^*) = \KK^*$
and finitely generated divisor class 
group, 
and let $\mathcal{R}$ be a Cox sheaf.
Suppose that $\mathcal{R}$ is locally of finite type,
e.g., $X$ is ${\mathbb Q}$-factorial
or $\mathcal{R}(X)$ is finitely generated.
Taking the relative spectrum gives 
an irreducible normal prevariety
\begin{eqnarray*}
\rq{X}
& := & 
\Spec_X(\mathcal{R}).
\end{eqnarray*}
The $\Cl(X)$-grading of the sheaf $\mathcal{R}$ 
defines an action of the diagonalizable group
$H_X := \Spec \, \KK[\Cl(X)]$ on $\rq{X}$,  
the canonical morphism $q_X \colon \rq{X} \to X$
is a good quotient for this action,
and we have an isomorphism of sheaves
\begin{eqnarray*}
\mathcal{R}
& \cong &
(q_X)_* (\mathcal{O}_{\rq{X}}). 
\end{eqnarray*}
We call $q_X \colon \rq{X} \to X$ the 
{\em characteristic space\/} associated 
to $\mathcal{R}$, 
and $H_X$ the
{\em characteristic quasitorus\/} of~$X$. 
\end{construction}

\begin{proof}
Everything is standard except
irreducibility and normality,
which follow from Theorem~\ref{prop:crnormal}.
\end{proof}

The Cox sheaf $\mathcal{R}$ was defined as 
the quotient of a sheaf $\mathcal{S}$ of 
divisorial algebras by a sheaf $\mathcal{I}$
of ideals.
Geometrically this means that the characteristic 
space comes embedded into the relative 
spectrum of a sheaf of divisorial algebras;
compare~\ref{constr:univtorssmooth} for
the case of a smooth $X$.
Before making this precise in the general
case, we have to relate local finite 
generation of the 
sheaves $\mathcal{R}$ and $\mathcal{S}$ 
to each other.

\begin{proposition}
\label{prop:crlocfg}
Let $X$ be a normal prevariety with only 
constant invertible functions,
finitely generated divisor class group
and Cox sheaf $\mathcal{R}$.
Moreover, let $\mathcal{S}$ be 
the sheaf of divisorial algebras
associated to a finitely generated subgroup 
$K \subseteq \WDiv(X)$ projecting 
onto $\Cl(X)$ and 
$U \subseteq X$ an open affine subset.
Then the algebra $\Gamma(U,\mathcal{R})$ is 
finitely generated if and only if 
the algebra $\Gamma(U,\mathcal{S})$ is 
finitely generated.
\end{proposition}

\begin{proof}
Lemma~\ref{lem:globfact} tells us that 
$\Gamma(U,\mathcal{R})$ is a factor 
algebra of  $\Gamma(U,\mathcal{S})$.
Thus, if $\Gamma(U,\mathcal{S})$ is 
finitely generated then the same holds 
for $\Gamma(U,\mathcal{R})$.
Moreover, Lemma~\ref{lem:gradedcomps} 
says that the projection
$\Gamma(U,\mathcal{S}) \to \Gamma(U,\mathcal{R})$
defines isomorphisms along the homogeneous
components.
Thus, Proposition~\ref{prop:squeeze} shows 
that finite generation of $\Gamma(U,\mathcal{R})$
implies finite generation of $\Gamma(U,\mathcal{S})$.
\end{proof}

\begin{construction}
\label{rem:torsemb}
Let $X$ be a normal prevariety with 
$\Gamma(X,\mathcal{O}^*) = \KK^*$ 
and finitely generated divisor class 
group,
and let $K \subseteq \WDiv(X)$ be a 
finitely generated subgroup projecting onto $\Cl(X)$.
Consider the sheaf of divisorial algebras 
$\mathcal{S}$ associated to $K$ and the Cox sheaf
$\mathcal{R} = \mathcal{S} / \mathcal{I}$ 
as constructed in~\ref{constr:crtorsion},
and suppose that one of these sheaves
is locally of finite type.
Then the projection 
$\mathcal{S} \to \mathcal{R}$
of $\mathcal{O}_X$-algebras defines 
a commutative diagram
\begin{equation*} 
\xymatrix{
{\rq{X}}
\ar[rr]^{\imath}
\ar[dr]_{q_X}
&&
{\t{X}}
\ar[dl]^{p}
\\
&
X
&
}
\end{equation*}
for the relative spectra $\rq{X} = \Spec_X \mathcal{R}$
and $\t{X} = \Spec_X \mathcal{S}$.
We have the actions of $H_X = \Spec \, \KK[\Cl(X)]$ 
on $\rq{X}$ and $H = \Spec \, \KK[K]$ 
on $\t{X}$.
The map $\imath \colon \rq{X} \to \t{X}$ is a closed 
embedding and it is $H_X$-invariant,
where $H_X$ acts on $\t{X}$ via the inclusion 
$H_X \subseteq H$ defined by the projection 
$K \to \Cl(X)$.
The image $\imath(\rq{X}) \subseteq \t{X}$ 
is precisely the zero set 
of the ideal sheaf $\mathcal{I}$.
\end{construction}

\begin{proposition}
\label{prop:etalequot}%
\label{prop:smallmap2}%
\label{cor:quasaf}%
Situation as in Construction~\ref{constr:univtors}.
\begin{enumerate}
\item
The inverse image $q_X^{-1}(X_{\reg}) \subseteq \rq{X}$
of the set of smooth points $X_{\reg} \subseteq X$ 
is smooth, the group $H_X$ acts freely on 
$q_X^{-1}(X_{\reg})$ and the restriction 
$q_X \colon q_X^{-1}(X_{\reg}) \to X_{\reg}$
is an \'etale $H_X$-principal bundle. 
\item
For any closed $A \subseteq X$ 
of codimension at least two, 
$q_X^{-1}(A) \subseteq \rq{X}$
is as well of codimension at least two.
\item
The prevariety $\rq{X}$ is $H_X$-factorial.
\item
If $X$ is of affine intersection,
then $\rq{X}$ is a quasiaffine variety.
\end{enumerate}
\end{proposition}

\begin{proof}
For~(i), we refer to 
the proof of Construction~\ref{constr:univtorssmooth}.
To obtain~(ii) consider 
an affine open set $U \subseteq X$ 
and $\rq{U} := q_X^{-1}(U)$.
By Lemma~\ref{rem:codim2divsheaf2},
the open set $\rq{U} \setminus q_X^{-1}(A)$
has the same regular functions as $\rq{U}$.
Since $\rq{X}$ is normal, we conclude 
that $\rq{U} \cap q_X^{-1}(A)$ is of codimension 
at least two in $\rq{U}$.
Now, cover $X$ by affine $U \subseteq X$
and Assertion~(ii) follows.
We turn to~(iii).
According to~(ii)
we may assume that $X$ is smooth.
In this case, the statement was proven in
Lemma~\ref{prop:torshfact}.
We show~(iv). We may assume that we are 
in the setting of 
Construction~\ref{rem:torsemb}.
Corollary~\ref{cor:quaff2} then ensures that 
$\t{X}$ is quasiaffine and 
Construction~\ref{rem:torsemb} gives that~$\rq{X}$ 
is a closed subvariety of $\t{X}$.
\end{proof}

The following statement relates the
divisor of a $[D]$-homogeneous function 
on~$\rq{X}$ to its $[D]$-divisor on $X$;
the smooth case was settled in 
Lemma~\ref{lem:smoothtorsionpull}.

\goodbreak

\begin{proposition}
\label{prop:divalgspec2}
In the situation of~\ref{constr:univtors},
consider the pullback homomorphism
$q_X^* \colon \WDiv(X) \to \WDiv(\rq{X})$
defined in~\ref{rem:weilpull}.
Then, for every $[D] \in \Cl(X)$ 
and every $f \in \Gamma(X,\mathcal{R}_{[D]})$,
we have 
\begin{eqnarray*}
\div(f) & = & q_X^*(\div_{[D]}(f)),
\end{eqnarray*}
where on the left hand side 
$f$ is a function on $\rq{X}$, 
and on the right hand side a 
function on $X$.
If $X$ is of affine intersection and 
$X \setminus \Supp(\div_{[D]}(f))$
is affine, then we have moreover
\begin{eqnarray*}
\Supp(\div(f)) & = & q_X^{-1}(\Supp(\div_{[D]}(f))).
\end{eqnarray*}
\end{proposition}

\begin{proof}
We may assume that we are in the setting of 
Construction~\ref{rem:torsemb}.
Let the divisor $D \in K$ represent 
the class $[D] \in \Cl(X)$,
and let $\t{f} \in \Gamma(X,\mathcal{S}_D)$ 
project to 
$f \in \Gamma(X,\mathcal{R}_{[D]})$.
The commutative diagram 
of~\ref{rem:torsemb} yields
$$ 
\div(f) 
\ = \ 
\imath^*(\div(\t{f}))
\ = \ 
\imath^*(p^*(\div_{D}(\t{f})))
\ = \ 
q_X^*(\div_{[D]}(f)),
$$
where $\imath \colon \rq{X} \to \t{X}$ denotes the inclusion
and the equality 
$\div(\t{f})= p^*(\div_{D}(\t{f}))$ 
was established
in Proposition~\ref{prop:divalgspec}.
Similarly, we have
\begin{eqnarray*}
\Supp(\div(f)) 
& = &
\imath^{-1}(\Supp(\div(\t{f})))
\\
& = &  
\imath^{-1}(p^{-1}(\Supp(\div_{D}(\t{f}))))
\\
& = &  
q_X^{-1}(\Supp(\div_{[D]}(f)))
\end{eqnarray*}
provided that $X$ is of affine intersection 
and $X \setminus \Supp(\div_{[D]}(f))$
is affine, 
because  Proposition~\ref{prop:divalgspec}
then ensures
$\Supp(\div(\t{f}))
=
p^{-1}(\Supp(\div_{D}(\t{f})))$.
\end{proof}

\begin{corollary}
\label{cor:homogzeroes2}
Situation as in Construction~\ref{constr:univtors}.
Let $\rq{x} \in \rq{X}$ be a point such that 
$H_X \mal \rq{x} \subseteq \rq{X}$ is closed,
and let $f \in \Gamma(X,\mathcal{R}_{[D]})$.
Then we have 
\begin{eqnarray*}
f(\rq{x}) \ = \ 0
& \iff &
q_X(\rq{x}) \ \in \ \Supp(\div_{[D]}(f)).
\end{eqnarray*}
\end{corollary}

\begin{proof}
The image $q_X(\Supp(\div(f))$ 
is contained in $\Supp(\div_{[D]}(f))$.
By the definition of the pullback and 
Proposition~\ref{prop:divalgspec2},
the two sets coincide in $X_{\reg}$.
Thus, $q_X(\Supp(\div(f))$ is dense
in $\Supp(\div_{[D]}(f))$.
By Theorem~\ref{prop:GclosGsep},
the image 
$q_X(\Supp(\div(f))$
is closed and thus we have
\begin{eqnarray*}
q_X(\Supp(\div(f)))
& = & 
\Supp(\div_{[D]}(f)).
\end{eqnarray*}
In particular, if $f(\rq{x}) = 0$
holds, then $q_X(\rq{x})$ lies in 
$\Supp(\div_{[D]}(f))$.
Conversely, if $q_X(\rq{x})$ belongs 
to $\Supp(\div_{[D]}(f))$, then 
some $\rq{x}' \in \Supp(\div(f))$
belongs to the fiber of $\rq{x}$.
Since $H_x \mal \rq{x}$ is closed,
Corollary~\ref{prop:goodquotfibers}
tells us that $\rq{x}$ is contained
in the orbit closure of $\rq{x}'$ and hence
belongs to $\Supp(\div(f))$.
\end{proof}

\begin{corollary}
\label{cor:clos2clos}
Situation as in Construction~\ref{rem:torsemb}
and suppose that $X$ is of affine 
intersection.
For $x \in X$, let $\rq{x} \in q_X^{-1}(x)$
such that $H_X \mal \rq{x}$ is closed in
$\rq{X}$. 
Then $H \mal \rq{x}$ is closed in
$\t{X}$. 
\end{corollary}

\begin{proof}
Assume that the orbit $H \mal \rq{x}$ is not closed in
$\t{X}$. 
Then there is a point $\t{x} \in p^{-1}(x)$ 
having a closed $H$-orbit in $\t{X}$, 
and $\t{x}$ lies in the closure of 
$H \mal \rq{x}$.
Since $\t{X}$ is quasiaffine, we find 
a function $\t{f} \in \Gamma(X,\mathcal{S}_D)$
with $\t{f}(\t{x}) = 0$ but 
$\t{f}(\rq{x}) \ne 0$.
Corollary~\ref{cor:homogzeroes} gives 
$p(\t{x}) \in \Supp(\div_D(\t{f}))$.
Since we have $q_{X}(\rq{x}) = p(\t{x})$, 
this contradicts 
Corollary~\ref{cor:homogzeroes2}.
\end{proof}

\subsection{Divisor classes and isotropy groups}
The aim of this subsection is to
relate local properties of a prevariety 
to properties of the characteristic quasitorus 
action on its characteristic space.
Again, everything takes places over an 
algebraically 
closed field~$\KK$ of characteristic zero.

\index{local class group!}%
For a normal prevariety $X$ and a point $x \in X$,
let $\PDiv(X,x) \subseteq \WDiv(X)$ denote the 
subgroup of all Weil divisors, which 
are principal on some neighbourhood of~$x$. 
We define the {\em local class group\/} 
of $X$ in $x$ to be the factor group
\begin{eqnarray*}
\Cl(X,x) 
& := & 
\WDiv(X) / \PDiv(X,x).
\end{eqnarray*}
Obviously the group $\PDiv(X)$ of principal 
divisors is contained in $\PDiv(X,x)$. 
Thus, there is a canonical epimorphism 
$\pi_x \colon \Cl(X) \to \Cl(X,x)$.
The%
\index{Picard group!}
{\em Picard group\/} 
of $X$ is the factor group of the group $\CDiv(X)$ 
of Cartier divisors by the subgroup of principal 
divisors:
$$
\Pic(X) 
\ = \ 
\CDiv(X) / \PDiv(X)
\ = \ 
\bigcap_{x \in X} \ker(\pi_x).
$$

\begin{proposition}
\label{prop:divmonoid}
Situation as in Construction~\ref{constr:univtors}.
For $x \in X$, let $\rq{x} \in q_X^{-1}(x)$ 
be a point with closed $H_X$-orbit. Define a 
submonoid
$$ 
S_x
\ := \
\{[D] \in \Cl(X); \; f(\rq{x}) \ne 0 
       \text{ for some } 
       f \in \Gamma(X,\mathcal{R}_{[D]}) \} 
\ \subseteq \ \Cl(X),
$$
and let $\Cl_x(X) \subseteq \Cl(X)$ denote the 
subgroup generated by $S_x$.
Then the  local class groups of $X$ and the 
Picard group are given by
$$
\Cl(X,x) \ = \ \Cl(X)/\Cl_x(X),
\qquad\qquad
\Pic(X) 
\ = \
\bigcap_{x \in X} \Cl_x(X).
$$
\end{proposition}

\begin{proof}
First observe that Corollary~\ref{cor:homogzeroes2} 
gives us the following description of the monoid 
$S_x$ in terms of the $[D]$-divisors:
\begin{eqnarray*}
S_x
& = & 
\{[D] \in \Cl(X); \; x \not\in \div_{[D]}(f) 
       \text{ for some } 
       f \in \Gamma(X,\mathcal{R}_{[D]}) \} 
\\
& = & 
\{[D] \in \Cl(X); \; D \ge 0, \, x \not\in \Supp(D)\},
\end{eqnarray*}
where the latter equation is due to the fact that 
the $[D]$-divisors are precisely the 
effective divisors with class $[D]$. 
The assertions thus follow from
\begin{eqnarray*}
\Cl_x(X) 
& = & 
\{[D] \in \Cl(X); \; x \not\in \Supp(D)\}
\\
& = &
\{[D] \in \Cl(X); \; D \text{ principal near } x\}.
\end{eqnarray*}
\end{proof}

\begin{proposition}
\label{prop:isogrp2divisors2}
Situation as in Construction~\ref{constr:univtors}.
Given~$x \in X$, let  
$\rq{x} \in q_X^{-1}(x)$ be a point with 
closed $H_X$-orbit.
Then the inclusion $H_{X,\rq{x}} \subseteq H_X$ 
of the isotropy group of $\rq{x} \in \rq{X}$ 
is given by the epimorphism 
$\Cl(X) \to \Cl(X,x)$ of character groups.
In particular, we have 
$$
H_{X,\rq{x}} \ = \ \Spec \, \KK[\Cl(X,x)],
\qquad\qquad
\Cl(X,x) \ = \ \Chi(H_{X,\rq{x}}).
$$
\end{proposition}

\begin{proof}
Let $U \subseteq X$ be any affine open neighbourhood 
of $x \in X$.
Then $U$ is of the form $X_{[D],f}$ with some 
$f \in \Gamma(X,\mathcal{R}_{[D]})$
and $\t{U} := q_X^{-1}(U)$ is affine.
According to Proposition~\ref{prop:crloc}, we have
$$
\Gamma(\t{U}, \mathcal{O})
\ = \ 
\Gamma(U, \mathcal{R})
\ = \ 
\Gamma(X,\mathcal{R})_f
\ = \ 
\Gamma(\rq{X},\mathcal{O})_f.
$$
Corollary~\ref{cor:homogzeroes2} shows that 
the group $\Cl_x(X)$ is generated by the classes 
$[E] \in \Cl(X)$ admitting a section
$g \in \Gamma(U,\mathcal{R}_{[E]})$
with $g(\rq{x}) \ne 0$.
In other words, $\Cl_x(X)$ is the orbit group of 
the point $\rq{x} \in \t{U}$.
Now Proposition~\ref{prop:isogrorbitlattice}
gives the assertion.
\end{proof}

\index{factorial! point}%
\index{point! factorial!}%
\index{$\QQ$-factorial! point}%
\index{point! $\QQ$-factorial!}%
A point $x$ of a normal prevariety $X$ is called 
{\em factorial\/} if near $x$ every divisor is principal.
Thus, $x \in X$ is factorial if and only if 
its local ring $\mathcal{O}_{X,x}$ admits unique 
factorization.
Moreover, a point $x \in X$
is called {\em $\QQ$-factorial\/} if near
$x$ for every divisor some multiple is principal.

\goodbreak

\begin{corollary}
Situation as in Construction~\ref{constr:univtors}.
\begin{enumerate}
\item
A point $x \in X$ is factorial if and only if 
the fiber $q_X^{-1}(x)$ is a single $H_X$-orbit with 
trivial isotropy.
\item
A point $x \in X$ is $\QQ$-factorial if and only if 
the fiber $q_X^{-1}(x)$ is a single $H_X$-orbit.
\end{enumerate}
\end{corollary}

\begin{proof}
The point $x \in X$ is factorial if and only if $\Cl(X,x)$ 
is trivial, and it is  $\QQ$-factorial if and only if 
$\Cl(X,x)$ is finite.
Thus, the statement follows from 
Proposition~\ref{prop:isogrp2divisors2} 
and Corollary~\ref{prop:goodquotfibers}.
\end{proof}

\goodbreak

\begin{corollary}
\label{cor:univtorquot}
Situation as in Construction~\ref{constr:univtors}.
\begin{enumerate}
\item
The action of $H_X$ on $\rq{X}$ is free if and only if 
$X$ is locally factorial.
\item
The good quotient $q_X \colon \rq{X} \to X$ is geometric 
if and only if $X$ is $\QQ$-factorial.
\end{enumerate}
\end{corollary}

\begin{corollary}
\label{prop:picviacoxrings}
Situation as in Construction~\ref{constr:univtors}.
Let $\rq{H}_X \subseteq H_X$ 
be the subgroup generated by
all isotropy groups 
$H_{X,\rq{x}}$, where $\rq{x} \in \rq{X}$. 
Then we have 
\begin{eqnarray*}
\ker \bigl(\Chi(H_X) \to \Chi(\rq{H}_X)\bigr)
& = & 
\bigcap_{\rq{x} \in \rq{X}} \ker\bigl(\Chi(H_X) \to \Chi(H_{X,\rq{x}})\bigr)
\end{eqnarray*}
and the projection $H_X \to H_X/\rq{H}_X$ 
corresponds to the inclusion 
$\Pic(X) \subseteq \Cl(X)$ of character 
groups. 
\end{corollary}


\begin{corollary}
Situation as in Construction~\ref{constr:univtors}.
If the variety $\rq{X}$ contains an $H_X$-fixed point,
then the Picard group $\Pic(X)$ is trivial.
\end{corollary}

\subsection{Total coordinate space and irrelevant ideal}
Here we consider the situation that the Cox 
ring is finitely generated.
This allows us to introduce the total coordinate 
space as the spectrum of the Cox ring.
As always, we work over an 
algebraically 
closed field~$\KK$ of characteristic zero.

\begin{construction}
\label{constr:totcoorspac}
\index{total coordinate space}%
Let $X$ be a normal prevariety of affine 
intersection with only constant invertible 
functions and 
finitely generated divisor class group.
Let $\mathcal{R}$ be a Cox sheaf and 
assume that the Cox ring $\mathcal{R}(X)$
is finitely generated. Then we have a diagram
$$ 
\xymatrix{
{\Spec_X \mathcal{R}}
\ar@{=}[r]
&
{\rq{X}}
\ar[d]_{q_X}
\ar[r]^{\imath}
&
{\b{X}}
\ar@{=}[r]
&
{\Spec(\mathcal{R}(X))}
\\
&
X
&
&
&
}
$$
where the canonical morphism $\rq{X} \to \b{X}$ 
is an $H_X$-equivariant open embedding,
the complement $\b{X} \setminus \rq{X}$ 
is of codimension at least two
and $\b{X}$ is an $H_X$-factorial affine variety.
We call the $H_X$-variety $\b{X}$ the
{\em total coordinate space\/} associated to 
$\mathcal{R}$.
\end{construction}

\begin{proof}
Cover $X$ by affine open sets
$X_{[D],f} = X \setminus \Supp(\div_{[D]}(f))$,
where $[D] \in \Cl(X)$ and 
$f \in \Gamma(X,\mathcal{R}_{[D]})$.
Then, according to Proposition~\ref{prop:divalgspec2}, 
the variety $\rq{X}$ is covered 
by the affine sets $\rq{X}_f = q_X^{-1}(X_{[D],f})$.
Note that we have 
$$ 
\Gamma(\rq{X}_f,\mathcal{O})
\ = \ 
\Gamma(\rq{X},\mathcal{O})_f
\ = \ 
\Gamma(\b{X},\mathcal{O})_f
\ = \ 
\Gamma(\b{X}_f,\mathcal{O}).
$$
Consequently, the canonical morphisms 
$\rq{X}_f \to \b{X}_f$ are isomorphisms.
Gluing them together gives the desired 
open embedding $\rq{X} \to \b{X}$.
\end{proof}

\begin{definition}
\label{def:irrelideal}
\index{ideal! irrelevant}%
\index{irrelevant ideal!}%
Situation as in Construction~\ref{constr:totcoorspac}.
The {\em irrelevant ideal\/} of the prevariety $X$ is 
the vanishing ideal of the complement 
$\b{X} \setminus \rq{X}$ in the Cox ring:
$$ 
\mathcal{J}_{\irr}(X)
\ := \ 
\{f \in \mathcal{R}(X); \; f_{\vert \b{X} \setminus \rq{X}} = 0\}
\ \subseteq \
\mathcal{R}(X).
$$
\end{definition}

\begin{proposition}
\label{prop:irrelideal}
Situation as in Construction~\ref{constr:totcoorspac}.
\begin{enumerate}
\item
For any section $f \in \Gamma(X,\mathcal{R}_{[D]})$, 
membership in the irrelevant ideal is characterized 
as follows:
$$ 
f \in \mathcal{J}_{\irr}(X)
\ \iff \ 
\b{X}_f = \rq{X}_f
\ \iff \ 
\rq{X}_f
\text{ is affine}.
$$
\item 
Let $0 \ne f \in \Gamma(X,\mathcal{R}_{[D]})$.
If the $[D]$-localization 
$X_{[D],f}$ is affine, then we have 
$f \in \mathcal{J}_{\irr}(X)$.
\item
Let $0 \ne f_i \in \Gamma(X,\mathcal{R}_{[D_i]})$, where
$1 \le i \le r$ be such that the sets $X_{[D_i],f_i}$
are affine and cover $X$. Then we have 
\begin{eqnarray*}
\mathcal{J}_{\irr}(X)
& = & 
\sqrt{\bangle{f_1, \ldots, f_r}}.
\end{eqnarray*}
\end{enumerate}
\end{proposition}

\begin{proof}
The first equivalence in~(i) is obvious and 
the second one follows from the fact that
$\b{X} \setminus \rq{X}$ is of codimension 
at least two in $\b{X}$.
Proposition~\ref{prop:divalgspec2} tells us 
that for affine $X_{[D],f}$ also $\rq{X}_f$ is affine,
which gives~(ii).
We turn to (iii).
Proposition~\ref{prop:divalgspec2} and (ii)
ensure that the functions $f_1, \ldots,  f_r$ 
have $\b{X} \setminus \rq{X}$ as their common
zero locus. 
Thus Hilbert's Nullstellensatz gives the assertion. 
\end{proof}

\begin{corollary}
Situation as in Construction~\ref{constr:totcoorspac}.
Then $X$ is affine if and only if $\rq{X} = \b{X}$ 
holds.
\end{corollary}

\begin{proof}
Take $f = 1$ in the characterization~\ref{prop:irrelideal}~(i).
\end{proof}

\begin{corollary}
\label{cor:irrelidealqfact}
Situation as in Construction~\ref{constr:totcoorspac}
and assume that $X$ is $\QQ$-factorial.
Then $0 \ne f \in \Gamma(X,\mathcal{R}_{[D]})$
belongs to $\mathcal{J}_{\irr}(X)$ if and only if 
$X_{[D],f}$ is affine. In particular, we have 
\begin{eqnarray*}
\mathcal{J}_{\irr}(X)
& = &
\lin_{\KK}(
f \in \Gamma(X,\mathcal{R}_{[D]}); \;
[D] \in \Cl(X), \, 
X_{[D],f} \text{ is affine}
).
\end{eqnarray*}
\end{corollary}

\begin{proof}
We have to show that for any $[D]$-homogeneous 
$f \in \mathcal{J}_{\irr}(X)$, the $[D]$-localization
$X_{[D],f}$ is affine.
Note that $\rq{X}_f$ is affine by 
Proposition~\ref{prop:irrelideal}~(i).
The assumption of $\QQ$-factoriality ensures
that $q_X \colon \rq{X} \to X$ is a geometric
quotient, see Corollary~\ref{cor:univtorquot}.
In particular, all $H_X$-orbits in $\rq{X}$ 
are closed and thus Corollary~\ref{cor:homogzeroes2} 
gives us $\rq{X}_f = q_X^{-1}(X_{[D],f})$.
Thus, as the good quotient space 
of the affine variety $\rq{X}_f$,
the set $X_{[D],f}$ is affine. 
\end{proof}

\index{ample divisor!}%
\index{divisor! ample}%
Recall that a divisor $D$ on a prevariety $X$ is called 
{\em ample\/} if it admits sections 
$f_1, \ldots, f_r \in \Gamma(X,\mathcal{O}_X(D))$ such that 
the sets $X_{D,f_i}$ are affine and cover $X$.

\begin{corollary}
\label{cor:irrelidealample}
Situation as in Construction~\ref{constr:totcoorspac}.
If $[D] \in \Cl(X)$ is the class of an ample divisor,
then we have 
\begin{eqnarray*}
\mathcal{J}_{\irr}(X)
& = &
\sqrt{\bangle{ \Gamma(X,\mathcal{R}_{[D]})}}.
\end{eqnarray*}
\end{corollary}

\subsection{Characteristic spaces via GIT}
As we saw, the characteristic space of a 
prevariety~$X$ of affine intersection is 
a quasiaffine variety~$\rq{X}$ with an 
action of the characteristic quasitorus 
$H_X$ having $X$ as a good quotient. 
Our aim is to characterize this situation
in terms of Geometric Invariant Theory.
The crucial notion is the following.

\begin{definition}
\label{def:stronglystable}
\index{strongly stable action!}%
\index{action! strongly stable}%
Let $G$ be an affine algebraic group and 
$W$ a 
$G$-prevariety. 
We say that the $G$-action on $W$ is 
{\em strongly stable\/} if there is an 
open invariant subset $W' \subseteq W$ 
with the following properties:
\begin{enumerate}
\item
the complement $W \setminus W'$ is of codimension
at least two in $W$,
\item 
the group $G$ acts freely, i.e.~with 
trivial isotropy groups, on $W'$,
\item
for every $x \in W'$ the orbit $G \mal x$ 
is closed in $W$.
\end{enumerate}
\end{definition}

\begin{remark}
Let $X$ be a normal prevariety 
as in Construction~\ref{constr:univtors}
and consider the characteristic space
$q_X \colon \rq{X} \to X$ introduced
there.
Then Proposition~\ref{prop:smallmap2}
shows that the subset $q_X^{-1}(X_{\reg}) \subseteq \rq{X}$
satisfies the properties of~\ref{def:stronglystable}.
\end{remark}

\index{characteristic space!}%
\index{space! characteristic}%
Let $X$ and $q_X \colon \rq{X} \to X$ be 
as in Construction~\ref{constr:univtors}.
In the sequel, we mean by a 
{\em characteristic space\/} for $X$
more generally a good quotient 
$q \colon \mathcal{X} \to X$ 
for an action of a diagonalizable group~$H$ 
on a prevariety $\mathcal{X}$ such that 
there is an equivariant isomorphism 
$(\mu,\t{\mu})$ making the following
diagram commutative
$$ 
\xymatrix{
{\mathcal{X}}
\ar[dr]_{q}
\ar[rr]^{\mu}
&&
{\rq{X}}
\ar[dl]^{q_{X}}
\\
&
X
&
}
$$
Recall that here $\mu \colon \mathcal{X} \to \rq{X}$ 
is an isomorphism of varieties and 
$\t{\mu} \colon H \to H_X$ is an isomorphism
of algebraic groups such that we always have
$\mu(h \mal x) = \t{\mu}(h) \mal \mu(x)$.
Note that a good quotient 
$q \colon \mathcal{X} \to X$ of a quasiaffine 
$H$-variety is a characteristic space if and only 
if we have 
an isomorphism of graded sheaves
$\mathcal{R} \to q_*(\mathcal{O}_{\mathcal{X}})$,
where $\mathcal{R}$ is a Cox sheaf on $X$.

\begin{theorem}
\label{bigfreequot}
Let a quasitorus $H$ act on a normal quasiaffine 
variety $\mathcal{X}$ with a good quotient
$q \colon \mathcal{X} \to X$.
Assume that 
$\Gamma(\mathcal{X},\mathcal{O}^*) = \KK^*$ holds, 
$\mathcal{X}$ is $H$-factorial
and the $H$-action is strongly stable.
Then $X$ is a normal prevariety of affine intersection,
$\Gamma(X,\mathcal{O}^*) = \KK^*$ holds,
$\Cl(X)$ is finitely generated,
the Cox sheaf of $X$ is locally of finite type,
and $q \colon \mathcal{X} \to X$
is a characteristic space for $X$.
\end{theorem}

The proof will be given later in this section.
\index{total coordinate space!}%
First we also generalize the concept 
of the {\em total coordinate space\/}
of a prevariety $X$ of affine intersection 
with finitely generated Cox ring 
$\mathcal{R}(X)$:
this is from now on any affine $H$-variety
isomorphic to the affine $H_X$-variety 
$\b{X}$ of Construction~\ref{constr:totcoorspac}.

\begin{corollary}
\label{prop:geomchartotcorsp}
Let $Z$ be a normal affine 
variety with an action of 
a quasitorus $H$.
Assume that every invertible function on
$Z$ is constant, 
$Z$ is $H$-factorial, and there exists an
open $H$-invariant subset $W\subseteq Z$ 
with $\codim_Z (Z \setminus W) \ge 2$ such that
the $H$-action on $W$ is strongly stable and admits a
good quotient $q \colon W \to X$.
Then $Z$ is a  total coordinate
space for $X$.
\end{corollary}

A first step in the proof of Theorem~\ref{bigfreequot} 
is to describe
the divisor class group of the quotient 
space. 
Let us prepare the corresponding statement.
Consider an irreducible prevariety 
$\mathcal{X}$ with an action of a 
quasitorus $H = \Spec \, \KK[M]$.
For any $H$-invariant morphism 
$q \colon \mathcal{X} \to X$
to an irreducible prevariety $X$, we 
have the push forward homomorphism
$$ 
q_* \colon \WDiv(\mathcal{X})^H \ \to \ \WDiv(X)
$$
from the invariant Weil divisors of $\mathcal{X}$ 
to the Weil divisors of $X$ 
sending an $H$-prime divisor $D \subseteq \mathcal{X}$ 
to the closure of its image $q(D)$ if the latter 
is of codimension one and to zero else.
\index{rational homogeneous function!}%
\index{function! rational homogeneous}%
By a {\em homogeneous rational function\/} 
we mean an element $f \in \KK(\mathcal{X})$ that is 
defined on an invariant open subset of 
$\mathcal{X}$ and is homogeneous there. 
We denote the multiplicative group of 
non-zero homogeneous rational functions 
on $\mathcal{X}$ by $E(\mathcal{X})$ and 
the subset of non-zero rational functions 
of weight 
$w \in M$ by $E(\mathcal{X})_w$.

\goodbreak

\begin{proposition}
\label{prop:strstab2clgroup}
Let a quasitorus $H = \Spec \, \KK[M]$ 
act on a normal quasiaffine 
variety $\mathcal{X}$ with a good quotient
$q \colon \mathcal{X} \to X$.
Assume that
$\Gamma(\mathcal{X},\mathcal{O}^*) = \KK^*$ holds, 
$\mathcal{X}$ is $H$-factorial
and the $H$-action is strongly stable.
Then $X$ is a normal prevariety of affine 
intersection and there is an epimorphism
$$ 
\delta
\colon 
E(\mathcal{X})
\ \to \
\WDiv(X),
\qquad
f
\ \mapsto \ 
q_*(\div(f)).
$$
We have $\div(f) = q^*(q_*(\div(f)))$ for every 
$f \in E(\mathcal{X})$.
Moreover, the epimorphism 
$\delta$ induces a well-defined 
isomorphism
$$ 
M \ \to \ \Cl(X),
\qquad
w \ \mapsto \ [\delta(f)],
\quad
\text{with any } 
f \in E(\mathcal{X})_w.
$$
Finally, for every $f \in E(\mathcal{X})_w$, 
and every open set $U \subseteq X$, 
we have an isomorphism
of $\Gamma(U,\mathcal{O})$-modules
$$ 
\Gamma(U,\mathcal{O}_X(\delta(f)))
\ \to \ 
\Gamma(q^{-1}(U),\mathcal{O}_{\mathcal{X}})_w,
\qquad
g \ \mapsto \ fq^*(g).
$$
\end{proposition}

\begin{proof}
First of all note that the good quotient
space $X$ inherits normality and the 
property to be of affine intersection from 
the normal quasiaffine variety $\mathcal{X}$.

Let $\mathcal{X}' \subseteq \mathcal{X}$ be 
as in Definition~\ref{def:stronglystable}.
Then, with $X' := q(\mathcal{X}')$,
we have $q^{-1}(X') =  \mathcal{X}'$.
Consequently, $X' \subseteq X$ is open.
Moreover, $X \setminus X'$ is of codimension 
at least two in $X$,
because $\mathcal{X} \setminus \mathcal{X}'$
is of codimension at least two in $\mathcal{X}$.
Thus, we may assume that $X=X'$ holds, 
which means in particular that $H$ acts freely.
Then we have homomorphisms of groups:
$$ 
\xymatrix{
E(\mathcal{X})
\ar[rr]^{f \mapsto \div(f) \qquad}
&&
{\WDiv(\mathcal{X})^H}
\ar@/^/[rr]^{q_*}
\ar@{<-}@/_/[rr]_{q^*}
&&
{\WDiv(X)}.
}
$$
The homomorphism from $E(\mathcal{X})$ 
to the group
of $H$-invariant Weil divisors
$\WDiv(\mathcal{X})^H$ 
is surjective, because $\mathcal{X}$ 
is $H$-factorial.
Moreover, $q^*$ and $q_*$ 
are inverse to each other,
which follows 
from the observation that 
$q \colon \mathcal{X} \to X$ is an 
\'{e}tale $H$-principal bundle.
This establishes the first part 
of the assertion.

We show that $\delta$ 
induces an isomorphism $M \to \Cl(X)$.
First we have to check that $[\delta(f)]$ 
does not depend on the choice of $f$.
So, let $f,g \in E(\mathcal{X})_w$.
Then $f/g$ is $H$-invariant, 
and hence defines a rational function 
on $X$.
We infer well-definedness of $w \mapsto [\delta(f)]$ 
from
$$ 
q_*(\div(f)) - q_*(\div(g))
\ = \ 
q_*(\div(f) - \div(g))
\ = \ 
q_*(\div(f/g))
\ = \ 
\div(f/g).
$$
To verify injectivity, let 
$\delta(f) = \div(h)$ 
for some $h \in \KK(X)^*$. 
Then we obtain $\div(f) = \div(q^*(h))$.
Thus, $f/q^*(h)$ is an invertible homogeneous
function on $\mathcal{X}$ 
and hence is constant.
This implies $w = \deg(f/q^*(h)) = 0$.
Surjectivity is clear, because
$E(\mathcal{X}) \to \WDiv(X)$ 
is surjective.

We turn to the last statement. 
First we note that
for every $g \in \Gamma(U,\mathcal{O}_X(\delta(f)))$
the function $fq^*(g)$ is regular on 
$q^{-1}(U)$, because we have 
$$ 
\div(fq^*(g))
\, = \, 
\div(f) + \div(q^*(g))
\, = \,
q^*(\delta(f)) + \div(q^*(g))
\, =  \,
q^*(\delta(f)+\div(g))
\, \ge \, 0. 
$$
Thus, the homomorphism
$\Gamma(U,\mathcal{O}_X(\delta(f))) 
\to 
\Gamma(q^{-1}(U),\mathcal{O}_{\mathcal{X}})_w$
sending $g$ to $fq^*(g)$ is well defined.
Note that $h \mapsto h/f$ defines an inverse.
\end{proof}

\begin{corollary}
Consider the characteristic space $q \colon \rq{X} \to X$ 
obtained from a Cox sheaf $\mathcal{R}$.
Then, for any non-zero 
$f \in \Gamma(X,\mathcal{R}_{[D]})$ 
the push forward $q_*(\div(f))$,
equals the $[D]$-divisor $\div_{[D]}(f)$.
\end{corollary}

\begin{proof}
Proposition~\ref{prop:strstab2clgroup} shows that 
$q^*(q_*(\div(f))$ equals $\div(f)$ and 
Proposition~\ref{prop:divalgspec2} tells 
us that $q^*(\div_{[D]}(f))$ equals $\div(f)$ 
as well.
\end{proof}

\begin{proof}[Proof of Theorem~\ref{bigfreequot}]
Writing $H = \Spec \, \KK[M]$ with the character
group~$M$ of $H$, we are in the setting 
of Proposition~\ref{prop:strstab2clgroup}.
Choose a finitely generated subgroup 
$K \subseteq \WDiv(X)$ mapping onto $\Cl(X)$,
and let $D_1, \ldots, D_s \in \WDiv(X)$ be a basis 
of $K$. 
By Proposition~\ref{prop:strstab2clgroup},
we have $D_i = \delta(h_i)$ 
with $h_i \in E(\mathcal{X})_{w_i}$.
Moreover, the isomorphism $M \to \Cl(X)$ 
given there identifies $w_i \in M$ 
with $[D_i] \in \Cl(X)$.
For $D = a_1D_1 + \ldots + a_sD_s$,
we have $D = \delta(h_D)$ with 
$h_D = h_{1}^{a_1} \cdots  h_{s}^{a_s}$. 

Let $\mathcal{S}$ be 
the sheaf of divisorial algebras associated to $K$
and for $D \in K$, let $w \in M$ correspond to 
$[D] \in \Cl(X)$.
Then, for any open set $U \subseteq X$ and any $D \in K$, 
Proposition~\ref{prop:strstab2clgroup} 
provides an
isomorphism of $\KK$-vector spaces
$$ 
\Phi_{U,D} \, \colon \,
\Gamma(U, \mathcal{S}_{D}) 
\ \to \  
\Gamma(q^{-1}(U), \mathcal{O})_w,
\qquad
g \ \mapsto \ q^*(g) \, h_D.
$$
The $\Phi_{U,D}$ fit together to an 
epimorphism of graded sheaves  
$\Phi \colon \mathcal{S} \to q_*(\mathcal{O}_{\mathcal{X}})$.
Once we know that $\Phi$ has the ideal $\mathcal{I}$
of Construction~\ref{constr:crtorsion} as its kernel,
we obtain an induced isomorphism 
$\mathcal{R} \to q_* \mathcal{O}_{\mathcal{X}}$,
where $\mathcal{R} = \mathcal{S} / \mathcal{I}$ 
is the associated Cox sheaf;
this shows that $\mathcal{R}$ is locally of 
finite type and gives an isomorphism
$\mu \colon \mathcal{X} \to \rq{X}$.

Thus we are left with showing that the 
kernel of $\Phi$ equals $\mathcal{I}$.
Consider a $\Cl(X)$-homogeneous  element 
$f \in \Gamma(U, \mathcal{S})$
of degree $[D]$, where $D\in K$. 
Let $K^0$ be the kernel of the surjection 
$K \to \Cl(X)$. Then we have 
$$
f \ = \ \sum_{E\in K^0} f_{D+E}, 
\qquad\qquad 
\Phi(f) \ = \ \sum_{E\in K^0} q^*(f_{D+E}) \, h_{D+E}. 
$$
With the character $\chi \colon K^0 \to \KK(X)^*$
defined by $q^*\chi(E) = h_{E}$, we may rewrite the 
image $\Phi(f)$ as 
$$
\Phi(f) 
\ = \ 
\sum_{E\in K^0} q^*(\chi(E)f_{D+E}) \, h_D
\ = \ 
q^*
\left(
\sum_{E\in K^0} \chi(E)  f_{D+E}
\right) 
\, h_D.
$$
So, $f$ lies in the kernel of $\Phi$ 
if and only if $\sum \chi(E)  f_{D+E}$
vanishes. 
Now observe that we have
\begin{eqnarray*}
f 
& = &
\sum_{E\in K^0} (1- \chi(E)) \, f_{D+E}  
\ + \ 
\sum_{E\in K^0} \chi(E) \, f_{D+E}. 
\end{eqnarray*}
The second summand is $K$-homogeneous, 
and thus we infer from 
Lemma~\ref{lem:secpres} that 
$f \in \mathcal{I}$ holds if and only if 
$\sum \chi(E) \, f_{D+E} = 0$
holds.
\end{proof}

\begin{remark}
\label{rem:bigfreequotdetails}
Consider the isomorphism $(\mu,\t{\mu})$
identifying the characteristic spaces 
$q \colon \mathcal{X} \to X$ and 
$q_X \colon \rq{X} \to X$ in the above proof.
Then the isomorphism $\t{\mu}$ identifying 
the quasitori $H$ and $H_X$ is given 
by the isomorphism $M \to \Cl(X)$ of their
character groups provided 
by Proposition~\ref{prop:strstab2clgroup}.
\end{remark}

\cleardoublepage


\cleardoublepage

\input{CoxRings_Chapter_1.ind}


\begin{thebibliography}{}%
%
\bibitem{AAC}
 A.~A'Campo-Neuen:
Note on a counterexample to Hilbert's fourteenth problem given by P. Roberts.  
Indag. Math. (N.S.) 5 (1994), no.~3, 253--257.
%
 \bibitem{ACHa1}
 A.~A'Campo-Neuen, J.~Hausen:
 Quotients of toric varieties by the action of a subtorus.
 Tohoku Math. J. (2) 51 (1999), no.~1, 1--12.
%
 \bibitem{ACHa}
 A.~A'Campo-Neuen, J.~Hausen:
 Toric prevarieties and subtorus actions.
 Geom. Dedicata 87  (2001),  no.~1-3, 35--64.
%
 \bibitem{ACHaS} 
 A.~A'Campo-Neuen, J.~Hausen, S.~Schr\"oer:
 Homogeneous coordinates and quotient presentations for toric varieties.
 Math. Nachr.~246/247 (2002), 5--19.
%
 \bibitem{AlNi}
 V.A.~Alekseev, V.V~Nikulin:
 Classification of log del Pezzo surfaces  of index $\le 2$.
 Memoirs of the Mathematical Society of Japan, vol.~15 (2006),
 Preprint Version arXiv:math.AG/0406536.
%
 \bibitem{AlNi2}
 V.A.~Alexeev, V.V.~Nikulin:
 Del Pezzo and K3-surfaces.
 MSJ Memoirs, 15. Mathematical Society of Japan,
 Tokyo, 2006.
%
 \bibitem{AlHa}
 K.~Altmann, J.~Hausen:
 Polyhedral divisors and algebraic torus actions.
 Math. Ann. 334 (2006), no.~3, 557--607.
%
 \bibitem{AlHaSu} 
 K.~Altmann, J.~Hausen, H.~S\"u{\ss}:
 Gluing affine torus actions via divisorial fans.
 Transformation Groups 13 (2008), no.~2, 215--242.
%
 \bibitem{Anders} 
D.F.~Anderson:
Graded Krull domains.
Comm. Algebra 7 (1979), no.~1, 79--106.
%
 \bibitem{ArCoGrHa}
 E.~Arbarello, M.~Cornalba, P.A.~Griffiths, J.~Harris:
 Geometry of algebraic curves.
 Grundlehren der Mathematischen Wissenschaften
 267, Springer-Verlag, New York, 1985
%
 \bibitem{ArHaLa}
 M.~Artebani, J.~Hausen, A.~Laface:
 On Cox rings of K3-surfaces.
 To appear in Compositio Math., arXiv:0901.0369.
%
 \bibitem{ArLa}
 M.~Artebani, A.~Laface:
 Cox rings of surfaces and the anticanonical Iitaka dimension.
 Preprint, arXiv:0909.1835.
%
 \bibitem{Art}
 M.~Artin:
 Some numerical criteria for contractability of curves on algebraic surfaces.
 Amer. J. Math. 84 (1962), 485--496.
%
 \bibitem{Ar1}
 I.V.~Arzhantsev
 Projective embeddings of homogeneous spaces with small boundary.
 Izvestiya RAN: Ser. Mat. 73 (2009), no.~3, 3--22 (Russian);
 English transl.: Izvestya Mathematics 73 (2009), no.~3, 437--453.
%
 \bibitem{Ar2} 
 I.V.~Arzhantsev: 
 On the factoriality of Cox rings.
 Mat. Zametki 85 (2009), no.~5, 643--651 (Russian);
 English transl.: Math. Notes 85 (2009), no.~5, 623--629.
%
 \bibitem{ArGa}
 I.V.~Arzhantsev, S.A.~Gaifullin:
 Cox rings, semigroups and automorphisms of affine algebraic varieties.
 Mat. Sbornik 201 (2010), no.~1, 3--24 (Russian);
 English transl.: Sbornik Math. 201 (2010), no.~1, 1--21.
%
 \bibitem{ArHa1} 
 I.V.~Arzhantsev, J.~Hausen:
 On embeddings of homogeneous spaces with small boundary.
 J.~Algebra~304 (2006), no.~2, 950--988.
%
 \bibitem{ArHa2} 
 I.V.~Arzhantsev, J.~Hausen:
 On the multiplication map of a multigraded algebra.
 Math. Res. Lett. 14 (2007), no.~1, 129--136.
%
 \bibitem{ArHa3} 
 I.V.~Arzhantsev, J.~Hausen:
 Geometric Invariant Theory via Cox rings.
 J. Pure Appl. Algebra 213 (2009), no.~1, 154--172.
%
%
 \bibitem{AtMa}
 M.F.~Atiyah, I.G.~Macdonald:
 Introduction to commutative algebra.
 Addison-Wesley Publishing Co., Reading, Mass.-London-Don Mills, Ont. 1969
%
 \bibitem{Au}
 M.~Audin:
 The topology of torus actions on symplectic manifolds.
 Prog. Math., 93. Birkh\"auser Verlag, Basel, 1991
%
\bibitem{Bae}
H.~B\"aker:
Good quotients of Mori Dream Spaces. 
Proc. Amer. Math. Soc. ? (2010), no.~?, ??-??.
%
 \bibitem{BaHuPeVa}
 W.P.~Barth, K.~Hulek, C.A.M.~Peters, A.~Van de Ven:
 Compact complex surfaces.
 Ergebnisse der Mathematik und ihrer Grenzgebiete. 3. Folge. 
 A Series of Modern Surveys in Mathematics 4, Second edition, 
 Springer-Verlag, Berlin, 2004
%
 \bibitem{Ba}
 V.V.~Batyrev:
 Quantum cohomology rings of toric manifolds.
 In: Journ\'ees de G\'eom\'etrie Alg\'ebrique d'Orsay, Ast\'eriques~218, 9--34 (1993).
%
\bibitem{BaHa}
 V.V.~Batyrev, F.~Haddad:
 The geometry of SL(2)-equivariant flips.
 Mosc. Math. J. 8 (2008), no.~4, 621--646.
%
 \bibitem{BaPo}
 V.V.~Batyrev, O.N.~Popov:
 The Cox ring of a Del Pezzo surface.
 In: Arithmetic of higher-dimensional algebraic varieties, Progr. Math.~226,
 85--103 (2004) 
%
 \bibitem{Be}
 A.~Beauville:
 Complex algebraic surfaces. Second edition. London Mathematical Society Student
 Texts 34. Cambridge University Press, Cambridge, 1996
%
 \bibitem{BeHa1}
 F.~Berchtold, J.~Hausen:
 Homogeneous coordinates for algebraic varieties.
 J. Algebra 266 (2003), no.~2, 636--670.
%
 \bibitem{BeHa2}
 F.~Berchtold, J.~Hausen:
 Bunches of cones in the divisor class group --
 a new combinatorial language for toric varieties.
 Inter. Math. Research Notices 6 (2004), 261--302.
%
\bibitem{BeHa3}
 F.~Berchtold, J.~Hausen:
 GIT-equivalence beyond the ample cone.
 Michigan Math. J. 54 (2006), no.~3, 483--515.
%
 \bibitem{BeHa4} 
 F.~Berchtold, J.~Hausen:
 Cox rings and combinatorics.
 Trans. Amer. Math. Soc.  359  (2007), no.~3, 1205--1252.
%
%
 \bibitem{BB2}
 A.~Bia\l ynicki-Birula:
 Algebraic quotients.
 In: Encyclopaedia of Mathematical Sciences 131.
 Invariant Theory and Algebraic Transformation Groups, II.
 Springer-Verlag, Berlin, 2002
%
 \bibitem{BBSw82}
 A.~Bia\l ynicki-Birula, J.~\'Swi\c{e}cicka:
Complete quotients by algebraic torus actions. 
In: Group actions and vector fields (Vancouver, B.C., 1981),  10--22, 
Lecture Notes in Math. 956, Springer, Berlin, 1982
%
\bibitem{BBSw4}  
A.~Bia\l ynicki-Birula, J.~\'Swi\c{e}cicka: 
Three theorems on existence of good quotients. 
Math. Ann. 307 (1997), 143--149.
%
 \bibitem{BBSw}
 A.~Bia\l ynicki-Birula, J.~\'Swi\c{e}cicka:
 A recipe for finding open subsets of vector spaces with a good quotient.
 Colloq. Math. 77 (1998), no.~1, 97--114.
%
 \bibitem{BiBoI}
 F.~Bien, A.~Borel:
 Sous-groupes \'epimorphiques de groupes lin\'eaires alg\'ebriques I.
 C. R. Acad. Sci. Paris, t.~315, S\'erie I (1992), 649--653.
%
 \bibitem{BBII}
 F.~Bien, A.~Borel:
 Sous-groupes \'epimorphiques de groupes lin\'eaires alg\'ebriques II.
 C. R. Acad. Sci. Paris, t.~315, S\'erie I (1992), 1341--1346.
%
 \bibitem{BiBoKo}
 F.~Bien, A.~Borel, J.~Kollar:
 Rationally connected homogeneous spaces.
 Invent. Math. 124 (1996), no.~1-3, 103--127.
%
\bibitem{Bo}
A.~Borel:
Linear algebraic groups.
Second edition.
Graduate Texts in Mathematics 126,
Springer-Verlag, New York, 1991
%
\bibitem{Bour}
N.~Bourbaki: Commutative algebra. Chapters 1--7.
Elements of Mathematics, Springer-Verlag, Berlin, 1998

%
\bibitem{Bou} J.~Boutot: Singularit\'{e}s rationelles et quotients par les
groupes r\'{e}ductifs. Invent. Math. 88, 65--68 (1987)
%
\bibitem{Bre}
 H.~Brenner:
 Rings of global sections in two-dimensional schemes.
 Beitr\"age Algebra Geom. 42 (2001), no.~2, 443--450.
%
 \bibitem{Br}
 M.~Brion:
 The total coordinate ring of a wonderful variety.
 J. Algebra 313 (2007), no.~1, 61--99.
%
%
 \bibitem{BuPa}
 V.M.~Buchstaber, T.E.~Panov:
 Torus actions and their applications in topology and combinatorics.
 Univ. Lecture Series~24, Providence R.I., AMS, 2002
%
 \bibitem{Ca}
 A.-M.~Castravet:
 The Cox ring of $\overline M\sb {0,6}$.
 Trans. Amer. Math. Soc. 361 (2009), no.~7, 3851--3878.
%
 \bibitem{CaTe}
 A.-M.~Castravet, E.A.~Tevelev: 
 Hilbert's 14th problem and Cox rings.
 Compos. Math. 142 (2006), no.~6, 1479--1498.
%
 \bibitem{CiMa} 
 R.~Chirivi, A.Muffei:
 The ring of sections of a complete symmetric variety.
 J. Algebra 261 (2003), no.~2, 310--326.
%
 \bibitem{CoDo}
 F.R. Cossec and I.V. Dolgachev:
 Enriques surfaces I.
 Birkh\"auser, Progress in Mathematics, Vol. 76, 1989
%
 \bibitem{Co}
 D.A.~Cox:
 The homogeneous coordinate ring of a toric variety.
 J. Alg. Geom. 4 (1995), no.~1, 17--50.
%
 \bibitem{Co2}
 D.A.~Cox:
 Recent developments in toric geometry.
 Proc. Symp. Pure Math. 62 (1997), 389--436.
%
 \bibitem{CoLiSch} 
 D.A.~Cox, J.~Little, H.~Schenk:
Toric varieties, draft.
%
 \bibitem{CrMa}
 A.~Craw, D.~Maclagan:
 Fiber fans and toric quotients.
 Discrete Comput. Geom. 37 (2007), no.~2, 251--266.
%
\bibitem{CTSaA}
J.L.~Colliot-Th\'el\`ene, J.-J.~Sansuc:  
Torseurs sous des groupes de type multiplicatif; 
applications \`a l'\' etude des points rationnels 
de certaines vari\'et\'es alg\'ebriques.  
C. R. Acad. Sci. Paris S\'er. A-B  282  (1976), 
no. 18, Aii, A1113--A1116.
%
\bibitem{CTSaB}
J.L.~Colliot-Th\'el\`ene, J.-J.~Sansuc:  
Vari\'et\'es de premi\`ere descente attach\'ees 
aux vari\'et\'es rationnelles. 
C. R. Acad. Sci. Paris S\'er. A-B  284  (1977), 
no. 16, A967--A970.
%
\bibitem{CTSaC}
J.L.~Colliot-Th\'el\`ene, J.-J.~Sansuc:  
La descente sur une vari\'et\'e rationnelle d\'efinie 
sur un corps de nombres. 
C. R. Acad. Sci. Paris S\'er. A-B  284  (1977), 
no. 19, A1215--A1218.
%
\bibitem{CTSa}
 J.L.~Colliot-Th\'el\`ene, J.-J.~Sansuc:
 La descente sur le vari\'et\'es rationnelles. II.
 Duke Math. J. 54 (1987), no.~2, 375--492.
%
 \bibitem{Da}
 V.I.~Danilov:
 The geometry of toric varieties.
 Uspekhi Mat. Nauk 33 (1978), no.~2(200), 85--134 (Russian);
 English transl: Russian Math. Surveys 33 (1978), no.~2, 97--154.
%
\bibitem{DLRS}
J.A.~De~Loera, J.~Rambau, F:~Santos:
Triangulations  -- Structures for Algorithms and Applications
Series: Algorithms and Computation in Mathematics, Vol. 25.
Springer Verlag, 2010.
%
 \bibitem{Del}
 T.~Delzant:
 Hamiltoniens p\'eriodiques et images convexes de l'application moment.
 Bull. Soc. Math. France 116 (1988), no.~3, 315--339.
%
%
 \bibitem{De3}
 U.~Derenthal:
 Singular Del Pezzo surfaces whose universal torsors are hypersurfaces.
 Preprint, arXiv:math.AG/0604194.
%
 \bibitem{De2}
 U.~Derenthal:
 Universal torsors of Del Pezzo surfaces and homogeneous spaces.
 Adv. Math. 213 (2007), no.~2, 849--864.
%
 \bibitem{DeTs}
 U.~Derenthal, Yu.~Tschinkel:
 Universal torsors over Del Pezzo surfaces and rational points.
 In "Equidistribution in Number theory, An Introduction", 
 (A. Granville, Z. Rudnick eds.),
 169-196, NATO Science Series II, 237, Springer, 2007
%
 \bibitem{Do}
 I.V.~Dolgachev:
 Rational surfaces with a pencil of elliptic curves.
 Izv. Akad. Nauk SSSR  Ser. Mat., 30 (1966), no.~5, 1073--1100 (Russian).
%
 \bibitem{DoHu}
 I.V.~Dolgachev, Y.~Hu:
 Variation of geometric invariant theory quotients.
 With an appendix by Nicolas Ressayre.
 Inst. Hautes Etudes Sci. Publ. Math. 87 (1998), 5--56.
%
 \bibitem{EiPo}
 D.~Eisenbud, S.~Popescu:
 The projective geometry of the Gale transform.
 J.~Algebra 230 (2000), no.~1, 127--173.
%
 \bibitem{Ei}
 D.~Eisenbud:
 Commutative Algebra with a view towards algebraic geometry.
 Graduate texts in Math. 150, Springer Verlag, New York, 1995
%
 \bibitem{El}
 E.J.~Elizondo:
 Chow varieties, the Euler-Chow series and the total coordinate ring.
 In: Transcendental aspects of algebraic cycles, 3--43, London Math. 
 Soc. Lecture Note Ser. 313, Cambridge Univ. Press, Cambridge, 2004
%
 \bibitem{ElKuWa}
 E.J.~Elizondo, K.~Kurano, K.~Watanabe:
 The total coordinate ring of a normal projective variety.
 J. Algebra 276 (2004), no.~2, 625--637.
%
 \bibitem{Ew}
 G.~Ewald:
 Polygons with hidden vertices.
 Beitr\"age Algebra Geom. 42 (2001), no.~2, 439--442.
%
 \bibitem{Fu}
 W.~Fulton:
 Introduction to Toric Varieties.
 Princeton Univ. Press, Princeton, 1993
%
 \bibitem{Fu2}
 W.~Fulton:
 Intersection theory. Second edition.
 Ergebnisse der Mathematik und ihrer Grenzgebiete. 
 3. Folge. A Series of Modern Surveys in Mathematics
 2. Springer-Verlag, Berlin, 1998
%
\bibitem{FuHa}
W.~Fulton, J.~Harris:
Representation theory. 
A first course. 
Graduate Texts in Mathematics, 129. 
Springer-Verlag, New York, 1991. 
%
 \bibitem{Ga}
 S.A.~Gaifullin:
 Affine toric SL(2)-embeddings.
 Mat. Sbornik 199 (2008), no.~3, 3--24 (Russian);
 English transl.: Sbornik Math. 199 (2008), no.~3-4, 319--339.
%
 \bibitem{GaMo}
 C.~Galindo, F.~ Monserrat:
 The total coordinate ring of a smooth projective surface.
 J. Algebra 284 (2005), no.~1, 91--101.
%
 \bibitem{GuMa}
 S.~Giuffrida, R.~Maggioni:
 The global ring of a smooth projective surface.
 Matematiche (Catania) 55 (2000), no.~1, 133--159.
%
 \bibitem{Gr}
 F.D.~Grosshans:
 Algebraic Homogeneous Spaces and Invariant Theory.
 LNM 1673, Springer-Verlag, Berlin, 1997
%
 \bibitem{Gu}
 N.~Guay:
 Embeddings of symmetric varieties.
 Transformation Groups 6 (2001), no.~4, 333--352.
%
 \bibitem{Har}
 R.~Hartshorne:
 Algebraic Geometry.
 GTM 52, Springer-Verlag, 1977
%
 \bibitem{Has}
 B.~Hassett: 
 Equations of universal torsors and Cox rings.  
 Math. Institut, Georg-August-Universit\"at G\"ottingen: 
Seminars Summer Term 2004,  135--143,
 Universit\"atsdrucke G\"ottingen, G\"ottingen (2004)
%
 \bibitem{HaTs} 
 B.~Hassett, Yu.~Tschinkel: 
 Universal torsors and Cox rings.
 In: Arithmetic of higher-dimensional algebraic varieties,  
 Progr. Math., 226, Birkh\"auser Boston, Boston, MA, 149--173 (2004)
%
%
 \bibitem{Ha0}
 J.~Hausen:
 Producing good quotients by embedding into toric varieties.  
 Geometry of toric varieties,  193--212, S\'emin. Congr., 6, 
 Soc. Math. France, Paris, 2002. 
%
 \bibitem{Ha02}
 J.~Hausen:
 Equivariant embeddings into smooth toric varieties.  
 Canad. J. Math. 54 (2002), no.~3, 554--570
%
\bibitem{Ha1}
 J.~Hausen:
 Geometric invariant theory based on Weil divisors.
 Compos. Math. 140 (2004), no.~6, 1518--1536.
%
 \bibitem{Ha2}
 J.~Hausen:
 Cox rings and combinatorics II.
 Mosc. Math. J. 8 (2008), no.~4,  711--757.
%
 \bibitem{HaSu}
 J.~Hausen, H.~S\"u\ss:
 The Cox ring of an algebraic varitey with torus action.
 To appear in Adv. Math., arXiv:0903.4789 
%
\bibitem{HR} M.~Hochster, J.~Roberts: Rings of invariants of reductive groups
  acting on regular rings are Cohen-Macaulay. Adv. Math. 13, 115--175 (1974)
%
 \bibitem{HuKe} 
 Y.~Hu, S.~Keel: 
 Mori dream spaces and GIT.
 Michigan Math. J. 48 (2000) 331--348.
%
\bibitem{Hu}
J.E.~Humphreys:
Linear algebraic groups. 
Graduate Texts in Mathematics, No. 21. 
Springer-Verlag, New York-Heidelberg, 1975
%
\bibitem{Kaj}
T.~Kajiwara: 
The functor of a toric variety with enough invariant effective Cartier divisors.
Tohoku Math. J. 50 (1998), 139--157.
%
 \bibitem{Ka} 
 M.M.~Kapranov: 
 Chow quotients of Grassmannians. I.
 Advances in Soviet Math. 16 (1993), Part 2, 29--110.
%
 \bibitem{King}
 A.D.~King: 
 Moduli of representations of finite-dimensional algebras.  
 Quart. J. Math. Oxford Ser. (2) 45 (1994), no.~180, 515--530.
%
 \bibitem{Ki}
 F.~Kirwan:
 Cohomology of quotients in symplectic and algebraic geometry.
 Math. Notes 31, Princeton University Press, Princeton, NJ, 1984
%
\bibitem{Kl}  
P.~Kleinschmidt: 
A classification of toric varieties with few generators. 
Aequationes Math.~35 (1988), no.~2-3, 254--266.
%
%
 \bibitem{Kn}
 F.~Knop: 
 \"Uber Hilberts vierzehntes Problem f\"ur 
 Variet\"aten mit Kompliziertheit eins. 
 Math. Z. 213 (1993), no.~ 1, 33--36.
%
 \bibitem{KKLV}
 F.~Knop, H.~Kraft, D.~Luna, Th.~Vust: 
 Local properties of algebraic group actions.  
 Algebraische Transformationsgruppen und Invariantentheorie,  
 63--75, DMV Sem., 13, Birkh\"auser, Basel (1989)
%
 \bibitem{KKV}
 F.~Knop, H.~Kraft, Th.~Vust: 
 The Picard group of a G-variety.
 Algebraische Transformationsgruppen und Invariantentheorie,  
 77--87, DMV Sem., 13, Birkh\"auser, Basel (1989)
%
 \bibitem{Kon} 
 S.~Kond\=o:
 Algebraic K3-surfaces with finite automorphism groups.
 Nagoya Math. J. 116, (1989), 1--15.
%
 \bibitem{Kov} 
 S.~Kov\'acs:
 The cone of curves of a K3-surface.
 Math. Ann. 300 (1994), 681--691.
%
 \bibitem{KoRu}
 M.~Koras, P.~Russell: 
 Linearization problems.  
 In Algebraic group actions and quotients,  
 91--107, Hindawi Publ. Corp., Cairo (2004)
%
 \bibitem{Kr} 
 H.~Kraft:
 Geometrische Methoden in der Invariantentheorie.
Vieweg Verlag, Braunschweig, 1984
%
 \bibitem{KrPo} 
 H.~Kraft, V.L.~Popov:
 Semisimple group actions on the three dimensional affine space are linear. 
 Comm. Math. Helv. 60 (1985), no.~3, 466--479.
%
 \bibitem{LaVe1}
 A.~Laface, M.~Velasco:
 A survey on Cox rings.
 Geom. Dedicata 139 (2009), 269--287. 
%
 \bibitem{LaVe2}
 A.~Laface, M.~Velasco:
 Picard-graded Betti numbers and the defining ideals of Cox rings.  
 J. Algebra 322 (2009), no.~2, 353--372.
%
 \bibitem{La}
 R.~Lazarsfeld:
 Positivity in algebraic geometry. I.
 Classical setting: line bundles and linear series. A Series of Modern Surveys in
 Mathematics 48. Springer-Verlag, Berlin, 2004
%
 \bibitem{Lu}
 D.~Luna:
 Slices etales. 
 Bull. Soc. Math. Fr., Suppl. M\'em. 33 (1973), 81--105.
%
%
\bibitem{LuVu}
D. Luna, Th. Vust: Plongements d'espaces homog\`{e}nes. 
Comment. Math. Helv.  58  (1983),  no. 2, 186--245.
%
\bibitem{Mat} K. Matsuki: Introduction to the Mori program. 
Berlin, Heidelberg, New York, Springer Verlag (2001)
%
%
 \bibitem{Ma}
 A.L.~Mayer: 
 Families of K3-surfaces.
 Nagoya Math. J. 48 (1972), 1--17.
%
 \bibitem{Mi}
 J.S.~Milne: 
 \'Etale Cohomology. Princeton University Press, Princeton, NJ, 1980
%
 \bibitem{MiSt} 
 E.~Miller, B.~Sturmfels: 
 Combinatorial Commutative Algebra.
 Graduate Texts in Math. 227, Springer-Verlag, New York, 2005
%
\bibitem{Mo}
 S.~Mori:
 Graded factorial domains.
 Japan J. Math. 3 (1977), no.~2, 223--238.
%
 \bibitem{Mor}
 D.~Morrison:
 On K3-surfaces with large Picard number.
 Invent. Math. 75, (1984), no.~1, 105--121.
%
\bibitem{MuRed}
D.~Mumford:
The red book of varieties and schemes.
Springer-Verlag
%
 \bibitem{MuFoKi}
 D.~Mumford, J.~Fogarty, F.~Kirwan:
 Geometric invariant theory. Third edition. 
 Ergebnisse der Mathematik und ihrer Grenzgebiete 34. 
 Springer-Verlag, Berlin, 1994
%
 \bibitem{Mus}
 I.M.~Musson:
 Differential operators on toric varieties. 
 J. Pure Appl. Algebra 95 (1994), no.~3, 303--315.
%
 \bibitem{Mu} 
 M.~Mustata: 
 Vanishing theorems on toric varieties. 
 Tohoku Math. J. (2) 54 (2002), no.~3, 451--470.
%
 \bibitem{Ni1}
 V.V.~Nikulin:
 A remark on algebraic surfaces with polyhedral Mori cone.
 Nagoya Math. J. 157 (2000), 73--92,
%
 \bibitem{Ni2} 
 V.V.~Nikulin:
 Factor groups of groups of automorphisms
 of hyperbolic forms with respect to subgroups
 generated by 2-reflections.
 J. Soviet Math. 22 (1983), 1401--1475.
%
 \bibitem{Ni3} 
 V.V.~Nikulin:
 Surfaces of type K3 with a finite automorphism 
 group and a Picard group of rank three.
 Trudy Mat. Inst. Steklov 165 (1984), 119--142 (Russian); English transl.:
 Proc. Steklov Math. Institute, Issue 3, 131--155 (1985)
%
\bibitem{Oda}
T.~Oda:
Convex bodies and algebraic geometry. 
An introduction to the theory of toric varieties. 
Ergebnisse der Mathematik und ihrer Grenzgebiete (3),
Springer-Verlag, Berlin, 1988
%
\bibitem{Oda-Pr}
T.~Oda:
Problems on Minkowski sums of convex polytopes.
The Oberwolfach Conference "Combinatorial Convexity 
and Algebraic Geometry", 26.10--02.11, 1997, 
arXiv:0812.1418
%
 \bibitem{OdPa}
 T.~Oda, H.S.~Park: 
 Linear Gale transforms and Gelfand-Kapranov-Zelevinskij decompositions.  
 Tohoku Math. J. (2) 43 (1991), no.~3, 375--399.
%
\bibitem{OnVi} A.L.~Onishchik, E.B.~Vinberg: Lie Groups and Algebraic Groups.
Springer-Verlag, Berlin Heidelberg, 1990
%
%
 \bibitem{Pa} 
 T.E.~Panov: Topology of Kempf-Ness sets for algebraic torus actions.
Proc. Steklov Math. Inst. 263 (2008), 1--13.
%
 \bibitem{Park}
 H.S.~Park: 
 The Chow rings and GKZ-decompositions for $\QQ$-factorial toric varieties. 
 Tohoku Math. J. (2) 45 (1993), no.~1, 109--145.
%
 \bibitem{Pe} 
 M.~Perling:
 Graded rings and equivariant sheaves on toric varieties.
 Math. Nachr. 265 (2004), 87--197.
%
 \bibitem{PeSu} 
 L.~Petersen, H.~S\"u{\ss}:
 Torus invariant divisors.
 Preprint, arXiv:0811.0517, to appear in Israel J. Math.
%
 \bibitem{PiSh}
 I.~Piatetski-Shapiro, I.R.~Shafarevich:
 A Torelli theorem  for algebraic surfaces of type K3.
 Math. USSR Izv. 5 (1971), 547--587.
%
\bibitem{Po}
 O.N.~Popov:
 The Cox ring of a Del Pezzo surface has rational singularities.
 Preprint, arXiv:math.AG/0402154.
%
 \bibitem{PoVi}
 V.L.~Popov, E.~Vinberg: 
 Invariant Theory.
 Encyclopedia Math. Sciences 55, 123--185,
 Springer-Verlag, Heidelberg, 1994
%
 \bibitem{Ren} 
 L.~Renner: 
 The cone of semi-simple monoids with the same factorial hull. 
 Preprint, arXiv:math.AG/0603222.
%
 \bibitem{Res}
 N.~Ressayre: 
 The GIT-equivalence for G-line bundles.  
 Geom. Dedicata 81 (2000), no.~1-3, 295--324.
%
 \bibitem{Ru} 
 P.~Russell: 
 Gradings of polynomial rings.  
 In Algebraic geometry and its applications (West Lafayette, IN, 1990),  
 365--373, Springer, New York (1994)
%
 \bibitem{Sai} 
 B.~Saint-Donat:
 Projective models of K3-surfaces.
 Amer. J. Math. 96 (1974), no.~4, 602--639.
%
 \bibitem{Sa} 
 P.~Samuel: 
 Lectures on unique factorization domains. 
 Tata Institute, Mumbay, 1964
%
 \bibitem{ScSt} 
 G.~Scheja, U.~Storch:
 Zur Konstruktion faktorieller graduierter Integrit\"atsbereiche.
 Arch. Math. (Basel) 42 (1984), no.~1, 45--52.
%
 \bibitem{SeSk} 
 V.V.~Serganova, A.N~Skorobogatov: 
 Del Pezzo surfaces and representation theory. 
 Algebra Number Theory  1  (2007),  no.~4, 393--419.
%
\bibitem{SeSk2}
[122] V.V.~Serganova, A.N.~Skorobogatov: 
On the equations for universal torsors 
over del Pezzo surfaces. 
J. Inst. Math. Jussieu 9 (2010), no.~1, 203-223.
%
 \bibitem{Se} 
 C.S.~Seshadri: 
 Quotient spaces modulo reductive algebraic groups. 
 Ann. Math. (2) 95 (1972), 511--556.
%
\bibitem{Sp}
T.A.~Springer:
Linear algebraic groups. 
Reprint of the 1998 second edition. 
Modern Birkh\"auser Classics. 
Birkh\"auser Boston, Inc., Boston, MA, 2009
%
\bibitem{Sk0}  
A.N.~Skorobogatov: 
On a theorem of Enriques-Swinnerton-Dyer.
Ann. Fac. Sci. Toulouse Math. (6) 2 (1993), 
no.~3, 429--440.
%
\bibitem{Sk} 
A.N.~Skorobogatov: 
Torsors and rational points. 
Cambridge Tracts in Math.~144, Cambridge University Press, 2001
%
\bibitem{Stein}
R.~Steinberg: Nagata's example.
Algebraic groups and Lie groups, 375--384, Austral. Math. Soc. Lect. Ser. 9, 
Cambridge Univ. Press, Cambridge, 1997
%
 \bibitem{StTeVe}
 M.~Stillman, D.~Testa, M.~Velasco: 
 Gr\"obner bases, monomial group actions, 
 and the Cox rings of del Pezzo surfaces.  
 J. Algebra 316 (2007), no.~2, 777--801.
%
 \bibitem{St}
 B.~Sturmfels: 
 Gr\"obner bases and convex polytopes. 
 University Lecture Series, vol.~8, American Mathematical Society, Providence, RI, 1996
%
\bibitem{StVe}
 B.~Sturmfels, M.~Velasco:
 Blow-ups of $\mathbb{P}^{n-3}$ at n points and spinor varieties.
 Preprint, arXiv:0906.5096.
%
 \bibitem{StXu}
 B.~Sturmfels, Zh.~Xu:
 SAGBI bases of Cox-Nagata rings.
 Preprint, arXiv:0803.0892.
 %
\bibitem{Sum}
 H.~Sumihiro:  
 Equivariant completion.  
 J. Math. Kyoto Univ. 14  (1974), 1--28.
%
 \bibitem{Su} 
 H.~S\"u{\ss}:
 Canonical divisors on T-varieties.
 Preprint, arXiv:0811.0626.
%
 \bibitem{Sw1}
 J.~\'Swi\c{e}cicka: 
 Quotients of toric varieties by actions of subtori.  
 Colloq. Math. 82 (1999), no.~1, 105--116.
%
\bibitem{Sw2}
J.~\'Swi\c{e}cicka: 
A combinatorial construction of sets with good quotients 
by an action of a reductive group.  
Colloq. Math. 87 (2001), no.~1, 85--102. 
%
 \bibitem{TeVaVe}
 D.~Testa, A.~ V\'arilly-Alvarado, M.~Velasco:
 Cox rings of degree one del Pezzo surfaces. 
 Preprint, arXiv:0803.0353.
%
 \bibitem{TeVaVe2}
 D.~Testa, A.~Varilly, M.~Velasco:
 Big rational surfaces.
 Preprint, arXiv:0901.1094.
%
 \bibitem{Th1}
 M.~Thaddeus: 
 Geometric invariant theory and flips.  
 J. Amer. Math. Soc. 9 (1996), no.~3, 691--723.
%
%
%
\bibitem{Tim}
D.A.~Timashev:
Homogeneous spaces and equivariant embeddings.
To appear in Encyclopaedia Math. Sci. {\bf 138}, Springer-Verlag, 2010; 
arXiv:math/0602228.
%
\bibitem{Vi-c}
E.B.~Vinberg:
Complexity of actions of reductive groups.
Func. Anal. Appl. 20 (1986), no.~1, 1--11.
%
 \bibitem{Vi} 
 E.B.~Vinberg: 
 On reductive algebraic semigroups.
 In: Lie Groups and Lie Algebras, E.B.~Dynkin's Seminar, Amer. Math. Soc. Transl.
 Ser.~2, vol.~169, Amer. Math. Soc., Providence, RI, 145--182 (1995)
%
 \bibitem{Wl}
 J.~W\l odarczyk:
 Embeddings in toric varieties and prevarieties.
 J. Alg. Geom. 2 (1993), no.~4, 705--726.
%
 \bibitem{Za}
 D.-Q.~Zhang: 
 Quotients of K3-surfaces modulo involutions.
 Japan. J. Math. (N.S.)  24  (1998), no.~2, 335--366.
%
\end{thebibliography}
\end{document}